\documentclass[10pt, a4paper, twoside, titlepage]{article}
\usepackage[a4paper, margin=1.25in, marginparwidth=1in, marginparsep=0.125in]{geometry}
\pdfoutput=1

\usepackage[ngerman, english]{babel} 
\usepackage{amsmath}
\usepackage{amssymb}
\usepackage{amsthm}
\usepackage{thmtools}
\usepackage{algorithm, algpseudocode}
\usepackage{bbm}
\usepackage{blkarray}
\usepackage{pgfplots}
\usetikzlibrary{graphs,arrows.meta,matrix,fit,patterns,shapes}
\pgfplotsset{compat=1.8}
\usepackage{makecell}
\usepackage{svg}
\usepackage{colonequals}
\usepackage{soul}
\usepackage{color}
\usepackage{setspace}
\usepackage{mathtools}
\usepackage{hyperref}
\usepackage[nameinlink,capitalize]{cleveref}
\usepackage[toc,acronym,symbols]{glossaries}
\usepackage{stmaryrd}
\usepackage{eso-pic}
\usepackage{transparent}
\usepackage{graphicx}
\usepackage[framemethod=TikZ]{mdframed}
\usepackage[shortlabels]{enumitem}
\usepackage[nottoc,numbib]{tocbibind}
\usepackage{ifoddpage}
\usepackage{lmodern} 
\usepackage[T1]{fontenc}
\usepackage{mathpazo} 
\usepackage{marginnote}

\usepackage{titlesec} 
\numberwithin{equation}{section} 

\usepackage{tabularx}
\usepackage{multicol}
\newcolumntype{Y}{>{\centering\arraybackslash}X}
\usepackage[column=C]{cellspace}
\addparagraphcolumntypes{X}
\definecolor{tumBlue}{RGB}{0,101,189} 
\definecolor{tumDarkBlue}{RGB}{0,82,147} 
\definecolor{tumLightBlue}{RGB}{100,160,200} 
\definecolor{tumLighterBlue}{RGB}{152,198,234} 
\definecolor{tumOrange}{RGB}{227,114,34} 
\definecolor{tumGreen}{RGB}{162,173,0} 
\definecolor{tumGray}{RGB}{153,153,153} 
\definecolor{tumLight}{RGB}{218,215,203} 

\definecolor{commentpurple}{RGB}{102, 0, 102} 
\colorlet{commentpurplemuted}[RGB]{commentpurple!20!white}

\colorlet{sectionblue}{tumBlue}
\definecolor{linkred}{RGB}{127,0,0} 
\definecolor{darklinkred}{RGB}{50,0,0} 
\colorlet{headingcolor}{sectionblue}
\colorlet{headingcolormuted}[RGB]{headingcolor!20!white}
\colorlet{linkcolor}{linkred}

\usepackage[labelfont={color=headingcolor,bf}]{caption}

\makeatletter
\def\moverlay{\mathpalette\mov@rlay}
\def\mov@rlay#1#2{\leavevmode\vtop{%
   \baselineskip\z@skip \lineskiplimit-\maxdimen
   \ialign{\hfil$\m@th#1##$\hfil\cr#2\crcr}}}
\newcommand{\charfusion}[3][\mathord]{
    #1{\ifx#1\mathop\vphantom{#2}\fi
        \mathpalette\mov@rlay{#2\cr#3}
      }
    \ifx#1\mathop\expandafter\displaylimits\fi}
\makeatother

\newcommand{\bigcupdot}{\charfusion[\mathop]{\bigcup}{\cdot}}

\newcommand{\midd}[1]{\mathrel{}\middle#1\mathrel{}}

\mdfsetup{ skipbelow=-0.2em }

\newlength{\spaceblength}
\declaretheoremstyle[
    headfont=\bfseries,
    notefont=\bfseries,
    notebraces={}{\\[\parskip]}, 
    bodyfont=\normalfont\upshape,
    headpunct={},
    postheadspace=\spaceblength,
    spacebelow=\parskip,
    spaceabove=\parskip,
    headformat=\color{headingcolor}\NAME\ \NUMBER,
    qed={$\begingroup\color{headingcolormuted}\blacktriangleleft\endgroup$}
]{boxstyle}
\declaretheoremstyle[
    headfont=\bfseries\itshape,
    notefont=\normalfont\bfseries,
    notebraces={}{\\[\parskip]}, 
    bodyfont=\normalfont,
    headpunct={},
    postheadspace=\spaceblength,
    spacebelow=\parskip,
    spaceabove=\parskip,
    headformat=\color{headingcolor}\NAME,
    qed=\qedsymbol
]{proofstyle}
\declaretheoremstyle[
    headfont=\bfseries\itshape,
    notefont=\bfseries\itshape\hypersetup{hidelinks},
    notebraces={}{}, 
    bodyfont=\normalfont,
    headpunct={},
    postheadspace=\spaceblength,
    spacebelow=\parskip,
    spaceabove=\parskip,
    headformat={\color{headingcolor}\NAME\ \NOTE},
    qed=\qedsymbol
]{oproofstyle}
\declaretheoremstyle[
    headfont=\bfseries,
    notefont=\bfseries,
    notebraces={}{\\[\parskip]}, 
    bodyfont=\normalfont,
    headpunct={},
    postheadspace=\spaceblength,
    spacebelow=\parskip,
    spaceabove=\parskip,
    headformat=\color{headingcolor}\NAME\ \NUMBER,
    qed={$\begingroup\color{headingcolormuted}\blacktriangleleft\endgroup$}
]{examplestyle}
\declaretheoremstyle[
    headfont=\normalfont,
    notefont=\bfseries,
    notebraces={}{\\[\parskip]}, 
    bodyfont=\normalfont,
    headpunct={},
    postheadspace=\spaceblength,
    spacebelow=\parskip,
    spaceabove=\parskip,
    headformat=\color{headingcolor}\NAME,
    qed={$\begingroup\color{headingcolormuted}\blacktriangleleft\endgroup$}
]{remarkstyle}

\declaretheorem[style=boxstyle]{definition}
\declaretheorem[style=boxstyle,sibling=definition]{theorem}

\declaretheorem[style=boxstyle,sibling=definition]{lemma}
\declaretheorem[style=boxstyle,sibling=definition]{corollary}

\declaretheorem[style=boxstyle,sibling=definition]{conjecture}

\declaretheorem[style=proofstyle,numbered=no,name=Proof]{tproof}

\declaretheorem[style=oproofstyle,numbered=no,name=Proof of]{oproof}

\addto\extrasenglish{

}

\newlength{\secskip}
\titleformat*{\section}{\large\bfseries\color{headingcolor}}
\titleformat*{\subsection}{\bfseries\color{headingcolor}}
\titleformat*{\subsubsection}{\color{headingcolor}}

\hypersetup{
    colorlinks=true,
    allcolors=.,
    urlcolor=linkcolor,
    citecolor=linkcolor,
    linkcolor=linkcolor,
}

\newtagform{roman}[]()

\newcommand{\Prb}[1] {\mathbb{P}\left[#1\right]}

\newcommand{\Ex}[1] {\mathbb{E}\left[#1\right]}

\newcommand{\dx}[1] {\;\mathrm{d}#1}

\newcommand{\bbP} {\mathbb{P}}
\newcommand{\bbE} {\mathbb{E}}

\newcommand{\bbR} {\mathbb{R}}

\newcommand{\bbN} {\mathbb{N}}
\newcommand{\bbZ} {\mathbb{Z}}

\newcommand{\calF} {\mathcal{F}}
\newcommand{\calG} {\mathcal{G}}

\newcommand{\calN} {\mathcal{N}}

\newcommand{\bP} {\mathbf{P}}
\newcommand{\bQ} {\mathbf{Q}}

\newcommand{\ivar} {\mathrm{\mathbf{i}}}
\newcommand{\1} {\mathbbm{1}}

\newcommand{\stickman}[2] {\begin{scope}[shift={#1}]
  \draw [#2] (-0.2, -0.5) -- (0,-0.25) -- (0,0.25);
  \draw [#2] (0.2, -0.5) -- (0,-0.25);
  \draw [#2] (-0.2, 0.15) -- (0,-0.1);
  \draw [#2] (0.2, 0.15) -- (0,-0.1);
  \fill [#2] (0,0.25) circle (0.125);
\end{scope}}

\tikzset{
  diagonal fill a/.style n args=2{path picture={%
  \draw[fill=#1, draw=none] (path picture bounding box.south west) --
              (path picture bounding box.north east) -- (path picture bounding box.south east) -- cycle; %
  \draw[fill=#2, draw=none] (path picture bounding box.south west) --
              (path picture bounding box.north east) -- (path picture bounding box.north west) -- cycle;}},
  diagonal fill b/.style n args=2{path picture={%
  \draw[fill=#1, draw=none] (path picture bounding box.south west) --
              (path picture bounding box.north east) -- (path picture bounding box.south east) -- cycle; %
  \draw[pattern=north west lines, pattern color=#2, draw=none] (path picture bounding box.south west) --
              (path picture bounding box.north east) -- (path picture bounding box.north west) -- cycle;}},
  diagonal fill c/.style n args=2{path picture={%
  \draw[pattern=north west lines, pattern color=#1, draw=none] (path picture bounding box.south west) --
              (path picture bounding box.north east) -- (path picture bounding box.south east) -- cycle; %
  \draw[fill=#2, draw=none] (path picture bounding box.south west) --
              (path picture bounding box.north east) -- (path picture bounding box.north west) -- cycle;}},
  diagonal fill d/.style n args=2{path picture={%
  \draw[pattern=north west lines, pattern color=#1, draw=none] (path picture bounding box.south west) --
              (path picture bounding box.north east) -- (path picture bounding box.south east) -- cycle; %
  \draw[pattern=north west lines, pattern color=#2, draw=none] (path picture bounding box.south west) --
              (path picture bounding box.north east) -- (path picture bounding box.north west) -- cycle;}},
  table nodes/.style={
    rectangle,
    draw=none,
    align=center,
    minimum height=7mm,
    text depth=0.5ex,
    text height=2ex,
    inner xsep=0pt,
    outer sep=0pt
  },      
  table/.style={
    matrix of nodes,
    row sep=-\pgflinewidth,
    column sep=-\pgflinewidth,
    nodes={
        table nodes
    }
  }
}

\usepackage{ifthen}
\newboolean{extendedversion}
\setboolean{extendedversion}{false}

\title{Interacting Edge-Reinforced Random Walks}

\author{Nina Gantert, Fabian Michel, Guilherme Reis}
\date{17/10/2022}

\begin{document}

\begin{center}
  \textcolor{sectionblue}{
    \LARGE\bfseries
    Interacting Edge-Reinforced Random Walks
  }\\[1em]
  Nina Gantert, Fabian Michel, Guilherme Reis
\end{center}

\section*{Abstract}

We consider the edge-reinforced random walk with multiple (but finitely many) walkers which influence the
edge weights together. The walker which moves at a given time step is chosen uniformly at random, or according to a fixed order. First, we
consider 2 walkers with linear reinforcement on a line graph comprising three nodes. We show that the edge weights
evolve similarly to the setting with a single walker which corresponds to a P\'olya urn. In particular, the left edge weight proportion is
a martingale at certain stopping times, showing that a (random) limiting proportion exists. We then look at an arbitrary number
of walkers on $\bbZ$ with very general reinforcement. We show that in this case, the behaviour is also the
same as for a single walker: either all walkers are recurrent or all walkers have finite range. In the particular case of reinforcements of ``sequence type'', we give a criterion for recurrence.

\ifthenelse{\boolean{extendedversion}}
{{
  \hypersetup{hidelinks}
  \tableofcontents
}}{}

\section{Introduction}

\subsection{Model}

In this paper, we study the edge-reinforced random walk with multiple walkers. The edge-reinforced random walk (ERRW) with a single walker is a stochastic process in discrete time defined on a graph. The edges in the graph are weighted, and the probability to leave a
node via one of the incident edges is proportional to the respective edge weight compared to
the weights of the other incident edges. Each time an edge is crossed, its weight is increased
according to some reinforcement scheme. Thus, it becomes more likely to visit parts of the
graph which have already been visited before and the process is not a Markov process. If the initial edge weights are the same everywhere and the edge weight increment
is $1$ (upon traversal), then the walk is called linearly edge-reinforced random walk (LERRW), the most commonly studied model of reinforced random walks.

We extend the above model to include multiple walkers. There is a common set of edge weights for all walkers,
and at every time step, one of the walkers moves (the walker to move may be selected at random, or according to some fixed order).
The probability to take a specific edge is still proportional to the weight, and we keep the reinforcement: the weight
of an edge is increased whenever any of the walkers crosses it.

\subsection{Literature}

The ERRW (as well as its counterpart, the vertex-reinforced random walk) has been studied extensively,
the first papers dating back to 1987, when the model was introduced by Coppersmith and Diaconis. Even before,
\cite{definettimarkovchain} showed that the LERRW (if certain assumption are satisfied) has a representation
as a mixture of Markov chains, i.e.~a Markov chain with random transition probabilities, also known as random walk in random environment.
Much later, \cite{rwforerrw} showed that this representation can be used for the
LERRW on any graph, and \cite{magicformula} even gave a formula for the so-called mixing measure on finite
graphs. The mixing measure simply is the distribution of the random transition probabilities in the mixture of
Markov chains.

Relatively early, results for the LERRW on trees were obtained. \cite{errwpemantle,rwrelyonspemantle,rwrelyonspemantlecorr}
showed that there is a phase transition between recurrence and transience in the initial edge weights. The interesting
case of $\bbZ^d$ remained open much longer until \cite{localizationerrw,errwandvrjp,transienceerrw} showed that for $d \geq 3$,
there is again a transition from recurrence to transience in the initial weights. \cite{monphasetrans} proved that this transition
is sharp, i.e.~there is a certain critical initial edge weight such that for smaller initial weights, the random walk is
recurrent, and transient for larger weights.

An overview of results on reinforced processes in general can be found in \cite{surveyreinf,rrwsurveykozma}. These surveys
also show that reinforced random walks are closely related to urn processes, which have a very similar reinforcement
component to the linearly reinforced walks: in most urn models, when a ball of a certain color is drawn, a fixed number of
balls of the same color is added to the urn. This is also a type of linear reinforcement, and urns have been used on many
occasions to analyze reinforced walks (see, for example, \cite{errwpemantle}).

Results for multiple interacting reinforced random walks are scarce. In~\cite{ErhardReis}, the authors consider multiple ant random walks.
However, the ant random walk is defined via a reinforcement rule on oriented edges instead of on non-oriented edges (as for the
ERRW). \cite{ErhardReis} also considers the case of superlinear reinforcement and compares the reinforced walk
with P\'olya urns. 

Edge-Reinforced walks with multiple walkers on non-oriented graphs have (to the best of our knowledge) not been studied yet.
For multiple walkers, there is no representation as a mixture of Markov chains, even in the case with
linear reinforcement. Therefore, most of the methods which have been used to study the LERRW cannot be applied
to the case of multiple walkers. However, our paper is partially based on methods which have been used to study
urn processes as there still is a similarity to urn models where multiple balls are drawn from
the urn. In particular, \cite{hypergrru} studied urns where multiple balls are drawn and replaced, and new balls
are added according to a random distribution depending on the balls which were drawn. The methods used to study the
limiting proportions of balls in the urn are also useful to study two edge-reinforced walkers on a 3-node segment (\begin{tikzpicture}
  \node[circle,fill=black,inner sep=0.7mm] (N0) at (0, 0) {};
  \node[circle,fill=black,inner sep=0.7mm] (N1) at (0.7, 0) {};
  \node[circle,fill=black,inner sep=0.7mm] (Nm1) at (-0.7, 0) {};
  \draw (Nm1) -- (N0) -- (N1);
\end{tikzpicture}) and
the limiting proportions of the corresponding edge weights.
In addition, \cite{errwz} studied the ERRW on $\bbZ$ with very general reinforcement (and not just linear reinforcement).
\cite{errwz} showed that a single walker is either recurrent or visits only a finite part of $\bbZ$.
The methods used there can be adapted for multiple walkers, but the method is specific to $\bbZ$.

\subsection{Results}

In the following, the random walkers will be denoted by sequences $\big(X^{(i)}_n\big)_{n \geq 0}$, where $X^{(i)}_n$ is equal to the position of walker $i$ at time $n$, so $X^{(i)}_n$ is a node of the graph.
We will show the following four main results:
\begin{itemize}
  \item \autoref{thm:main-alternating}: consider two walkers on a 3-node segment with linear reinforcement where the walkers move alternately. Then, the proportion of the left edge weight converges to a random limit which has no atoms.
  Specifically, if we denote the weights of the left and right edges at time $n$ by $w(n,-1)$ and $w(n,0)$, respectively, then the fraction $\frac{w\left(n,-1\right)}{w\left(n,-1\right) + w\left(n,0\right)}$ converges to a random limit $M_\infty$
  which satisfies $\Prb{M_\infty = a} = 0$ for any $a \in \left[0, 1\right]$.
  \item \autoref{thm:randwalkersmain}: consider two walkers on a 3-node segment with linear reinforcement where the walker to move is selected uniformly at random at every step. Then, the proportion of the left edge weight converges to a random limit. Furthermore, the left edge weight proportion
  is a martingale if looked at at certain stopping times. Specifically, the fraction $\frac{w\left(n,-1\right)}{w\left(n,-1\right) + w\left(n,0\right)}$ converges to a random limit, which is identical to the limit of the martingale
  $\frac{w\left(\tau_n,-1\right)}{w\left(\tau_n,-1\right) + w\left(\tau_n,0\right)}$, where $\tau_n = \inf\left\{ m > \tau_{n-1} : X^{\left(1\right)}_m \textrm{ and } X^{\left(2\right)}_m \textrm{ are both in the center node}\right\}$.
  \item \autoref{thm:allRecurrentOrAllFiniteRange}: consider $K$ walkers on $\bbZ$ with very general reinforcement, but where all but finitely many initial weights are $1$, and where the walker to move is selected uniformly at random. Then, either all walkers are recurrent or all walkers have finite range (visit only a finite part of $\bbZ$) a.s. That is, under these assumptions,
  \begin{align*}
    \Prb{\forall i: X^{\left(i\right)} \textrm{ is recurrent}}
    + \Prb{\forall i: X^{\left(i\right)} \textrm{ has finite range}} &= 1\, . \qedhere
  \end{align*}
  \item \autoref{thm:seqtypefinrange}: consider $K$ walkers on $\bbZ$ with sequence type reinforcement where the walker to move is selected uniformly at random. Then, we can characterize the two possible behaviors depending on the sequence of edge weight increments.
  ``Sequence type'' means that
  all initial edge weights are $1$, and that the edge weights are then increased by a fixed sequence of increments $\mathbf{a} = \left(a_k\right)_{k \geq 1}$ which do not
  depend on the location of the edge. Specifically, if we set
  \begin{align*}
    \phi\left(\mathbf{a}\right) = \sum_{k = 1}^{\infty} \left(1 + \sum_{l = 1}^k a_l\right)^{-1} \,,
  \end{align*}
  then, if $\phi\left(\mathbf{a}\right) = \infty$, all walkers are recurrent a.s., and if $\phi\left(\mathbf{a}\right) < \infty$,
  all walkers have finite range a.s.
\end{itemize}

\subsection{Methods}

\autoref{ssec:multiple_walkers_segment} treats the case of two walkers on the 3-node segment. When considering alternating walkers, the analysis of the left edge weight
proportion $\frac{w\left(n,-1\right)}{w\left(n,-1\right) + w\left(n,0\right)}$ mostly relies on relatively straightforward calculations which show that the proportion is a martingale at times $4n$.
Since this martingale is bounded, the existence of a limit can be derived. 
We will use the notation $M_n = \frac{w\left(4n,-1\right)}{w\left(4n,-1\right) + w\left(4n,0\right)}$, und use $M_\infty$ to denote the (random) limit.
Subsequently, methods taken from \cite{hypergrru} allow us to show the following central limit type result in \autoref{thm:main-alternating} \ref{thm:altwalkersclt}:
\begin{align*}
  \sqrt{n} \left(M_n - M_{\infty}\right) \to \calN\left(0, \frac{1}{2}M_\infty\left(1 - M_\infty\right)\right)
\end{align*}
where $\calN\left(\mu, \sigma^2\right)$ denotes a normal distribution. This result, and the comparison to the case with a single walker
(which is in 1-to-1 correspondence to a P\'olya urn where two balls of the drawn color are added after every draw) allow us to show that
$M_\infty$ has no atoms. The question of whether $M_\infty$ has a density remains open.

For the 3-node segment with two walkers where the walker to move is selected randomly, the straightforward calculations from the previous case
become much more complicated. Using a recursive formula for the conditional expectation, we can show that
$\frac{w\left(\tau_n,-1\right)}{w\left(\tau_n,-1\right) + w\left(\tau_n,0\right)}$ is a martingale, where $\tau_n = \inf\left\{ m > \tau_{n-1} : X^{\left(1\right)}_m \textrm{ and } X^{\left(2\right)}_m \textrm{ are both in the center node}\right\}$.
This again results in the existence of a limit random variable.  The method we used for the
alternating walkers to show that no atoms exist might be applicable in this case as well, but we could not carry through some of the necessary calculations.

\autoref{ssec:multiple_walkers_z} treats a finite number of walkers $K$ on $\bbZ$. To show that all walkers are recurrent or all have finite range, we use a martingale technique which originates from \cite{errwz}.
In particular, if you sum up the inverse of the edge weights between $0$ and the position of one of the walkers, the resulting sum is a nonnegative supermartingale.
By convergence of this supermartingale, we can deduce that every walker either reaches $0$ at some point or only visits finitely many nodes which
have not been visited before by any of the other walkers. In addition, you can exchange tails of the paths of walkers after they met in a certain node without affecting the law
of the overall process. These two key lemmata allow us to show that all walkers will have the same behavior.

Finally, in order to characterize the two possible behaviors for sequence type reinforcements, we only have to adapt the method from \cite{errwz} slightly in
order to apply it to multiple walkers simultaneously. The proof idea here is to show that all walkers get stuck on a single edge with positive probability if the
weight increments grow fast enough. On the other hand, one can show that the neighbor of a node which is visited infinitely often is visited infinitely often as well if
the increments grow slower.

\section{Preliminaries}
\label{sec:prelim}

\subsection{Graphs}

Throughout this paper, random walks on $\bbZ$ and on the 3-node segment will be considered.
Recall that $\bbZ$ has node set $V = \bbZ$ and edge set
$E = \left\{ \left\{u, v\right\} \midd| u,v \in V \textrm{ and } \left|u - v\right| = 1 \right\}$. The 3-node segment is the subgraph $P_3$ of $\bbZ$ with vertex set $\{-1,0,1\}$ and edges $\{\{-1,0\},\{0,+1\}\}$.

\subsection{Edge-reinforced random walk}
We define the edge-reinforced random walk, which is the main object of this paper. 

\begin{definition}[Edge-Reinforced Random Walk]
  \label{def:errw}
  For this definition we consider only the graph $\bbZ$. We define the crossing number induced by
  the \textbf{edge-reinforced random walk} (ERRW) $\left(X_n\right)_{n \geq 0}$ on the edge $\{x,x+1\}$
  (which we will identify with $x$ for more concise notation)
  as follows: $c(0, x) = 0$, and
  \begin{align*}
    \textrm{for } n \geq 1: \textrm{ }
    c\left(n, x\right) = \sum_{i=1}^n \mathbbm{1}_{\{\{X_{i-1},\,X_i\} =\{ x,x+1\}\} }\,, \quad w(n, x) := W_x(c(n, x))\,,
  \end{align*}
  i.e.~$c$ counts the number of edge traversals and $w$ describes the evolution of the edge weights. For $x \in \bbZ$, $W_x : \bbN \cup \{0\} \to \bbR_{>0}$
  is a non-decreasing function, called the reinforcement function for the edge $\{x, x+1\}$.
  The choice of reinforcement $w(n,x) = 1 + c(n,x)$ is called
  \textbf{linearly edge-reinforced random walk} (LERRW).

  Given the history $\mathcal H_n=\sigma(X_m\,:\,m\leq n)$, we have
  \begin{align*}
    \Prb{X_{n+1} = X_n+1 \midd| \mathcal{H}_n} &= \dfrac{w(n, X_n)}{w\left(n, X_n-1\right)+ w\left(n, X_n\right) }\\
    \Prb{X_{n+1} = X_n-1 \midd| \mathcal{H}_n} &= \dfrac{w(n, X_n-1)}{w\left(n, X_n-1\right)+ w\left(n, X_n\right) }\,.
  \end{align*}
  In words, the probability to take the transition from $X_n$ to $X_n+1$ at time $n$ is proportional to the weight $w(n,X_n)$ associated
  to the edge $\left\{X_n, X_n+1\right\}$, which depends on the number of traversals of the edge. 
\end{definition}

In the following, we will generalize the edge-reinforced random walk to include multiple walkers which move
on the same graph and influence the edge weights together. We start by looking at a 3-node segment with
two walkers in \autoref{sec:2RW3nodes}, and then consider $\bbZ$ with an arbitrary but finite number of walkers in \autoref{ssec:multiple_walkers_z}.

\subsection{Almost sure conditional convergence}\label{sec:as-cond-conv}

\begin{definition}\label{def:conditionalas}
  Let $\left(\Omega, \calG, \bbP\right)$ be a probability space.
  Consider a real random variable $X$ on $\Omega$ and a sub-$\sigma$-field $\calF$ of $\calG$. We define the \textbf{conditional distribution} $\mu_{X|\calF}$ of $X$ given $\calF$ as follows:
  $\mu_{X|\calF}$ is a mapping from $\Omega$ to the set of probability measures on $\bbR$ (with the Borel $\sigma$-algebra). Let $A \subseteq \bbR$ be a measurable set. Then
  \begin{align*}
    \mu_{X|\calF}\left( A \right) = \Ex{\mathbbm{1}_{X \in A} \midd| \calF}
  \end{align*}
  Note that $\omega \mapsto \mu_{X|\calF}\left( A \right)\left(\omega\right)$ is $\calF$-measurable for every measurable set $A$.

  Consider now a sequence $X_n$ of random variables on $\Omega$ and a sequence of sub-$\sigma$-fields $\calF_n$. Denote by $\mu_n$ the conditional distribution of $X_n$ given $\calF_n$. Let $\mu$ be a random probability measure on $\bbR$,
  where $\mu$ maps from $\Omega$ to the set of probability measures on $\bbR$, and $\omega \mapsto \mu\left(A\right)\left(\omega\right)$ is $\calG$-measurable for every set $A \subseteq \bbR$ which is Borel-measurable.
  We say that $X_n$ converges to $\mu$ in the sense of \textbf{almost sure conditional convergence} with respect to the conditioning system $\calF_n$ if
  $\mu_n$ converges weakly to $\mu$ for almost every $\omega \in \Omega$. Notice that $X_n\to \mu$ in the sense of almost sure conditional convergence with respect to $(\calF_n)_{n\geq 0}$ if, and only if, for all $f:\mathbb{R}\to \mathbb{R}$ continuous and bounded it holds that
  \begin{equation}\label{eq:def:conditionalweakconv}
      \mathbb{E}[f(X_n)|\mathcal{F}_n](\omega)\to \int f(\bar\omega)\mu_\omega(d\bar\omega) \,, 
  \end{equation}
  where $\mu_\omega$ represents a realization of the random measure $\mu$. 
\end{definition}

The notion of almost sure conditional convergence is stronger than the notion of weak convergence: in the setting of \autoref{def:conditionalas},  using the dominated convergence theorem, \eqref{eq:def:conditionalweakconv} implies that, for all continuous and bounded function $f$,
\begin{equation*}
    \mathbb{E}[f(X_n)]\to \int\int f(\bar\omega)\mu_\omega(d\bar\omega)\bbP(d\omega)\,,
\end{equation*}
that is, the sequence $(X_n)_{n\geq 0}$ converges weakly to the distribution defined through
\begin{equation*}
    \bar \mu(A):=\int \mu_\omega(A)\bbP(d\omega)\,.
\end{equation*}

The following lemma will be useful.
\begin{lemma}\label{lemma:weakconv} Let $\mu$ be a random probability measure and 
\begin{equation*}
    \bar \mu(A):=\int \mu_\omega(A)\bbP(d\omega)\,.
\end{equation*}
Assume that $(X_n)_{n\geq 0}$ converges to $\mu$ in the sense of almost sure conditional convergence. Let $(Y_n)_{n\geq 0}$ be a sequence of random variables such that
\begin{equation*}
    \lim_{n\to\infty}|X_n-Y_n|=0\,,\,\text{ almost surely.}
\end{equation*}
Then the laws of $Y_n$ converge weakly to $\bar \mu$ for $n\to\infty$.
\end{lemma}
\begin{tproof}
    Observe that
    \begin{equation*}
        |\mathbb{E}[f(X_n)]-\mathbb{E}[f(Y_n)]|\leq \mathbb{E}[|f(X_n)-f(Y_n)|]\,,
    \end{equation*}
    which goes to zero using the dominated convergence theorem.
\end{tproof}

In this work $\calN\left(0, V\right)$ denotes the distribution of a normal random variable with mean zero and (random) variance $V$.  We write $\overline\calN\left(0, V\right)$ for the mixture of the laws $\calN\left(0, V\right)$ when we average on the randomness of $V$:
    \begin{equation*}
        \overline\calN\left(0, V\right)(A):=\int \calN\big(0, V)(A)d\mathcal L(V)\,.
    \end{equation*}

\begin{theorem}
  \label{thm:a3}
  Let $M_n$ be a bounded martingale with respect to a filtration $\calF_n$ and let $V$ be a random variable. Let $M_\infty$ be the random variable such that $M_n \to M_\infty$ a.s.~and in $L^1$. Assume that
  \begin{enumerate}[(i)]
    \item \label{thm:a3cond1} $\Ex{\displaystyle\sup_{n \geq 1} \sqrt{n} \left|M_{n-1} - M_n\right|} < \infty$
    \item \label{thm:a3cond2} $n \displaystyle\sum_{j \geq n} \left(M_{j-1} - M_j\right)^2 \to V$ a.s. 
  \end{enumerate}
  Then, the following convergence holds in the sense of almost sure conditional convergence
  with respect to $\calF_n$:
  \begin{equation*}
      \sqrt{n} \left(M_n - M_\infty\right)\to \calN(0,V)\,.\qedhere
  \end{equation*} 
\end{theorem}

\begin{tproof}
  See \cite[Theorem A.1 (A.3 in the arxiv version)]{hypergrru}. The stated version is just a simplified version of it.
\end{tproof}

\begin{lemma}
  \label{lem:a2}
  Let $X_n$ be a sequence of random variables adapted to the filtration $\calF_n$ and let $X$ be a random variable. Assume that
  \begin{enumerate}[(i)]
      \item \label{lem:a2cond1} $\displaystyle\sum_{n \geq 1} n^{-2} \Ex{X_n^2} < \infty$ and
      \item \label{lem:a2cond2} $\Ex{X_n \midd| \calF_{n-1}} \to X$ a.s.
  \end{enumerate}
  Then
  \begin{align*}
    &n \sum_{j \geq n} \frac{X_j}{j^2} \to X \textrm{ a.s.} \qedhere
  \end{align*}
\end{lemma}

\begin{tproof}
  See \cite[Lemma A.2]{hypergrru}.
\end{tproof}

\section{Two walkers on a segment of 3 nodes}
\label{sec:2RW3nodes}

\subsection{The model}
\label{ssec:multiple_walkers_segment}

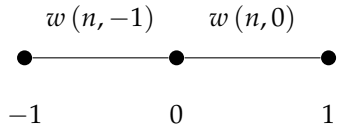
\begin{figure}[H]
  \begin{center}
    \begin{tikzpicture}
      \node[circle,fill=black,inner sep=0.7mm] (N0) at (0, 0) {};
      \node[circle,fill=black,inner sep=0.7mm] (N1) at (2, 0) {};
      \node[circle,fill=black,inner sep=0.7mm] (Nm1) at (-2, 0) {};
      \node[below=5mm] at (N0) {$0$};
      \node[below=5mm] at (Nm1) {$-1$};
      \node[below=5mm] at (N1) {$1$};
      \draw (Nm1) -- node[above=2mm] {$w\left(n,-1\right)$} (N0) -- node[above=2mm] {$w\left(n,0\right)$} (N1);
    \end{tikzpicture}
    \caption{Edge weights of the line segment at time $n$}
    \label{fig:errw_segment_multiple}
  \end{center}
\end{figure}

We start with a very simple model, the linearly edge-reinforced random walk with two walkers on
a segment of $\bbZ$ with three nodes. We use the notation shown in \autoref{fig:errw_segment_multiple}.
We can interpret this as a two-player urn: the left and right edge weights correspond to the number of balls in an urn, and the two
walkers correspond to two players which draw balls from the urn and add more balls of the drawn color.

We adapt \autoref{def:errw} for two walkers on the 3-node segment as follows:
\begin{itemize}
  \item There are two walkers $X^{\left(1\right)}$ and $X^{\left(2\right)}$, which both start at the
  node in the center, i.e.~$X^{\left(1\right)}_0 = 0, X^{\left(2\right)}_0 = 0$.
  \item Initially, both edge weights are $1$. We denote the edge weight of the left edge at time $n$
  by $w\left(n,-1\right)$, the weight of the right edge by $w\left(n,0\right)$. 
  \item Whenever an edge is crossed by either of the walkers, its weight is increased by $1$, which corresponds to
  the linear reinforcement function $W\left(n\right) = 1 + n$ in \autoref{def:errw}. More specifically, we define, for $\ell\in\{-1,0\}$, $c(0, \ell) = 0$ and 
  \begin{equation*}
          \textrm{for } n \geq 1: \textrm{ }
    c\left(n, \ell\right) = \sum_{k=1}^{2}\sum_{i=1}^n \mathbbm{1}_{\{\{X^{(k)}_{i-1},\,X^{(k)}_i\} =\{ \ell,\ell+1\}\} }\,,\,
  \quad \textrm{ and } \quad
      w(n,\ell)=1+c\left(n, \ell\right)\,.
  \end{equation*}
  \item When a walker at the node in the center is about to move, he chooses the edge to traverse with
  probability proportional to the respective edge weight. It remains to decide which walker jumps at each time.
  \item We will consider two variants of the order in which the walkers move: \begin{itemize}
    \item \textbf{Alternating walkers}: the walkers move alternately, i.e.~at odd time steps, walker $1$ moves (in particular, walker
    $1$ moves first at step $1$) and at even time steps, walker $2$ moves. Notice that the walkers
    will meet at the node in the center every four steps. As an example, given the history $\mathcal H_{4n}=\sigma(X^{(j)}_m\,:\,m\leq 4n,j=1,2)$, $X_{4n+1}^{(2)}=X_{4n}^{(2)}$ and $X_{4n+1}^{(1)}=X_{4n}^{(1)}+1$ with probability
  \begin{align*}
     \dfrac{w(4n, 0)}{w\left(4n, -1\right)+ w\left(4n, 0\right) }\,.
  \end{align*}
    \item \textbf{Random walker selection}: at every step, we choose uniformly at random (independently of all other steps)
    which of the two walkers moves. As an example, given the history $\mathcal H_{n}=\sigma(X^{(j)}_m\,:\,m\leq n,j=1,2)$, we select the walker $1$ to jump with probability $\frac{1}{2}$, and assuming that $X_n^{(1)}=0$, we have that $X_{n+1}^{(2)}=X_{n}^{(2)}$ and $X_{n+1}^{(1)}=X_{n}^{(1)}+1$ with probability
  \begin{align*}
     \dfrac{w(n, 0)}{w\left(n, -1\right)+ w\left(n, 0\right) }\,.
  \end{align*}
  \end{itemize}
\end{itemize}

\subsection{Results}

We define
\begin{equation*}
    N_n:=\dfrac{w\left(n,-1\right)}{w\left(n,-1\right) + w\left(n,0\right)}
\end{equation*}
and let $\calG_n:= \sigma\left(N_0, \ldots, N_{n}\right)$. Our goal is to prove limit theorems for the sequence $(N_n)_{n\geq 0}$.

To illustrate our methods, consider first the special case of just one walker on the 3-node graph. The walker decides to jump either to the right or to the left and afterwards comes back again to the center. As a consequence, the number of visits of each edge at times $0$, $2$, $4$, ... has the same distribution as the number of balls in a P\'olya urn where we add $2$ more balls of the drawn color at every step. It is a classical result that the proportion of balls in the P\'olya urn is a martingale. 

In the case of two walkers which move alternately, the walkers meet at node $0$ after every $4$ steps. As we will prove, the proportion of the left edge weight at times $4n$ is a martingale. This motivates the definition
\begin{equation*}
    M^{\rm alt}_n:=N_{4n}\,,\,\forall\,n\geq 0\,.
\end{equation*}
$M^{\rm alt}_n$ is the proportion of the weight of the left edge $\{-1,0\}$ at times $(4n)_{n \geq 0}$.

\begin{theorem}\label{thm:main-alternating} Consider the model where the walkers jump alternately with the correspondig sequence $(N_n)_{n\geq 0}$  of random variables. Then
\begin{enumerate}[(i)]
    \item \label{thm:altwalkersmart} $(M_n^{\rm alt})_{n\geq 0}$ is a bounded martingale and we define 
    \begin{equation*}
        M^{\rm alt}_\infty:=\lim_{n\to\infty}M_n^{\rm alt}\,.
    \end{equation*}
    We have that $\lim_{n\to\infty} |N_n-M_\infty^{\rm alt}|=0$ almost surely.
    \item \label{thm:altwalkersclt} The following convergence holds in the sense of almost sure conditional convergence w.r.t.~the filtration $(\calG_n)_{n\geq 0}$:
    \begin{align}
      \label{eq:altwalkersclt}
    \sqrt{n} \left(M^{\rm alt}_n - M^{\rm alt}_{\infty}\right) \to \calN\left(0, \frac{1}{2}M^{\rm alt}_\infty\left(1 - M^{\rm alt}_\infty\right)\right)\,.
    \end{align}
    As a consequence of \autoref{lemma:weakconv}, the following convergence holds in the weak sense:
    \begin{align*}
    \sqrt{n} \left(N_n - M^{\rm alt}_{\infty}\right) \to \overline \calN\left(0, \frac{1}{2}M^{\rm alt}_\infty\left(1 - M^{\rm alt}_\infty\right)\right) \,.
    \end{align*}
    \item \label{thm:altwalkersatoms} $\Prb{M^{\rm alt}_\infty = a} = 0$ for any $a \in \left[0, 1\right]$. As a consequence, $M^{\rm alt}_\infty$ is not deterministic.\qedhere
\end{enumerate}
\end{theorem}
We furthermore believe that the following is true.
\begin{conjecture}
  \label{conj:two_pl_urn_limit}
  $M_{\infty}^{\rm alt} \in \left[0, 1\right]$ has a density w.r.t.~the Lebesgue measure on $\left[0, 1\right]$.
  $M_{\infty}^{\rm alt}$ is not Beta-distributed.
\end{conjecture}

For the case of two walkers and random walker selection, we define the following sequence of stopping times: $\tau_0 = 0$ and, recursively,
\begin{equation*}
\tau_n = \inf\left\{ m > \tau_{n-1} : X^{\left(1\right)}_m = X^{\left(2\right)}_m = 0 \right\}\,.    
\end{equation*}
Furthermore, we define the sequence $(M_n^{\rm rd})_{n \geq 0}$ as follows:
\begin{equation*}
    M^{\rm rd}_n=N_{\tau_n}\,,\,\forall\,n\geq 0\,.
\end{equation*}
 
\begin{theorem}\label{thm:randwalkersmain} Consider the model where the walker to move is selected randomly with the correspondig sequence $(N_n)_{n\geq 0}$  of random variables. Then
\begin{enumerate}[(i)]
    \item \label{thm:randwalkersmart} $(M_n^{\rm rd})_{n\geq 0}$ is a bounded martingale and we define 
    \begin{equation*}
        M^{\rm rd}_\infty:=\lim_{n\to\infty}M_n^{\rm rd}\,.
    \end{equation*}
    \item \label{thm:randwalkersmeet} The sequence of stopping times $(\tau_n)_{n\geq 0}$ satisfies the following:
    \begin{align*}
      \Prb{\tau_{n+1} - \tau_n = 2l} &= 2^{-l} \textrm{ for all } l \geq 1  \,, \\
      \Ex{\tau_{n+1} - \tau_n } &= 4 \,.
    \end{align*}
    \item \label{thm:randwalkersconv} $\lim_{n\to\infty} N_n=M_\infty^{\rm rd}$ almost surely.\qedhere
\end{enumerate}
\end{theorem}
We prove \autoref{thm:main-alternating} and \autoref{thm:randwalkersmain}
in \autoref{sec:proof:alternating} and \autoref{sec:proof:random_selection}.

\subsection{Proofs: alternating walkers (\autoref{thm:main-alternating})} \label{sec:proof:alternating}

\begin{oproof}[\autoref{thm:main-alternating} \ref{thm:altwalkersmart}]\label{sec:alternating:martingale}
    In this section, we write $M_n$ instead of $M_n^{\rm alt}$.
  We first prove that $(M_n)_{n\geq 0}$ is a martingale. We assume that at time $4n$, the edge weights are given by $w\left(4n, -1\right) = a, w\left(4n, 0\right) = b$.
  Recall that the walkers move alternately and meet at the center node every four steps. Therefore, on the time interval $\{4n+1,4n+2,4n+3,4n+4\}$, the following possibilities occur:
  \begin{itemize}
      \item $X^{(1)}$ and $X^{(2)}$ move alternately from $0$ to $-1$ and back with probability
      \begin{equation*}
          \frac{a}{a + b} \frac{a + 1}{a + b + 1}
      \end{equation*}
      and the new edge weights are $w\left(4n+4, -1\right) = a+4, w\left(4n+4, 0\right) = b$.
      \item $X^{(1)}$ and $X^{(2)}$ move alternately from $0$ to $1$ and back with probability
      \begin{equation*}
          \frac{b}{a + b} \frac{b + 1}{a + b + 1}
      \end{equation*}
      and the new edge weights are $w\left(4n+4, -1\right) = a, w\left(4n+4, 0\right) = b+4$.
      \item $X^{(1)}$ moves from $0$ to $1$ and $X^{(2)}$ moves from $0$ to $-1$. $X^{(1)}$ and $X^{(2)}$ move back to $0$ with probability $1$.  This event happens with probability
      \begin{equation*}
          \frac{b}{a + b} \frac{a}{a + b + 1}
      \end{equation*}
      and the new edge weights are $w\left(4n+4, -1\right) = a+2, w\left(4n+4, 0\right) = b+2$.
      \item $X^{(1)}$ moves from $0$ to $-1$ and $X^{(2)}$ moves from $0$ to $1$. After that $X^{(1)}$ and $X^{(2)}$ move back to $0$ with probability $1$.  This event happens with probability
      \begin{equation*}
          \frac{a}{a + b} \frac{b}{a + b + 1}
      \end{equation*}
      and the new edge weights are $w\left(4n+4, -1\right) = a+1, w\left(4n+4, 0\right) = b+1$.
  \end{itemize}
  Therefore,
  \begin{align*}
    &\Ex{\frac{w\left(4n + 4,-1\right)}{w\left(4n + 4,-1\right) + w\left(4n + 4,0\right)} \midd| w\left(4n, -1\right) = a, w\left(4n, 0\right) = b} \\
    &= \frac{a}{a + b} \frac{a + 1}{a + b + 1} \frac{a + 4}{a + b + 4}
    + \frac{a}{a + b} \frac{b}{a + b + 1} \frac{a + 2}{a + b + 4} \\
    &\hphantom{\;=\;} + \frac{b}{a + b} \frac{a}{a + b + 1} \frac{a + 2}{a + b + 4}
    + \frac{b}{a + b} \frac{b + 1}{a + b + 1} \frac{a}{a + b + 4} \\
    &= \frac{a}{a + b} \cdot \frac{\left(a + 1\right)\left(a + 4\right) + 2b\left(a + 2\right) + b\left(b + 1\right)}{\left(a + b + 1\right)\left(a + b + 4\right)}
    =: \frac{a}{a + b}E_{ab}\,. 
  \end{align*}
  To see that $E_{ab}=1$, it is enough to count how many terms of the form $a^2$ (respectively, $b^2$, $b$, $a$, and $ab$) appear in the numerator and in the denominator.

  The martingale $(M_n)_{n\geq 0}$ is bounded and therefore $M_n=\frac{w\left(4n,-1\right)}{4n+2}$
  converges a.s.\ to a random variable $M_\infty$. We claim that $N_n = \frac{w\left(n,-1\right)}{n+2}$ also converges to $M_\infty$. Indeed, for $i \in \left\{1,2,3\right\}$, we have
  \begin{align*}
  \lim_{n\to\infty}\left|\frac{w\left(4n+i,-1\right)}{4n+i+2}-\frac{w\left(4n,-1\right)}{4n+2}\right| = 0\,,
  \end{align*}
  because $|w\left(4n+i,-1\right)-w\left(4n,-1\right)|\leq i\,.$
\end{oproof}

\begin{oproof}[\autoref{thm:main-alternating} \ref{thm:altwalkersclt}]\label{sec:alternating:item2}
In this section, we write again $M_n$ instead of $M_n^{\rm alt}$. The proof of \autoref{thm:main-alternating} \ref{thm:altwalkersclt} uses \autoref{thm:a3}. Indeed, we only have to verify
the two assumptions of \autoref{thm:a3}, with $V = \frac{1}{2}M_\infty(1 - M_\infty)$.

\noindent\textbf{First step / condition \ref{thm:a3cond1}:} it holds that
  \begin{align}\label{eq:first_step}
    \Ex{\sup_{n > 1} \sqrt{n}\left|M_{n-1} - M_n\right|} &\leq 1< \infty \,.
  \end{align}
  Since $\left|N_n - N_{n-1}\right|$ is bounded by $\frac{1}{\left(n-1\right) + 2}$, we get that
  $\left|M_{n-1} - M_n\right| \leq \frac{4}{4\left(n-1\right) + 2}$. Thus
  \begin{align*}
    \sqrt{n}\left|M_{n-1} - M_n\right| \leq \frac{4\sqrt{n}}{4\left(n-1\right) + 2} < 1
  \end{align*}
  which implies that~\eqref{eq:first_step} is true. 
  
  \noindent\textbf{Second step / condition \ref{thm:a3cond2}:} the following limit holds almost surely:
    \begin{equation}\label{eq:secon_step}
        n\sum_{j\geq n}(M_{j-1}-M_j)^2\to \dfrac12 M_\infty(1-M_\infty)\,.
    \end{equation}
To prove~\eqref{eq:secon_step}, we use~\autoref{lem:a2} with $X_n := n^2\left(M_{n-1} - M_n\right)^2$. Condition \ref{lem:a2cond1} of \autoref{lem:a2} is a direct consequence of
$\left|M_{n-1} - M_n\right| \leq \frac{4}{4\left(n-1\right) + 2}$.

  \noindent\textbf{Our goal (showing condition \ref{lem:a2cond2} of \autoref{lem:a2}):} to show that the following holds almost surely, where $\calF_n := \calG_{4n} = \sigma(N_0, \ldots, N_{4n})$:
    \begin{equation*}
        \Ex{n^2\left(M_{n-1} - M_n\right)^2\midd|\calF_{n-1}} \to \frac{1}{2}M_\infty\left(1 - M_\infty\right)\,.
    \end{equation*}
  We claim that
  \begin{align}\label{eq:convclaim}
    \Ex{M_n^2 \midd| \calF_{n-1}} = \left(1 - \frac{1}{2n^2}\right) M_{n-1}^2 + \frac{1}{2n^2} M_{n-1} + \Delta_n \quad \textrm{ with } \quad n^2\Delta_n \overset{\textrm{a.s.}}{\longrightarrow} 0\,.
  \end{align}
  Assuming the claim it follows that
  \begin{align*}
    \Ex{n^2\left(M_{n-1} - M_n\right)^2\midd|\calF_{n-1}}
    &= n^2 \Ex{M_n^2 \midd| \calF_{n-1}} - n^2 M_{n-1}^2 \\
    &= \left(n^2 - \frac{1}{2}\right) M_{n-1}^2 + \frac{1}{2} M_{n-1} + n^2\Delta_n - n^2 M_{n-1}^2 \\
    &= \frac{1}{2} M_{n-1} \left(1 - M_{n-1}\right) + n^2\Delta_n
    \overset{\textrm{a.s.}}{\longrightarrow} \frac{1}{2} M_{\infty} \left(1 - M_{\infty}\right)\,,
  \end{align*}
  as we wanted to prove.

  We proceed to show~\eqref{eq:convclaim}. By definition,
  \begin{align*}
    n^2\Delta_n &= n^2\left(\Ex{M_n^2 \midd| \calF_{n-1}} - \left(1 - \frac{1}{2n^2}\right) M_{n-1}^2 - \frac{1}{2n^2} M_{n-1}\right)\,.
  \end{align*}
  We condition on $w\left(4\left(n-1\right), -1\right) = a, w\left(4\left(n-1\right), 0\right) = b$, similarly to the proof of \autoref{thm:main-alternating}~\ref{thm:altwalkersmart}, and get
  \begin{align*}
    \Ex{M_n^2 \midd| \calF_{n-1}}
    =\; &\frac{a}{a+b} \cdot \frac{a+1}{a+b+1} \cdot \frac{\left(a+4\right)^2}{\left(a+b+4\right)^2}
    + 2 \cdot \frac{a}{a+b} \cdot \frac{b}{a+b+1} \cdot \frac{\left(a+2\right)^2}{\left(a+b+4\right)^2} \\
    &+ \frac{b}{a+b} \cdot \frac{b+1}{a+b+1} \cdot \frac{a^2}{\left(a+b+4\right)^2} \\
    =\; & \frac{a}{a+b} \cdot \left(
      \frac{a+1}{a+b+1} \cdot \frac{\left(a+4\right)^2}{\left(a+b+4\right)^2}
      + \frac{2b}{a+b+1} \cdot \frac{\left(a+2\right)^2}{\left(a+b+4\right)^2} \right. \\
    & \hphantom{\frac{a}{a+b} \cdot ()}\left.
      + \; \frac{b\left(b+1\right)}{a+b+1} \cdot \frac{a}{\left(a+b+4\right)^2}
    \right)\,.
  \end{align*}
  Furthermore, under the same conditioning, and using that $n = \frac{a + b + 2}{4}$,
  \begin{align*}
    \left(1 - \frac{1}{2n^2}\right) M_{n-1}^2
    &= \left(1 - \frac{8}{\left(a + b + 2\right)^2}\right) \frac{a^2}{\left(a + b\right)^2} \\
    \frac{1}{2n^2} M_{n-1}
    &= \frac{8}{\left(a + b + 2\right)^2} \cdot \frac{a}{a + b} \,.
  \end{align*}
  Putting this together (see \href{https://bit.ly/3I8rEWL}{https://bit.ly/3I8rEWL}), we arrive at
  \begin{align*}
    \Delta_n &= - \underbrace{\frac{8}{\left(a + b + 2\right)^2}}_{\frac{1}{2n^2}} \cdot \underbrace{\frac{a}{a + b} \cdot \frac{b}{a + b}}_{M_{n-1}\left(1 - M_{n-1}\right)} \cdot \frac{3a^2 + 6ab + 12a + 3b^2 + 12b + 8}{\left(a + b + 1\right)\left(a + b + 4\right)^2}\,, \\
    n^2 \Delta_n &= -\frac{1}{2} M_{n-1}\left(1 - M_{n-1}\right) \cdot \frac{3a^2 + 6ab + 12a + 3b^2 + 12b + 8}{\left(a + b + 1\right)\left(a + b + 4\right)^2} \,.
  \end{align*}
  Using again $n = \frac{a + b + 2}{4}$, we can rewrite the fraction on the right (see \href{http://bit.ly/3ICUJLz}{http://bit.ly/3ICUJLz})  to get to the final conclusion
  \begin{align*}
    n^2\Delta_n &= n^2\left(\Ex{M_n^2 \midd| \calF_{n-1}} - \left(1 - \frac{1}{2n^2}\right) M_{n-1}^2 - \frac{1}{2n^2} M_{n-1}\right) \\
    &= -\frac{1}{2} M_{n-1}\left(1 - M_{n-1}\right) \left(\frac{7}{9\left(n + \frac{1}{2}\right)} - \frac{1}{6\left(n + \frac{1}{2}\right)^2} - \frac{1}{36\left(n - \frac{1}{4}\right)}\right) \overset{n \to \infty}{\longrightarrow} 0\,. \qedhere
  \end{align*}
\end{oproof}

\begin{oproof}[\autoref{thm:main-alternating} \ref{thm:altwalkersatoms}]\label{sec:proof:thm:main-alternating-3}
In this section, we write again $M_n$ instead of $M_n^{\rm alt}$. Assuming \autoref{thm:main-alternating} \ref{thm:altwalkersclt}, we now prove \autoref{thm:main-alternating} \ref{thm:altwalkersatoms}.

We first prove that $\mathbb P(M_\infty=a)=0$ if $a\in (0,1)$. Afterwards we prove that $\mathbb P(M_\infty=0)=0$. The proof of $\mathbb P(M_\infty=1)=0$ follows by symmetry since
we can swap the labels of the left and the right node to see that the law of the proportion of the right edge weight is
equal to the law of the proportion of the left edge weight.

We follow the proof of \cite[Corollary 3]{hypergrru}.
Set $A = \left\{M_\infty = a\right\}$ for $a \in (0, 1)$. First observe that, as a consequence of \autoref{thm:main-alternating} \ref{thm:altwalkersclt} (recall that $\calF_n := \calG_{4n} = \sigma(N_0, \ldots, N_{4n})$),
\begin{equation*}
  \Ex{\exp\left(\ivar t \sqrt{n} \left(M_n - M_\infty\right) \right) \midd| \calF_{n}} \overset{\textrm{a.s.}}{\longrightarrow} \exp\left(-\frac{t^2 \cdot \frac{1}{2}M_\infty\left(1 - M_\infty\right)}{2}\right)\,.
\end{equation*}
We also have that $I_n := \Ex{\1_A \midd| \calF_{n}} \to \1_A$ a.s..
Therefore, for any $t\in \mathbb R$,
\begin{gather}
\label{eq:conv_in_charact}
\begin{split}
  \Ex{I_n \exp\left(\ivar t \sqrt{n} \left(M_n - M_\infty\right) \right) \midd| \calF_{n}}
  &= I_n \Ex{\exp\left(\ivar t \sqrt{n} \left(M_n - M_\infty\right) \right) \midd| \calF_{n}} \\
  &\overset{\textrm{a.s.}}{\longrightarrow} \1_A \cdot \exp\left(-\frac{t^2 \cdot \frac{1}{2}M_\infty\left(1 - M_\infty\right)}{2}\right)\,.
\end{split}
\end{gather}

Since $\1_A - I_n\to 0$ almost surely, it follows from the dominated convergence theorem for conditional expectation (see \cite[Lemma A.1]{hypergrru}) that
\begin{align}
\label{eq:conv_no_diff}
\Ex{\left(\1_A - I_n\right)\exp\left(\ivar t \sqrt{n} \left(M_n - M_\infty\right) \right) \midd| \calF_{n}} \overset{\textrm{a.s.}}{\longrightarrow} 0\,.
\end{align}
Combining \eqref{eq:conv_in_charact} and \eqref{eq:conv_no_diff} yields
\begin{align}\label{eq:1A}
\Ex{\1_A \exp\left(\ivar t \sqrt{n} \left(M_n - M_\infty\right) \right) \midd| \calF_{n}} \overset{\textrm{a.s.}}{\longrightarrow} \1_A \cdot \exp\left(-\frac{t^2 \cdot \frac{1}{2}M_\infty\left(1 - M_\infty\right)}{2}\right)\,.
\end{align}
Therefore,
\begin{align*}
\1_A \cdot \exp\left(-\frac{t^2 \cdot \frac{1}{2}a\left(1 - a\right)}{2}\right)
&= \1_A \cdot \exp\left(-\frac{t^2 \cdot \frac{1}{2}M_\infty\left(1 - M_\infty\right)}{2}\right) \\
\text{(using~\eqref{eq:1A}) }&= \lim_n \Ex{\1_A \exp\left(\ivar t \sqrt{n} \left(M_n - M_\infty\right) \right) \midd| \calF_{n}} \\
\text{(definition of $A$) }&= \lim_n \Ex{\1_A \exp\left(\ivar t \sqrt{n} \left(M_n - a\right) \right) \midd| \calF_{n}} \\
\text{(measurability of $M_n$ w.r.t.~$\calF_n$) } &= \lim_n I_n \exp\left(\ivar t \sqrt{n} \left(M_n - a\right) \right)\\
\text{(using that $I_n \to \1_A$)}&= \lim_n \1_A \exp\left(\ivar t \sqrt{n} \left(M_n - a\right) \right)\,.
\end{align*}
As a consequence,
\begin{align*}
\1_A &=  \lim_n \1_A \left|\exp\left(\ivar t \sqrt{n} \left(M_n - a\right) \right) \right| \\
&=\left| \lim_n \1_A \exp\left(\ivar t \sqrt{n} \left(M_n - a\right) \right) \right|\\
&= \1_A \cdot \underbrace{\exp\left(-\frac{t^2 \cdot \frac{1}{2}a\left(1 - a\right)}{2}\right)}_{< 1 \textrm{ since } a(1-a) > 0}\,.
\end{align*}
Therefore, $\Prb{A} = 0$.

Now, it remains to prove that $\mathbb P(M_\infty=0)=0$. The idea of the proof is to fix $\epsilon>0$ and to split the event $\{M_\infty=0\}$ depending on whether $M_n\leq \epsilon\sqrt{n}$ or $M_n> \epsilon\sqrt{n}$. When $M_n> \epsilon\sqrt{n}$ we use the CLT convergence to prove that the probability goes to zero and when  $M_n\leq \epsilon\sqrt{n}$ we use explicit bounds to show that the probability goes to zero. 

More specifically, for any $\epsilon>0$,
  \begin{equation}
    \label{eq:prb_atom_0}
      \Prb{M_\infty = 0} = \Prb{M_\infty = 0, M_n \leq \frac{\epsilon}{\sqrt{n}}} + \Prb{M_\infty = 0, M_n > \frac{\epsilon}{\sqrt{n}}}
  \end{equation}
  By \autoref{thm:main-alternating} \ref{thm:altwalkersclt}, i.e.~\eqref{eq:altwalkersclt}, we know that $\sqrt{n}\left(M_n - M_\infty\right)$ converges in law to the Dirac measure in $0$ if we condition on $M_\infty = 0$.
  As a consequence, $\Prb{M_\infty = 0, \sqrt{n}\left|M_n - M_\infty\right| > \epsilon} \to 0$ for $n \to \infty$. This implies in particular that the rightmost term in \eqref{eq:prb_atom_0} converges to $0$ for $n \to \infty$.
    On the other hand,
  \begin{equation*}
       \Prb{M_\infty = 0, M_n \leq \frac{\epsilon}{\sqrt{n}}}\leq  \Prb{M_n \leq \frac{\epsilon}{\sqrt{n}}}=: p_{n,\epsilon}\,.
  \end{equation*}
 The proof is then finished once we prove the following claim:
  \begin{align}
    \label{eq:pnc_lim_inequ}
    \limsup_{n \to \infty} p_{n,\epsilon} &\leq 9 \sqrt{\epsilon}\,.
  \end{align}
  \eqref{eq:pnc_lim_inequ} implies that the first term in \eqref{eq:prb_atom_0} is upper bounded by a number arbitrarily close to $0$ in the limit. We dedicate the remaining section to the proof of~\eqref{eq:pnc_lim_inequ}.

  Notice that
  \begin{align*}
    p_{n,\epsilon} := \Prb{M_n \leq \frac{\epsilon}{\sqrt{n}}} = \Prb{\frac{w(4n, -1)}{w(4n, -1) + w(4n, 0)} \leq \frac{\epsilon}{\sqrt{n}}} = \Prb{w(4n, -1) \leq \frac{\epsilon\left(4n+2\right)}{\sqrt{n}}}\,.
  \end{align*}
  Denote by $l_n$ the number of crossings of the left edge in the outwards direction in the first $4n$ steps,
  and by $r_n$ the number of crossings of the right edge in the outwards direction. It follows that $l_n + r_n = 2n$, and $w(4n, -1) = 1 + 2l_n$. Observe also that $\frac{\epsilon}{\sqrt{n}} \leq \frac{1}{2}$ for $n$ large enough. As a consequence,
  \begin{align*}
    w(4n, -1) \leq \frac{\epsilon\left(4n+2\right)}{\sqrt{n}} &\implies 2l_n \leq \frac{\epsilon\left(4n+2\right)}{\sqrt{n}} - 1 \implies l_n \leq \frac{\epsilon\left(2n+1\right)}{\sqrt{n}} - \frac{1}{2} \\
    &\implies l_n \leq \left\lfloor 2\epsilon\sqrt{n} \right\rfloor\,.
  \end{align*}
  Therefore,
  \begin{equation*}
      p_{n,\epsilon} \leq \sum_{j = 0}^{\left\lfloor 2\epsilon\sqrt{n} \right\rfloor} \Prb{l_n = j}\,.
  \end{equation*}

  We claim that
  \begin{equation}\label{eq:claim1}
      \Prb{l_n = j}\leq \binom{2n}{j}
      \cdot
      \frac{1 \cdot 3 \cdot 5 \cdot \ldots \cdot (2j - 1) \cdot 1 \cdot 3 \cdot 5 \cdot \ldots \cdot (2(2n-j) - 1)}{2 \cdot 3 \cdot 6 \cdot 7 \cdot 10 \cdot 11 \cdot \ldots \cdot (4n - 2) \cdot (4n - 1)}
      % \dfrac{\prod_{i=1}^{j}(2i-1)\prod_{i=1}^{2n-j}(2i-1)}{\prod_{i=1}^n(4i-2)(4i-1)}\,.
  \end{equation}
The reader is advised to write examples to convince herself/himself that~\eqref{eq:claim1} holds true. The binomial coefficient counts the
number of possible walker movement sequences of length $4n$ which end with $l_n = j$. The sequence in the denominator
of the following fraction is just the sequence of total edge weights observed at the points in time
at which one of the walkers moves from the node in the center to one of the outer nodes. For the
sequence in the numerator, we just take an upper bound on the possible
edge weights of the left and right edges before a walker crosses the respective
edge in outwards direction.

Reordering the terms, and writing out the binomial coefficient, we get
  \begin{align*}
    p_{n,\epsilon} &\leq \sum_{j = 0}^{\left\lfloor 2\epsilon\sqrt{n} \right\rfloor} \frac{1 \cdot 2 \cdot 3 \cdot 4 \cdot \ldots \cdot 2n}{2 \cdot 3 \cdot 6 \cdot 7 \cdot \ldots \cdot \left(4n - 1\right)} \cdot \frac{1 \cdot 3 \cdot 5 \cdot \ldots \cdot \left(2j - 1\right)}{1 \cdot 2 \cdot 3 \cdot \ldots \cdot j} \cdot \frac{1 \cdot 3 \cdot 5 \cdot \ldots \cdot \left(2\left(2n - j\right) - 1\right)}{1 \cdot 2 \cdot 3 \cdot \ldots \cdot \left(2n-j\right)} \\
    &= 2^{-2n} \cdot \prod_{i=1}^{n} \frac{4i}{4i - 1} \cdot \sum_{j = 0}^{\left\lfloor 2\epsilon\sqrt{n} \right\rfloor} \prod_{i = 1}^{j} \frac{2i - 1}{i} \cdot \prod_{i = 1}^{2n - j} \frac{2i - 1}{i} \\
    &= \prod_{i=1}^{n} \frac{4i}{4i - 1} \cdot \sum_{j = 0}^{\left\lfloor 2\epsilon\sqrt{n} \right\rfloor} \prod_{i = 1}^{j} \frac{2i - 1}{2i} \cdot \prod_{i = 1}^{2n - j} \frac{2i - 1}{2i}\,.
  \end{align*}

  To bound the products, we do the following: let $m \in \left\{2,4\right\}$ and note that $\log\left(1-x\right) \geq -x-x^2$ for $x \in \left[0,\frac{1}{2}\right]$.
  \begin{align*}
    \prod_{i=1}^{j} \frac{m \cdot i}{m \cdot i - 1}
    &= \exp\left(-\sum_{i=1}^{j} \log\left(\frac{m \cdot i - 1}{m \cdot i}\right)\right)
    = \exp\left(-\sum_{i=1}^{j} \log\left(1 - \frac{1}{m \cdot i}\right)\right) \\
    &\leq \exp\left(\sum_{i=1}^{j} \left(\frac{1}{m \cdot i} + \frac{1}{m^2 \cdot i^2}\right)\right)
    \leq \exp\left(\frac{1}{m}\left(\log\left(j\right) + 1\right) + \frac{\pi^2}{6m^2}\right)
    \leq 3j^{\frac{1}{m}}\,.
  \end{align*}
  Furthermore, using that $\log\left(1-x\right) \leq -x$, we get that
  \begin{align*}
    \prod_{i=1}^{j} \frac{m \cdot i}{m \cdot i - 1}
    &\geq \exp\left(\sum_{i=1}^{j} \frac{1}{m \cdot i}\right) \geq \exp\left(\frac{1}{m}\log\left(j\right)\right) = j^\frac{1}{m} \\
    \implies \prod_{i=1}^{j} \frac{m \cdot i - 1}{m \cdot i} &\leq j^{-\frac{1}{m}}\,.
  \end{align*}
  Hence (we let again $n$ large enough such that $\frac{\epsilon}{\sqrt{n}} \leq \frac{1}{2}$)
  \begin{align*}
    p_{n,\epsilon}
    &\leq \prod_{i=1}^{n} \frac{4i}{4i - 1} \cdot \sum_{j = 0}^{\left\lfloor 2\epsilon\sqrt{n} \right\rfloor} \prod_{i = 1}^{j} \frac{2i - 1}{2i} \cdot \prod_{i = 1}^{2n - j} \frac{2i - 1}{2i}
    \leq 3n^{\frac{1}{4}} \cdot \sum_{j = 0}^{\left\lfloor 2\epsilon\sqrt{n} \right\rfloor} \frac{1}{\sqrt{j + \mathbbm{1}_{j=0}}} \cdot \frac{1}{\sqrt{2n - j}} \\
    &\leq 3n^{\frac{1}{4}} \cdot \frac{1}{\sqrt{2n - \left\lfloor 2\epsilon\sqrt{n} \right\rfloor}} \left(1 + \sum_{j = 1}^{\left\lfloor 2\epsilon\sqrt{n} \right\rfloor} \frac{1}{\sqrt{j}}\right)
    \leq 3n^{\frac{1}{4}} \cdot \frac{1}{\sqrt{n}} \left(1 + \int_{0}^{2\epsilon\sqrt{n}} \frac{1}{\sqrt{j}} \dx{j}\right) \\
    &= 3n^{\frac{1}{4}} \cdot \frac{1}{\sqrt{n}} \left(1 + 2\sqrt{2\epsilon\sqrt{n}}\right)
    \leq 3n^{-\frac{1}{4}} + 6\sqrt{2}\sqrt{\epsilon} \cdot n^{\frac{1}{4}} \cdot n^{-\frac{1}{2}} \cdot n^\frac{1}{4} \leq 3n^{-\frac{1}{4}} + 9 \sqrt{\epsilon}\,.
  \end{align*}
  Therefore, $\limsup_{n \to \infty} p_{n,\epsilon} \leq 9 \sqrt{\epsilon}$.
\end{oproof}

\subsection{Proofs: random walker selection (\autoref{thm:randwalkersmain})}\label{sec:proof:random_selection}

Consider next the case with random walker selection. \autoref{lem:two_pl_urn_rand_meeting_time} shows that the expected time to meet again in the
middle, if both walkers start in the center, is $4$, just as in the case of alternating walkers. Of course, the difference
now is that the next meeting time is random. 

We want to prove \autoref{thm:randwalkersmain} \ref{thm:randwalkersconv}, i.e.~that the proportion of the left edge weight,
$N_n$,
converges a.s.~for $n \to \infty$, in analogy to the case with alternating walkers, as well as in analogy with the case of a single walker which corresponds to a standard
P\'olya urn. The proof is based on the fact that the proportion of the left edge weight is a martingale at the times
$\tau_n$ (recall that $\tau_n = \inf\{ m > \tau_{n-1} : X^{\left(1\right)}_m = X^{\left(2\right)}_m = 0 \}$), as we will see in \autoref{lem:two_pl_urn_mart_rand} (which also proves \autoref{thm:randwalkersmain} \ref{thm:randwalkersmart}).

In order to prove \autoref{thm:randwalkersmain} \ref{thm:randwalkersconv}, we need two lemmata. The main reason for convergence is the
martingale property of the $N_{\tau_n}$, and the remaining section will mostly
be devoted to showing \autoref{lem:two_pl_urn_mart_rand}. In order to make the step from the times $\tau_n$ to $n$,
it is enough to see that the probability of $\left\{\tau_{n+1} - \tau_n = 2l\right\}$ decays exponentially in $l$
in \autoref{lem:two_pl_urn_rand_meeting_time} (which proves \autoref{thm:randwalkersmain} \ref{thm:randwalkersmeet}), which allows us to conclude with a Borel-Cantelli argument.

In this section, we write $M_n$ instead of $M^{\rm rd}_n$, i.e.
\begin{equation*}
    M_n = M^{\rm rd}_n = N_{\tau_n} = \frac{w\left(\tau_n,-1\right)}{w\left(\tau_n,-1\right) + w\left(\tau_n,0\right)}\,.
\end{equation*}
\begin{lemma}
  \label{lem:two_pl_urn_mart_rand}
  The sequence of random variables $(M_n)_{n\geq 0}$ together with its natural filtration is a (bounded) martingale.
\end{lemma}

\begin{lemma}
  \label{lem:two_pl_urn_rand_meeting_time}
  For any $n \geq 0$ and any $l \geq 1$, it holds that $\Prb{\tau_{n+1} - \tau_n = 2l} = 2^{-l}$
  and $\Ex{\tau_{n+1} - \tau_n} = 4$.
\end{lemma}

Assuming \autoref{lem:two_pl_urn_mart_rand} and \autoref{lem:two_pl_urn_rand_meeting_time}, we can prove \autoref{thm:randwalkersmain} \ref{thm:randwalkersconv}:

\begin{oproof}[\autoref{thm:randwalkersmain} \ref{thm:randwalkersconv}]
  Recall that
  \begin{align*}
    N_n := \frac{w\left(n,-1\right)}{w\left(n,-1\right) + w\left(n,0\right)}\,, \qquad
    M_n := N_{\tau_n} = \frac{w\left(\tau_n,-1\right)}{w\left(\tau_n,-1\right) + w\left(\tau_n,0\right)}\,.
  \end{align*}
  By \autoref{lem:two_pl_urn_mart_rand}  there exists a random variable $M_\infty$ such that
  \begin{equation}
      M_{\infty} := \lim\limits_{n \to \infty} M_n\,.
  \end{equation} 
  Let $\varepsilon > 0$ and consider the events
  \begin{align*}
    A_n := \left\{ \left|N_m - M_{\infty}\right| > \varepsilon \textrm{ for some } m \in \left[\tau_n, \tau_{n+1}\right] \right\}
  \end{align*}
  It is sufficient to show that only finitely many of the events $A_n$ can occur a.s.
  Further set
  \begin{align*}
    B_n := \left\{ \left|N_m - M_n\right| > \frac{\varepsilon}{2} \textrm{ for some } m \in \left[\tau_n, \tau_{n+1}\right] \right\}
  \end{align*}
  Since $M_n$ converges a.s.~by \autoref{lem:two_pl_urn_mart_rand}, there is some (random) $N \in \bbN$ such that for all $n \geq N$,
  it holds that $\left|M_n - M_{\infty}\right| < \frac{\varepsilon}{2}$. For $n \geq N$, the occurrence of $A_n$
  implies that $B_n$ occurs as well, so it is sufficient to show that only finitely many of the events
  $B_n$ can occur.
  
  Now, at time $\tau_n$, the random walkers must have moved at least $2n$ times, so
  $w\left(\tau_n,-1\right) + w\left(\tau_n,0\right) \geq 2n$. If $\left|N_m - M_n\right| = \left|N_m - N_{\tau_n}\right| > \frac{\varepsilon}{2}$
  for some $m \in \left[\tau_n, \tau_{n+1}\right]$, it is therefore necessary that at least $\varepsilon n$
  steps were made by the walkers between time $\tau_n$ and time $m$, since every step changes the value
  of $N_m$ by at most $\frac{1}{2n}$. Thus 
  \begin{align*}
    B_n \subseteq \left\{ \tau_{n+1} - \tau_n \geq \varepsilon n \right\}
  \end{align*}
  We have
  \begin{align*}
    \Prb{\tau_{n+1} - \tau_n \geq \varepsilon n}
    &= \sum_{l \geq \left\lceil \frac{\varepsilon n}{2} \right\rceil} \Prb{\tau_{n+1} - \tau_n = 2l}
    \overset{\textrm{\autoref{lem:two_pl_urn_rand_meeting_time}}}{=} \sum_{l \geq \left\lceil \frac{\varepsilon n}{2} \right\rceil} 2^{-l} \\
    &\leq 2^{-\frac{\varepsilon n}{2}} \sum_{l \geq 0} 2^{-l} = 2^{1 - \frac{\varepsilon n}{2}} \\
    \implies \sum_{n \geq 1} \Prb{\tau_{n+1} - \tau_n \geq \varepsilon n}
    &\leq \sum_{n \geq 1} 2^{1 - \frac{\varepsilon n}{2}} = 2 \cdot \sum_{n \geq 1} \left(2^{-\frac{\varepsilon}{2}}\right)^n < \infty
  \end{align*}
  By the Borel-Cantelli lemma, it follows that only finitely many of the events $\left\{ \tau_{n+1} - \tau_n \geq \varepsilon n \right\}$
  can occur, and therefore also only finitely many of the events $B_n$. This concludes the proof.
\end{oproof}

We proceed by first showing the exponential decay of $\Prb{\tau_{n+1} - \tau_n = 2l}$ in $l$. This is easy to
see since this probability can be calculated independently of any edge weights.

\begin{oproof}[\autoref{lem:two_pl_urn_rand_meeting_time} and \autoref{thm:randwalkersmain} \ref{thm:randwalkersmeet}]
  We consider a Markov chain consisting of three states and coupled with the edge-reinforced random walk.
  The Markov chain is in state $s_{\textrm{center}}$ if both walkers are in the center, in state $s_{\textrm{mixed}}$
  if one walker is in the center and the other in either of the two outer nodes, and in state $s_{\textrm{none}}$
  if none of the walkers is in the center. It is easy to verify that this is indeed a Markov chain
  with the following transition probabilities:
  
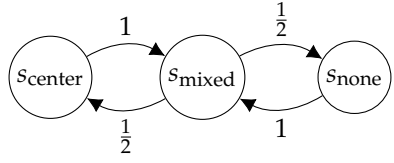
\begin{center}
    \begin{tikzpicture}
      \node[circle,draw=black,inner sep=0.7mm] (cen) at (-2, 0) {$s_{\textrm{center}}$};
      \node[circle,draw=black,inner sep=0.7mm] (mix) at (0, 0) {$s_{\textrm{mixed}}$};
      \node[circle,draw=black,inner sep=0.7mm] (non) at (2, 0) {$s_{\textrm{none}}$};
      \draw (cen) edge[-{Latex[length=2mm,width=2mm]},bend left] node[above] {$1$} (mix);
      \draw (mix) edge[-{Latex[length=2mm,width=2mm]},bend left] node[below] {$\frac{1}{2}$} (cen);
      \draw (mix) edge[-{Latex[length=2mm,width=2mm]},bend left] node[above] {$\frac{1}{2}$} (non);
      \draw (non) edge[-{Latex[length=2mm,width=2mm]},bend left] node[below] {$1$} (mix);
    \end{tikzpicture}
  \end{center}

  Let $n \geq 0$. At time $\tau_n$, both walkers are in the center and the Markov chain is therefore
  in state $s_{\textrm{center}}$. The time $\tau_{n+1} - \tau_n$ corresponds to the time needed to
  return again to the state $s_{\textrm{center}}$. At uneven time steps, the Markov chain will always
  be in state $s_{\textrm{mixed}}$, so $\Prb{\tau_{n+1} - \tau_n = 2l}$ corresponds to the probability
  that at every uneven time step, starting with $\tau_n + 1$ and up to $\tau_{n+1} - 3$, the chain transitions
  to the state $s_{\textrm{none}}$ and that the chain will finally go to $s_{\textrm{center}}$ in step $\tau_{n+1} - 1$.
  Since the probabilities of these transitions are $\frac{1}{2}$, we get $\Prb{\tau_{n+1} - \tau_n = 2l} = 2^{-l}$. Consequently,
  \begin{align*}
    \Ex{\tau_{n+1} - \tau_n} &= \sum_{l \geq 1} 2^{-l} \cdot 2l = \sum_{l \geq 1} l \cdot \left(\frac{1}{2}\right)^{l-1}
    = \frac{1}{\left(1 - \frac{1}{2}\right)^2} = 4 \qedhere
  \end{align*}
\end{oproof}

We next want to prove \autoref{lem:two_pl_urn_mart_rand}, i.e.~that the proportion of the left edge weight, that is
$\frac{w\left(\tau_n,-1\right)}{w\left(\tau_n,-1\right) + w\left(\tau_n,0\right)}$,
is a martingale, as it was in the previous case with alternating walkers and as in the P\'olya urn with the difference that
we look at the proportion not at every time step, but at certain stopping times. The proof is a bit more involved, but relies
only on basic calculations. The main idea is to find a recursive formula for the following quantities:
\begin{enumerate}[(1)]
  \item \label{term:two_pl_urn_rand_expec} We look at the expectation of the proportion of the left edge weight
  multiplied by the indicator of the event that the walkers need $2l$ steps to meet again:
  \begin{align*}
    \bbE_{a,b,l} := \Ex{\frac{w\left(\tau_{n+1},-1\right)}{w\left(\tau_{n+1},-1\right) + w\left(\tau_{n+1},0\right)} \cdot \mathbbm{1}_{\left\{\tau_{n+1} - \tau_n = 2l\right\}} \midd| w\left(\tau_n, -1\right) = a, w\left(\tau_n, 0\right) = b}
  \end{align*}
  \item \label{term:two_pl_urn_rand_proba} We also consider the probability that the last walker which
  returns to the center comes from the left node, again intersected with the event that the walkers need $2l$ steps to meet again:
  \begin{align*}
    L_n &:= \left\{ \textrm{the walker returning to the center at time } \tau_{n+1} \textrm{ comes from the left node} \right\} \\
    q_{a,b,l} &:= \Prb{L_n \textrm{ occurs and } \tau_{n+1} - \tau_n = 2l \midd| w\left(\tau_n, -1\right) = a, w\left(\tau_n, 0\right) = b}
  \end{align*}
\end{enumerate}

\begin{lemma}
  \label{lem:two_pl_urn_rand_recursion}
  We have the following recursive equations:
  \begin{align*}
    \bbE_{a,b,l+1} &= \frac{1}{2} \bbE_{a,b,l} + \frac{1}{\left(a+b+2l-1\right)\left(a+b+2l+2\right)} \left(\bbE_{a,b,l} - q_{a,b,l}\right) \\
    q_{a,b,l+1} &= \frac{1}{2} q_{a,b,l} + \frac{1}{4} \cdot \frac{a+b+2l}{a+b+2l-1} \left(\bbE_{a,b,l} - q_{a,b,l}\right) \qedhere
  \end{align*}
\end{lemma}

We will prove \autoref{lem:two_pl_urn_rand_recursion} later. First, we can complete the proof of \autoref{lem:two_pl_urn_mart_rand} as follows:

\begin{corollary}
  \label{cor:two_pl_urn_rand_recursion_solved}
  The expectation of the proportion of the left edge weight and the probability that the last walker
  returning to the center comes from the left coincide, and:
  \begin{align*}
    \bbE_{a,b,l} &= q_{a,b,l} = \frac{1}{2^l} \cdot \frac{a}{a + b} \qedhere
  \end{align*}
\end{corollary}

\begin{tproof}
  We use \autoref{lem:two_pl_urn_rand_recursion}. Let us first calculate $\bbE_{a,b,1}$ and $q_{a,b,1}$.
  There are four possible paths of length $2$ which end again with both walkers in the center: first,
  we choose which of the two walker moves, and this walker can then either move left or right and then
  back to the center. Since it is irrelevant which walker we choose in the beginning, we can disregard
  which walker moves. We thus get:
  \begin{itemize}
    \item With probability $\frac{a}{a + b} \cdot \frac{1}{2}$ the walker which moves in the first step
    moves left and is then chosen again to move back to the center in the next step. In this case, $L_n$
    occurs and the resulting edge weight ratio is $\frac{a + 2}{a + b + 2}$.
    \item With probability $\frac{b}{a + b} \cdot \frac{1}{2}$ the walker which moves in the first step
    moves right and is then chosen again to move back to the center in the next step. In this case, $L_n$
    does not occur and the resulting edge weight ratio is $\frac{a}{a + b + 2}$.
  \end{itemize}
  We see directly that $q_{a,b,1} = \frac{1}{2} \cdot \frac{a}{a + b}$ and that
  \begin{align*}
    \bbE_{a,b,1} = \frac{1}{2} \cdot \left(\frac{a}{a + b} \cdot \frac{a + 2}{a + b + 2} + \frac{b}{a + b} \cdot \frac{a}{a + b + 2}\right)
    = \frac{1}{2} \cdot \frac{a}{a + b}
  \end{align*}
  The remaining proof is now a simple induction using \autoref{lem:two_pl_urn_rand_recursion}, where
  it should be noted that $\bbE_{a,b,l} - q_{a,b,l} = 0$ under the induction assumption.
\end{tproof}

\begin{oproof}[\autoref{lem:two_pl_urn_mart_rand} and \autoref{thm:randwalkersmain} \ref{thm:randwalkersmart}]
  We have by \autoref{cor:two_pl_urn_rand_recursion_solved}:
  \begin{align*}
    &\Ex{\frac{w\left(\tau_{n+1},-1\right)}{w\left(\tau_{n+1},-1\right) + w\left(\tau_{n+1},0\right)} \cdot \mathbbm{1}_{\left\{\tau_{n+1} - \tau_n = 2l\right\}} \midd| w\left(\tau_n, -1\right) = a, w\left(\tau_n, 0\right) = b}
    = 2^{-l} \cdot \frac{a}{a + b} \\
    \implies
    &\Ex{\frac{w\left(\tau_{n+1},-1\right)}{w\left(\tau_{n+1},-1\right) + w\left(\tau_{n+1},0\right)} \midd| w\left(\tau_n, -1\right) = a, w\left(\tau_n, 0\right) = b}
    = \frac{a}{a+b} \cdot \sum_{l \geq 1} 2^{-l} = \frac{a}{a+b}
  \end{align*}
  as we wanted to prove.
\end{oproof}

This concludes the proof of \autoref{thm:randwalkersmain}, up to the proof of \autoref{lem:two_pl_urn_rand_recursion}, which follows now.
In order to prove \autoref{lem:two_pl_urn_rand_recursion}, we will use a recursive path construction technique.
We define a path of the two walkers as a sequence of the symbols
$1_{\textrm{l}}, 1_{\textrm{r}}, 2_{\textrm{l}}, 2_{\textrm{r}}$ which correspond to the first (respectively
second) walker moving left and right, where we assume that both walkers start in the center.
The set $\textrm{Path}_{2l}$ contains all the paths of length $2l$ (a sequence of $2l$ symbols) such
that the first time at which both walkers are in the center at the same time again is at the end of
the path. Note that any such path must be of even length since each walker can only be in the center
after having made an even number of movements.

For $\rho \in \textrm{Path}_{2l}$, we set $d\left(\rho\right)$ to be the number of traversals of the
left edge when the path $\rho$ is taken and we can now write
\begin{align*}
  \bbE_{a,b,l} &= \sum_{\rho \in \textrm{Path}_{2l}} \underbrace{\Prb{\textrm{the walkers move according to } \rho \midd| w\left(\tau_n, -1\right) = a, w\left(\tau_n, 0\right) = b}}_{=:\bbP_{\rho,a,b}} \cdot \frac{a + d\left(\rho\right)}{a + b + 2l}
\end{align*}

\begin{figure}[H]
  \begin{center}
    \begin{tabular}{cc}
      new step: \textcolor{sectionblue}{center} walker goes left & new step: \textcolor{sectionblue}{center} walker goes right \\
      ~ & ~ \\
      \begin{tikzpicture}
        \node[circle,fill=black,inner sep=0.7mm] (N0) at (0, 0) {};
        \node[circle,fill=black,inner sep=0.7mm] (N1) at (3, 0) {};
        \node[circle,fill=black,inner sep=0.7mm] (Nm1) at (-3, 0) {};
        \stickman{(-3,0.8)}{tumOrange}
        \stickman{(0,0.8)}{sectionblue}
        \node[below=1mm] at (N0) {$0$};
        \node[below=1mm] at (Nm1) {$-1$};
        \node[below=1mm] at (N1) {$1$};
        \draw (Nm1) -- node[above] {$a + d\left(\rho\right) - 1$} (N0) -- node[above] {$b + 2l - d\left(\rho\right)$} (N1);
        \draw[sectionblue] ([shift={(-0.4,0.7)}]N0.center) edge[bend right,-{Latex[length=2mm,width=2mm]}] ([shift={(0.4,0.7)}]Nm1.center);
      \end{tikzpicture}
      & \begin{tikzpicture}
        \node[circle,fill=black,inner sep=0.7mm] (N0) at (0, 0) {};
        \node[circle,fill=black,inner sep=0.7mm] (N1) at (3, 0) {};
        \node[circle,fill=black,inner sep=0.7mm] (Nm1) at (-3, 0) {};
        \stickman{(-3,0.8)}{tumOrange}
        \stickman{(0,0.8)}{sectionblue}
        \node[below=1mm] at (N0) {$0$};
        \node[below=1mm] at (Nm1) {$-1$};
        \node[below=1mm] at (N1) {$1$};
        \draw (Nm1) -- node[above] {$a + d\left(\rho\right) - 1$} (N0) -- node[above] {$b + 2l - d\left(\rho\right)$} (N1);
        \draw[sectionblue] ([shift={(0.4,0.7)}]N0.center) edge[bend left,-{Latex[length=2mm,width=2mm]}] ([shift={(-0.4,0.7)}]N1.center);
      \end{tikzpicture}
      \\ or & or \\
      \begin{tikzpicture}
        \node[circle,fill=black,inner sep=0.7mm] (N0) at (0, 0) {};
        \node[circle,fill=black,inner sep=0.7mm] (N1) at (3, 0) {};
        \node[circle,fill=black,inner sep=0.7mm] (Nm1) at (-3, 0) {};
        \stickman{(3.1,0.8)}{tumOrange}
        \stickman{(-0.1,0.8)}{sectionblue}
        \node[below=1mm] at (N0) {$0$};
        \node[below=1mm] at (Nm1) {$-1$};
        \node[below=1mm] at (N1) {$1$};
        \draw (Nm1) -- node[above] {$a + d\left(\rho\right)$} (N0) -- node[above] {$b + 2l - d\left(\rho\right) - 1$} (N1);
        \draw[sectionblue] ([shift={(-0.4,0.7)}]N0.center) edge[bend right,-{Latex[length=2mm,width=2mm]}] ([shift={(0.4,0.7)}]Nm1.center);
      \end{tikzpicture}
      & \begin{tikzpicture}
        \node[circle,fill=black,inner sep=0.7mm] (N0) at (0, 0) {};
        \node[circle,fill=black,inner sep=0.7mm] (N1) at (3, 0) {};
        \node[circle,fill=black,inner sep=0.7mm] (Nm1) at (-3, 0) {};
        \stickman{(3.1,0.8)}{tumOrange}
        \stickman{(-0.1,0.8)}{sectionblue}
        \node[below=1mm] at (N0) {$0$};
        \node[below=1mm] at (Nm1) {$-1$};
        \node[below=1mm] at (N1) {$1$};
        \draw (Nm1) -- node[above] {$a + d\left(\rho\right)$} (N0) -- node[above] {$b + 2l - d\left(\rho\right) - 1$} (N1);
        \draw[sectionblue] ([shift={(0.4,0.7)}]N0.center) edge[bend left,-{Latex[length=2mm,width=2mm]}] ([shift={(-0.4,0.7)}]N1.center);
      \end{tikzpicture}
    \end{tabular}
    \caption{Possible walker locations after $2l - 1$ steps of a path $\rho$ of length $2l$}
    \label{fig:new_paths}
  \end{center}
\end{figure}
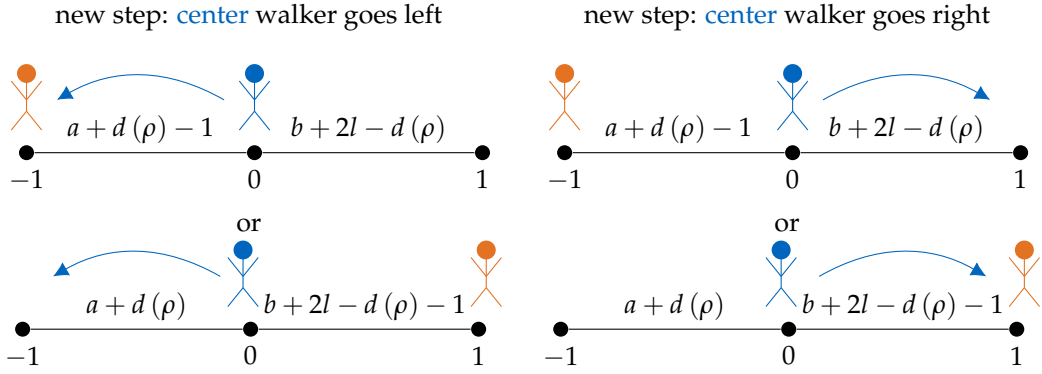

Consider $\rho \in \textrm{Path}_{2l}$. If the walkers move according to $\rho$, then, at any uneven time step, there
will be one walker which is in the \textcolor{sectionblue}{center} and one which is in one of the \textcolor{tumOrange}{outer} nodes
(this corresponds to $s_{\textrm{mixed}}$ in the proof of \autoref{lem:two_pl_urn_rand_meeting_time}). In addition,
at any even time step, except for the beginning and end, both walkers have to be in
the outer nodes, not necessarily the same one (state $s_{\textrm{none}}$). The possible walker locations and edge weights after $2l - 1$
steps of the path $\rho$ are depicted in \autoref{fig:new_paths}, where we assume that we start with edge weights $a$ and $b$. They will be relevant for the recursive
path construction which we present now.

For $\rho \in \textrm{Path}_{2l}$, we define $p_k\left(\rho\right) \in \textrm{Path}_{2\left(l + 1\right)}, 1 \leq k \leq 4$ as follows:
\begin{itemize}
  \item delete the last symbol from $\rho$,
  \item add three new symbols to $\rho$, using the following rules. After the $2l - 1$ steps of $\rho$ which we did not delete, the
  walkers are in two possible configurations: there is one walker which is in the \textcolor{sectionblue}{center}, and the other walker
  can be in either of the \textcolor{tumOrange}{outer} nodes. We now add a symbol to the path such that the \textcolor{sectionblue}{center}
  walker moves either left or right in the next step, and then we add two more symbols such that both walkers return to the center in either
  order. This leaves us a total of $2 \times 2 = 4$ choices for our new path ($2$ choices for which way the \textcolor{sectionblue}{center}
  walker moves, and $2$ choices for the order in which the walkers return). We can therefore construct $4$ new paths, which we will label as follows:
  \begin{itemize}
    \item $p_1\left(\rho\right)$ is the new path where the \textcolor{sectionblue}{center} walker moves left, and then returns immediately to the center before the other walker
    \item $p_2\left(\rho\right)$ is the new path where the \textcolor{sectionblue}{center} walker moves left, and then returns to the center only after the other walker
    \item $p_3\left(\rho\right)$ is the new path where the \textcolor{sectionblue}{center} walker moves right, and then returns immediately to the center before the other walker
    \item $p_4\left(\rho\right)$ is the new path where the \textcolor{sectionblue}{center} walker moves right, and then returns to the center only after the other walker
  \end{itemize}
  This construction is summarized in \autoref{fig:new_paths}.
\end{itemize}

\begin{lemma}
  \label{lem:rec_path_union}
  We can write the set $\textrm{Path}_{2\left(l + 1\right)}$ as a \emph{disjoint} union over the modifications of the paths $\rho$ in $\textrm{Path}_{2l}$:
  \begin{align*}
    \textrm{Path}_{2\left(l + 1\right)} &= \bigcupdot_{\rho \in \textrm{Path}_{2l}} \left\{p_1\left(\rho\right), p_2\left(\rho\right), p_3\left(\rho\right), p_4\left(\rho\right)\right\} \qedhere
  \end{align*}
\end{lemma}

\begin{tproof}
  To see that we get indeed all paths of length $2\left(l + 1\right)$ if we take the given union, note that
  the paths in $\textrm{Path}_{2l}$ must fulfil the condition that at any even time step (except for the beginning and end),
  both walkers have to be in the outer nodes, not necessarily the same one, and that the walkers are in the center
  at the beginning and at the end of the path. This is the only condition, so we will get all possible prefixes of
  paths of length $2\left(l + 1\right)$ by taking all $\rho \in \textrm{Path}_{2l}$ and deleting the last symbol. We then
  only have to complete the paths by all possible suffixes, but the possible suffixes are just the $4$ given in the
  construction above by the constraints imposed on the paths. The union is disjoint because we do not construct any path twice: any two
  different paths in $\textrm{Path}_{2l}$ must already differ somewhere in the first $2l - 1$ steps.
\end{tproof}

The edge weights at time $2l - 1$ depend on where the \textcolor{tumOrange}{outer} walker is at step $2l - 1$. Therefore, we will also
use the notation $\textrm{left}\left(\rho\right) = 1$ if the \textcolor{tumOrange}{outer} walker
is on the left after step $2l - 1$, and $\textrm{left}\left(\rho\right) = 0$ otherwise. Finally, in order to state \autoref{lem:rec_path_construc},
we also need to define $f_{\rho,a,b}$ which is the left edge weight divided by the sum of the edge weights, where we plug in the
edge weights which we get if we start with weights $a$ and $b$ and then execute the path $\rho$.

\begin{lemma}
  \label{lem:rec_path_construc}
  The path probabilities and edge weight ratios of the paths $p_k\left(\rho\right)$ change as follows compared to $\rho$:
  \begin{itemize}
    \item if $\textrm{left}\left(\rho\right) = 1$, then:
    \begin{align*}
      \bbP_{p_1\left(\rho\right),a,b} = \bbP_{p_2\left(\rho\right),a,b} &= \bbP_{\rho,a,b} \cdot \frac{1}{4} \cdot \frac{a + d\left(\rho\right) - 1}{a + b + 2l - 1} \\
      f_{p_1\left(\rho\right),a,b} = f_{p_2\left(\rho\right),a,b} &= \frac{a + d\left(\rho\right) + 2}{a + b + 2l + 2} \\
      \bbP_{p_3\left(\rho\right),a,b} = \bbP_{p_4\left(\rho\right),a,b} &= \bbP_{\rho,a,b} \cdot \frac{1}{4} \cdot \frac{b + 2l - d\left(\rho\right)}{a + b + 2l - 1} \\
      f_{p_3\left(\rho\right),a,b} = f_{p_4\left(\rho\right),a,b} &= \frac{a + d\left(\rho\right)}{a + b + 2l + 2}
    \end{align*}
    \item and if $\textrm{left}\left(\rho\right) = 0$, then:
    \begin{align*}
      \bbP_{p_1\left(\rho\right),a,b} = \bbP_{p_2\left(\rho\right),a,b} &= \bbP_{\rho,a,b} \cdot \frac{1}{4} \cdot \frac{a + d\left(\rho\right)}{a + b + 2l - 1} \\
      f_{p_1\left(\rho\right),a,b} = f_{p_2\left(\rho\right),a,b} &= \frac{a + d\left(\rho\right) + 2}{a + b + 2l + 2} \\
      \bbP_{p_3\left(\rho\right),a,b} = \bbP_{p_4\left(\rho\right),a,b} &= \bbP_{\rho,a,b} \cdot \frac{1}{4} \cdot \frac{b + 2l - d\left(\rho\right) - 1}{a + b + 2l - 1} \\
      f_{p_3\left(\rho\right),a,b} = f_{p_4\left(\rho\right),a,b} &= \frac{a + d\left(\rho\right)}{a + b + 2l + 2} \qedhere
    \end{align*}
  \end{itemize}
\end{lemma}

\begin{tproof}
  We look at how following the path $p_k\left(\rho\right)$ changes the outcome (i.e.~the final edge weights) and the probability of the path,
  compared to the original path $\rho \in \textrm{Path}_{2l}$.
  \begin{itemize}
    \item If the \textcolor{tumOrange}{outer} walker is on the left after $2l - 1$ steps (this is a condition on $\rho$: $\textrm{left}\left(\rho\right) = 1$):
    \begin{itemize}
      \item If the \textcolor{sectionblue}{center} walker should move left in step $2l$ (this is the choice for creating the new paths $p_1\left(\rho\right), p_2\left(\rho\right)$):
      the probability $\bbP_{\rho,a,b}$ of $\rho$ is a product over the probability that the walker indicated in $\rho$ is chosen
      at the respective step (which is always $\frac{1}{2}$) and the probability that the walker moves
      in the direction indicated by $\rho$ (this can either be $1$, if the walker moves back to the center,
      or a fraction depending on the edge weights). In our new modified paths $p_1\left(\rho\right), p_2\left(\rho\right)$ of length $2\left(l + 1\right)$,
      the first difference is that we choose a different walker to move in step $2l$. This event has
      probability $\frac{1}{2}$, but this is the same as chosing the original walker, hence this part
      is already included in the product giving the probability $\bbP_{\rho,a,b}$ of $\rho$.

      Next, the probability for the center walker to go left, if he is chosen to move at step $2l$,
      is given by $\frac{a + d\left(\rho\right) - 1}{a + b + 2l - 1}$ (compare with \autoref{fig:new_paths}).
      This factor is new and has to be added to the product. In the next step, one of the two walkers
      is chosen and will move back to the center. We get a new factor of $\frac{1}{2}$, since the probability
      that the walker indicated by $p_1\left(\rho\right)$ (or $p_2\left(\rho\right)$ respectively) is chosen is $\frac{1}{2}$.
      Finally, in step $2\left(l + 1\right)$, the walker which is still in an
      outer node has to be chosen to move back into the center. This happens with probability $\frac{1}{2}$,
      and this is the final factor to be added to the product. Our new paths therefore have probability
      \begin{align*}
        \bbP_{p_1\left(\rho\right),a,b} = \bbP_{p_2\left(\rho\right),a,b} &= \bbP_{\rho,a,b} \cdot \frac{1}{2} \cdot \frac{1}{2} \cdot \frac{a + d\left(\rho\right) - 1}{a + b + 2l - 1}
      \end{align*}
      The new paths will end with the following ratio of the left edge weight divided by the total edge weights:
      \begin{align*}
        f_{p_1\left(\rho\right),a,b} = f_{p_2\left(\rho\right),a,b} = \frac{a + d\left(\rho\right) + 2}{a + b + 2l + 2}
      \end{align*}
      since the left edge will be traversed twice more by the center walker.
      \item If the \textcolor{sectionblue}{center} walker should move right in step $2l$ (paths $p_3\left(\rho\right), p_4\left(\rho\right)$):
      \begin{align*}
        \textrm{new probability: }\bbP_{\rho,a,b} \cdot \frac{1}{4} \cdot \frac{b + 2l - d\left(\rho\right)}{a + b + 2l - 1},
        \quad \textrm{new outcome: }\frac{a + d\left(\rho\right)}{a + b + 2l + 2}
      \end{align*}
    \end{itemize}
    \item If the \textcolor{tumOrange}{outer} walker is on the right after $2l - 1$ steps ($\textrm{left}\left(\rho\right) = 0$):
    \begin{itemize}
      \item If the \textcolor{sectionblue}{center} walker should move left in step $2l$ (paths $p_1\left(\rho\right), p_2\left(\rho\right)$):
      \begin{align*}
        \textrm{new probability: }\bbP_{\rho,a,b} \cdot \frac{1}{4} \cdot \frac{a + d\left(\rho\right)}{a + b + 2l - 1},
        \quad \textrm{new outcome: }\frac{a + d\left(\rho\right) + 2}{a + b + 2l + 2}
      \end{align*}
      \item If the \textcolor{sectionblue}{center} walker should move right in step $2l$ (paths $p_3\left(\rho\right), p_4\left(\rho\right)$):
      \begin{align*}
        \textrm{new probability: }\bbP_{\rho,a,b} \cdot \frac{1}{4} \cdot \frac{b + 2l - d\left(\rho\right) - 1}{a + b + 2l - 1},
        \quad &\textrm{new outcome: }\frac{a + d\left(\rho\right)}{a + b + 2l + 2} \qedhere
      \end{align*}
    \end{itemize}
  \end{itemize}
\end{tproof}

\begin{oproof}[\autoref{lem:two_pl_urn_rand_recursion}]
  By \autoref{lem:rec_path_union} and \autoref{lem:rec_path_construc}, we have
  \begin{align*}
    &\bbE_{a,b,l+1} = \sum_{\rho \in \textrm{Path}_{2\left(l+1\right)}} \bbP_{\rho,a,b} \cdot \frac{a + d\left(\rho\right)}{a + b + 2l + 2}
    = \sum_{\rho \in \textrm{Path}_{2l}} \sum_{k=1}^4 \bbP_{p_k\left(\rho\right),a,b} \cdot f_{p_k\left(\rho\right),a,b} =\\
    &\sum_{\rho \in \textrm{Path}_{2l}} \bbP_{\rho,a,b} \cdot \frac{1}{2} \cdot
    \left(
      \textrm{left}\left(\rho\right) \cdot \left(
        \frac{a + d\left(\rho\right) - 1}{a + b + 2l - 1} \cdot \frac{a + d\left(\rho\right) + 2}{a + b + 2l + 2}
        + \frac{b + 2l - d\left(\rho\right)}{a + b + 2l - 1} \cdot \frac{a + d\left(\rho\right)}{a + b + 2l + 2}
      \right) \right. \\
    & \qquad\;\left. +
      \left(1 - \textrm{left}\left(\rho\right)\right) \cdot \left(
        \frac{a + d\left(\rho\right)}{a + b + 2l - 1} \cdot \frac{a + d\left(\rho\right) + 2}{a + b + 2l + 2}
        + \frac{b + 2l - d\left(\rho\right) - 1}{a + b + 2l - 1} \cdot \frac{a + d\left(\rho\right)}{a + b + 2l + 2}
      \right)
    \right)
  \end{align*}
  Using this expression for $\bbE_{a,b,l+1}$, we can calculate (a simple, but longer calculation which
  we skip here, see
  \href{https://www.wolframalpha.com/input?i=%28a%2Bd-1%29%28a%2Bd%2B2%29%28a%2Bb%2B2l%29+%2B+%28b%2B2l-d%29%28a%2Bd%29%28a%2Bb%2B2l%29+-+%28a%2Bd%29%28a%2Bb%2B2l-1%29%28a%2Bb%2B2l%2B2%29}{https://bit.ly/3BlMorx}
  and
  \href{https://www.wolframalpha.com/input?i=%28a%2Bd%29%28a%2Bd%2B2%29%28a%2Bb%2B2l%29+%2B+%28b%2B2l-d-1%29%28a%2Bd%29%28a%2Bb%2B2l%29+-+%28a%2Bd%29%28a%2Bb%2B2l-1%29%28a%2Bb%2B2l%2B2%29}{http://bit.ly/3PgLlPq} --
  the calculation is done by expanding the involved fractions to the same denominator):
  \begin{align*}
    \bbE_{a,b,l+1} - \frac{1}{2} \bbE_{a,b,l} &= \bbE_{a,b,l+1} - \sum_{\rho \in \textrm{Path}_{2l}} \frac{1}{2} \cdot \bbP_{\rho,a,b} \cdot \frac{a + d\left(\rho\right)}{a + b + 2l}\\
    &= \sum_{\rho \in \textrm{Path}_{2l}} \bbP_{\rho,a,b} \cdot
    \frac{1}{\left(a + b + 2l - 1\right)\left(a + b + 2l + 2\right)} \;\cdot \\
    &\hphantom{=\qquad} 
    \left(
      \textrm{left}\left(\rho\right) \cdot \frac{d\left(\rho\right) - b - 2l}{a + b + 2l}
      +  \left(1 - \textrm{left}\left(\rho\right)\right) \cdot \frac{a + d\left(\rho\right)}{a + b + 2l}
    \right) \\
    &= \frac{1}{\left(a + b + 2l - 1\right)\left(a + b + 2l + 2\right)} \;\cdot \\
    &\hphantom{=\qquad} \sum_{\rho \in \textrm{Path}_{2l}} \bbP_{\rho,a,b} \cdot
    \left(
      \textrm{left}\left(\rho\right) \cdot \frac{- a - b - 2l}{a + b + 2l}
      +  \frac{a + d\left(\rho\right)}{a + b + 2l}
    \right) \\
    &= \frac{1}{\left(a + b + 2l - 1\right)\left(a + b + 2l + 2\right)} \cdot \left(
      \bbE_{a,b,l} - \sum_{\rho \in \textrm{Path}_{2l}} \bbP_{\rho,a,b} \cdot \textrm{left}\left(\rho\right)
    \right) \\
    &= \frac{1}{\left(a + b + 2l - 1\right)\left(a + b + 2l + 2\right)} \cdot \left(
      \bbE_{a,b,l} - q_{a,b,l}
    \right)
  \end{align*}
  This proves the first equation in \autoref{lem:two_pl_urn_rand_recursion}. For the second equation, we use the same strategy. For our newly
  constructed paths, we already know if $L_n$ occurs or not (whether the last walker to return to the
  center comes from the left node):
  \begin{itemize}
    \item if $\textrm{left}\left(\rho\right) = 1$, then: $L_n$ occurs for $p_1\left(\rho\right), p_2\left(\rho\right), p_4\left(\rho\right)$,
    and does not occur for $p_3\left(\rho\right)$
    \item and if $\textrm{left}\left(\rho\right) = 0$, then: $L_n$ occurs for $p_2\left(\rho\right)$,
    and does not occur for $p_1\left(\rho\right), p_3\left(\rho\right), p_4\left(\rho\right)$.
  \end{itemize}
  Therefore, again by \autoref{lem:rec_path_union} and \autoref{lem:rec_path_construc},
  \begin{align*}
    q_{a,b,l+1} - \frac{1}{2} q_{a,b,l} &=
    \sum_{\rho \in \textrm{Path}_{2\left(l + 1\right)}} \bbP_{\rho,a,b} \cdot \textrm{left}\left(\rho\right)
    - \frac{1}{2} \cdot \sum_{\rho \in \textrm{Path}_{2l}} \bbP_{\rho,a,b} \cdot \textrm{left}\left(\rho\right) \\
    &= \sum_{\rho \in \textrm{Path}_{2l}} \left(\bbP_{p_1\left(\rho\right),a,b} + \bbP_{p_2\left(\rho\right),a,b} + \bbP_{p_4\left(\rho\right),a,b}\right) \cdot \textrm{left}\left(\rho\right) \\
    &\hphantom{=\qquad}+ \bbP_{p_2\left(\rho\right),a,b} \cdot \left(1 - \textrm{left}\left(\rho\right)\right)
    - \frac{1}{2} \cdot \bbP_{\rho,a,b} \cdot \textrm{left}\left(\rho\right) \\
    &=\sum_{\rho \in \textrm{Path}_{2l}} \bbP_{\rho,a,b} \cdot \frac{1}{2} \;\cdot \\
    &\hphantom{=\qquad}\left(
      \frac{a + d\left(\rho\right) - 1}{a + b + 2l - 1} \cdot \textrm{left}\left(\rho\right)
      + \frac{1}{2} \cdot \frac{b + 2l - d\left(\rho\right)}{a + b + 2l - 1} \cdot \textrm{left}\left(\rho\right) \right.\\
    &\hphantom{=\qquad(}\left.  + \;\frac{1}{2} \cdot \frac{a + d\left(\rho\right)}{a + b + 2l - 1} \cdot \left(1 - \textrm{left}\left(\rho\right)\right)
      - \textrm{left}\left(\rho\right)
    \right)\\
    &= \frac{1}{2} \cdot \sum_{\rho \in \textrm{Path}_{2l}} \bbP_{\rho,a,b} \;\cdot \\
    &\hphantom{=\qquad}\left(
      -\frac{1}{2} \cdot \frac{a + b + 2l}{a + b + 2l - 1} \cdot \textrm{left}\left(\rho\right)
      + \frac{1}{2} \cdot \frac{a + d\left(\rho\right)}{a + b + 2l - 1}
    \right) \\
    &= \frac{1}{4} \cdot \frac{a + b + 2l}{a + b + 2l - 1} \cdot \sum_{\rho \in \textrm{Path}_{2l}} \bbP_{\rho,a,b} \cdot
    \left(
      \frac{a + d\left(\rho\right)}{a + b + 2l} - \textrm{left}\left(\rho\right)
    \right) \\
    &= \frac{1}{4} \cdot \frac{a + b + 2l}{a + b + 2l - 1} \cdot \left(\bbE_{a,b,l} - q_{a,b,l}\right) \qedhere
  \end{align*}
\end{oproof}

This proof of \autoref{lem:two_pl_urn_rand_recursion} completes the proof of \autoref{thm:randwalkersmain}.
As a next step, one would analyze the (random) limit
of the fraction $\frac{w\left(n,-1\right)}{w\left(n,-1\right) + w\left(n,0\right)}$ and compare it to the single walker case where the
corresponding limit is Beta-distributed. We believe that the following holds.

\begin{conjecture}
  \label{conj:limit_rand_atoms}
  \autoref{thm:main-alternating} \ref{thm:altwalkersclt} and \ref{thm:altwalkersatoms} (the CLT and the no atoms property) are also true for the case where we randomly select the walker which moves.
\end{conjecture}

% We were able to prove some of the necessary lemmata, but showing that we also have
% \begin{align*}
%   \Ex{n^2\left(M_{n-1} - M_n\right)^2\midd|\calF_{\tau_{n-1}}} \to \frac{1}{2}M_\infty\left(1 - M_\infty\right)
% \end{align*}
% (which was done in the proof of \autoref{thm:main-alternating} \ref{thm:altwalkersclt} for alternating walkers)
% seems to be much more complicated for random walker selection.

\begin{conjecture}
  \label{conj:two_pl_urn_limit_rand}
  $M^{\rm rd}_\infty \in \left[0, 1\right]$ (the limit of the left edge weight proportion, with random walker selection) has a density w.r.t.~the Lebesgue measure on $\left[0, 1\right]$.
  $M^{\rm rd}_{\infty}$ is not Beta-distributed, and its distribution is different from the distribution of the limit
  in \autoref{conj:two_pl_urn_limit} (where we considered alternating walkers).
\end{conjecture}

\ifthenelse{\boolean{extendedversion}}{
\subsubsection{Simulations}

In order to support \autoref{conj:two_pl_urn_limit} and \autoref{conj:two_pl_urn_limit_rand}, we present some
simulations. First, for comparison, the case with a single walker was simulated (i.e.~the standard LERRW on the
three-node segment). It is known that the fraction of the left edge weight converges to a random limit with
distribution $B\left(\frac{1}{2},\frac{1}{2}\right)$ in this case. \autoref{fig:errw_segment_onelimit} shows the CDF of this
Beta distribution, compared to the appropriately normalized observed frequencies of the left edge weight fraction
after $4003$ steps of the walker. $1000001$ samples were collected and placed into $49$ bins of equal size.
The observed frequencies fit the theoretical limit quite well.

\autoref{fig:errw_segment_twolimit_alt} and \autoref{fig:errw_segment_twolimit} show the same simulations for the
case with $2$ walkers, first when they move alternately and then with random player selection. Using the method of
moments, the parameters of an assumed Beta distribution were then estimated from the sampled data, and the corresponding
CDFs were then plotted. The oservations do not swerve far from the Beta distributions, but compared to \autoref{fig:errw_segment_onelimit},
one can see that the observations in \autoref{fig:errw_segment_twolimit_alt} and \autoref{fig:errw_segment_twolimit} are
noticeably below the Beta CDF in the range from $0.05$ to $0.25$, while this was not the case for \autoref{fig:errw_segment_onelimit}.
In \autoref{fig:errw_segment_onelimit_detail}, \autoref{fig:errw_segment_twolimit_alt_detail} and \autoref{fig:errw_segment_twolimit_detail},
this range is shown in detail for easier comparison.

Note that in the case where the two walkers move alternately, a much more precise calculation can be made since the actual
distribution of the left edge weight fraction can be calculated with much less computational effort at a given time step.
\autoref{fig:errw_segment_twolimit_alt} and \autoref{fig:errw_segment_twolimit_alt_detail} where therefore generated using
a different method which exploits that the two walkers are back in the center every four steps. This method could also be
applied to the case with a single walker, but this was not done here in order to provide a better comparison to the 2 walker
case with random walker selection.

\begin{figure}[H]
  \begin{center}
    \begin{tikzpicture}
  \begin{axis}[
    xlabel={Limit of normalized left edge weight},
    ylabel={Frequency},
    width=12cm,
    grid=major,
    legend pos=north west
  ]
    \addplot[color=tumBlue] coordinates {
      (0.0102, 4.4834) (0.0306, 1.8757) (0.0510, 1.4560) (0.0714, 1.2231) (0.0918, 1.1025) (0.1122, 1.0183) (0.1327, 0.9348) (0.1531, 0.8812) (0.1735, 0.8419) (0.1939, 0.7988) (0.2143, 0.7816) (0.2347, 0.7408) (0.2551, 0.7250) (0.2755, 0.7212) (0.2959, 0.7007) (0.3163, 0.6775) (0.3367, 0.6747) (0.3571, 0.6704) (0.3776, 0.6662) (0.3980, 0.6462) (0.4184, 0.6404) (0.4388, 0.6369) (0.4592, 0.6321) (0.4796, 0.6432) (0.5000, 0.6383) (0.5204, 0.6375) (0.5408, 0.6447) (0.5612, 0.6430) (0.5816, 0.6485) (0.6020, 0.6539) (0.6224, 0.6428) (0.6429, 0.6656) (0.6633, 0.6774) (0.6837, 0.6916) (0.7041, 0.7042) (0.7245, 0.7131) (0.7449, 0.7353) (0.7653, 0.7549) (0.7857, 0.7701) (0.8061, 0.8014) (0.8265, 0.8442) (0.8469, 0.8775) (0.8673, 0.9385) (0.8878, 1.0120) (0.9082, 1.1153) (0.9286, 1.2054) (0.9490, 1.4548) (0.9694, 1.8762) (0.9898, 4.4791)
    };
    \addlegendentry{$1000001$ simulations at $4003$ steps};
    \addplot[color=tumOrange] coordinates {
      (0.0102, 3.1672) (0.0306, 1.8478) (0.0510, 1.4466) (0.0714, 1.2360) (0.0918, 1.1022) (0.1122, 1.0084) (0.1327, 0.9384) (0.1531, 0.8841) (0.1735, 0.8406) (0.1939, 0.8052) (0.2143, 0.7757) (0.2347, 0.7511) (0.2551, 0.7302) (0.2755, 0.7125) (0.2959, 0.6974) (0.3163, 0.6845) (0.3367, 0.6735) (0.3571, 0.6643) (0.3776, 0.6566) (0.3980, 0.6503) (0.4184, 0.6453) (0.4388, 0.6414) (0.4592, 0.6388) (0.4796, 0.6372) (0.5000, 0.6366) (0.5204, 0.6372) (0.5408, 0.6388) (0.5612, 0.6414) (0.5816, 0.6453) (0.6020, 0.6503) (0.6224, 0.6566) (0.6429, 0.6643) (0.6633, 0.6735) (0.6837, 0.6845) (0.7041, 0.6974) (0.7245, 0.7125) (0.7449, 0.7302) (0.7653, 0.7511) (0.7857, 0.7757) (0.8061, 0.8052) (0.8265, 0.8406) (0.8469, 0.8841) (0.8673, 0.9384) (0.8878, 1.0084) (0.9082, 1.1022) (0.9286, 1.2360) (0.9490, 1.4466) (0.9694, 1.8478) (0.9898, 3.1672)
    };
    \addlegendentry{CDF of $B\left(\frac{1}{2},\frac{1}{2}\right)$};
  \end{axis}
\end{tikzpicture}

    \caption{Simulated distribution of the limit in the $1$ walker case compared to actual limit}
    \label{fig:errw_segment_onelimit}
  \end{center}
\end{figure}

\begin{figure}[H]
  \begin{center}
    \begin{tikzpicture}
  \begin{axis}[
    xlabel={Limit of normalized left edge weight},
    ylabel={Frequency},
    width=12cm,
    grid=major,
    legend pos=north west
  ]
    \addplot[color=tumGreen] coordinates {
      (0.0102, 3.8145) (0.0306, 1.6672) (0.0510, 1.3240) (0.0714, 1.1367) (0.0918, 1.0568) (0.1122, 0.9922) (0.1327, 0.9498) (0.1531, 0.8857) (0.1735, 0.8623) (0.1939, 0.8391) (0.2143, 0.8155) (0.2347, 0.7860) (0.2551, 0.7962) (0.2755, 0.7731) (0.2959, 0.7634) (0.3163, 0.7505) (0.3367, 0.7413) (0.3571, 0.7493) (0.3776, 0.7397) (0.3980, 0.7297) (0.4184, 0.7327) (0.4388, 0.7324) (0.4592, 0.7203) (0.4796, 0.7312) (0.5000, 0.7314) (0.5204, 0.7264) (0.5408, 0.7308) (0.5612, 0.7352) (0.5816, 0.7380) (0.6020, 0.7397) (0.6224, 0.7302) (0.6429, 0.7506) (0.6633, 0.7603) (0.6837, 0.7581) (0.7041, 0.7601) (0.7245, 0.7803) (0.7449, 0.7921) (0.7653, 0.8069) (0.7857, 0.8084) (0.8061, 0.8386) (0.8265, 0.8639) (0.8469, 0.8929) (0.8673, 0.9441) (0.8878, 0.9904) (0.9082, 1.0703) (0.9286, 1.1381) (0.9490, 1.3178) (0.9694, 1.6707) (0.9898, 3.8330)
    };
    \addlegendentry{$1000001$ simulations at $4003$ steps};
    \addplot[color=tumBlue] coordinates {
      (0.0102, 3.8292) (0.0306, 1.6733) (0.0510, 1.3232) (0.0714, 1.1611) (0.0918, 1.0500) (0.1122, 0.9791) (0.1327, 0.9331) (0.1531, 0.8883) (0.1735, 0.8629) (0.1939, 0.8333) (0.2143, 0.8185) (0.2347, 0.7973) (0.2551, 0.7839) (0.2755, 0.7774) (0.2959, 0.7633) (0.3163, 0.7601) (0.3367, 0.7490) (0.3571, 0.7436) (0.3776, 0.7437) (0.3980, 0.7355) (0.4184, 0.7372) (0.4388, 0.7305) (0.4592, 0.7335) (0.4796, 0.7281) (0.5000, 0.7278) (0.5204, 0.7326) (0.5408, 0.7290) (0.5612, 0.7350) (0.5816, 0.7327) (0.6020, 0.7355) (0.6224, 0.7437) (0.6429, 0.7436) (0.6633, 0.7536) (0.6837, 0.7555) (0.7041, 0.7680) (0.7245, 0.7727) (0.7449, 0.7839) (0.7653, 0.8022) (0.7857, 0.8135) (0.8061, 0.8385) (0.8265, 0.8578) (0.8469, 0.8883) (0.8673, 0.9331) (0.8878, 0.9791) (0.9082, 1.0567) (0.9286, 1.1544) (0.9490, 1.3232) (0.9694, 1.6733) (0.9898, 3.8292)
    };
    \addlegendentry{$1$ precise calculation at $16012$ steps};
    \addplot[color=tumOrange] coordinates {
      (0.0102, 2.7152) (0.0306, 1.7256) (0.0510, 1.4045) (0.0714, 1.2303) (0.0918, 1.1173) (0.1122, 1.0367) (0.1327, 0.9759) (0.1531, 0.9281) (0.1735, 0.8896) (0.1939, 0.8579) (0.2143, 0.8315) (0.2347, 0.8092) (0.2551, 0.7902) (0.2755, 0.7740) (0.2959, 0.7602) (0.3163, 0.7484) (0.3367, 0.7383) (0.3571, 0.7298) (0.3776, 0.7227) (0.3980, 0.7168) (0.4184, 0.7122) (0.4388, 0.7086) (0.4592, 0.7061) (0.4796, 0.7046) (0.5000, 0.7041) (0.5204, 0.7046) (0.5408, 0.7061) (0.5612, 0.7086) (0.5816, 0.7122) (0.6020, 0.7168) (0.6224, 0.7227) (0.6429, 0.7298) (0.6633, 0.7383) (0.6837, 0.7484) (0.7041, 0.7602) (0.7245, 0.7740) (0.7449, 0.7902) (0.7653, 0.8092) (0.7857, 0.8315) (0.8061, 0.8579) (0.8265, 0.8896) (0.8469, 0.9281) (0.8673, 0.9759) (0.8878, 1.0367) (0.9082, 1.1173) (0.9286, 1.2303) (0.9490, 1.4045) (0.9694, 1.7256) (0.9898, 2.7152)
    };
    \addlegendentry{CDF of $B\left(0.579,0.579\right)$};
  \end{axis}
\end{tikzpicture}

    \caption{Calculated distribution of the limit in the $2$ walkers, moving alternately case compared to Beta distribution}
    \label{fig:errw_segment_twolimit_alt}
  \end{center}
\end{figure}

\begin{figure}[H]
  \begin{center}
    \begin{tikzpicture}
  \begin{axis}[
    xlabel={Limit of normalized left edge weight},
    ylabel={Frequency},
    width=12cm,
    grid=major,
    legend pos=north west
  ]
    \addplot[color=tumBlue] coordinates {
      (0.0102, 4.0041) (0.0306, 1.7524) (0.0510, 1.3815) (0.0714, 1.1668) (0.0918, 1.0764) (0.1122, 0.9802) (0.1327, 0.9402) (0.1531, 0.8880) (0.1735, 0.8619) (0.1939, 0.8412) (0.2143, 0.8159) (0.2347, 0.7834) (0.2551, 0.7661) (0.2755, 0.7624) (0.2959, 0.7398) (0.3163, 0.7351) (0.3367, 0.7364) (0.3571, 0.7231) (0.3776, 0.7169) (0.3980, 0.6977) (0.4184, 0.6959) (0.4388, 0.7058) (0.4592, 0.6913) (0.4796, 0.6889) (0.5000, 0.6934) (0.5204, 0.7043) (0.5408, 0.6862) (0.5612, 0.7100) (0.5816, 0.6891) (0.6020, 0.7101) (0.6224, 0.7179) (0.6429, 0.7176) (0.6633, 0.7399) (0.6837, 0.7383) (0.7041, 0.7361) (0.7245, 0.7710) (0.7449, 0.7623) (0.7653, 0.7898) (0.7857, 0.8042) (0.8061, 0.8292) (0.8265, 0.8553) (0.8469, 0.8881) (0.8673, 0.9444) (0.8878, 1.0119) (0.9082, 1.0841) (0.9286, 1.1488) (0.9490, 1.3686) (0.9694, 1.7354) (0.9898, 4.0106)
    };
    \addlegendentry{$1000001$ simulations at $4003$ steps};
    \addplot[color=tumOrange] coordinates {
      (0.0102, 2.8629) (0.0306, 1.7689) (0.0510, 1.4214) (0.0714, 1.2350) (0.0918, 1.1148) (0.1122, 1.0296) (0.1327, 0.9656) (0.1531, 0.9155) (0.1735, 0.8751) (0.1939, 0.8421) (0.2143, 0.8145) (0.2347, 0.7913) (0.2551, 0.7716) (0.2755, 0.7549) (0.2959, 0.7405) (0.3163, 0.7283) (0.3367, 0.7179) (0.3571, 0.7091) (0.3776, 0.7017) (0.3980, 0.6957) (0.4184, 0.6909) (0.4388, 0.6872) (0.4592, 0.6846) (0.4796, 0.6831) (0.5000, 0.6826) (0.5204, 0.6831) (0.5408, 0.6846) (0.5612, 0.6872) (0.5816, 0.6908) (0.6020, 0.6956) (0.6224, 0.7016) (0.6429, 0.7090) (0.6633, 0.7177) (0.6837, 0.7281) (0.7041, 0.7403) (0.7245, 0.7546) (0.7449, 0.7714) (0.7653, 0.7910) (0.7857, 0.8142) (0.8061, 0.8417) (0.8265, 0.8747) (0.8469, 0.9150) (0.8673, 0.9650) (0.8878, 1.0290) (0.9082, 1.1141) (0.9286, 1.2340) (0.9490, 1.4202) (0.9694, 1.7671) (0.9898, 2.8590)
    };
    \addlegendentry{CDF of $B\left(0.553,0.554\right)$};
  \end{axis}
\end{tikzpicture}

    \caption{Simulated distribution of the limit in the $2$ walkers, random selection case compared to Beta distribution}
    \label{fig:errw_segment_twolimit}
  \end{center}
\end{figure}

\begin{figure}[H]
  \begin{center}
    \begin{tikzpicture}
  \begin{axis}[
    xlabel={Limit of normalized left edge weight},
    ylabel={Frequency},
    width=12cm,
    grid=major,
    legend pos=north west,
    x tick label style={
      /pgf/number format/precision=2,
      /pgf/number format/fixed}
  ]
    \addplot[color=tumBlue] coordinates {
      (0.0510, 1.4560) (0.0714, 1.2231) (0.0918, 1.1025) (0.1122, 1.0183) (0.1327, 0.9348) (0.1531, 0.8812) (0.1735, 0.8419) (0.1939, 0.7988) (0.2143, 0.7816) (0.2347, 0.7408)
    };
    \addlegendentry{$1000001$ simulations at $4003$ steps};
    \addplot[color=tumOrange] coordinates {
      (0.0510, 1.4466) (0.0714, 1.2360) (0.0918, 1.1022) (0.1122, 1.0084) (0.1327, 0.9384) (0.1531, 0.8841) (0.1735, 0.8406) (0.1939, 0.8052) (0.2143, 0.7757) (0.2347, 0.7511)
    };
    \addlegendentry{CDF of $B\left(\frac{1}{2},\frac{1}{2}\right)$};
  \end{axis}
\end{tikzpicture}

    \caption{Simulated distribution of the limit in the $1$ walker case compared to actual limit}
    \label{fig:errw_segment_onelimit_detail}
  \end{center}
\end{figure}

\begin{figure}[H]
  \begin{center}
    \begin{tikzpicture}
  \begin{axis}[
    xlabel={Limit of normalized left edge weight},
    ylabel={Frequency},
    width=12cm,
    grid=major,
    legend pos=north west,
    x tick label style={
      /pgf/number format/precision=2,
      /pgf/number format/fixed}
  ]
    \addplot[color=tumGreen] coordinates {
      (0.0510, 1.3240) (0.0714, 1.1367) (0.0918, 1.0568) (0.1122, 0.9922) (0.1327, 0.9498) (0.1531, 0.8857) (0.1735, 0.8623) (0.1939, 0.8391) (0.2143, 0.8155) (0.2347, 0.7860)
    };
    \addlegendentry{$1000001$ simulations at $4003$ steps};
    \addplot[color=tumBlue] coordinates {
      (0.0510, 1.3232) (0.0714, 1.1611) (0.0918, 1.0500) (0.1122, 0.9791) (0.1327, 0.9331) (0.1531, 0.8883) (0.1735, 0.8629) (0.1939, 0.8333) (0.2143, 0.8185) (0.2347, 0.7973)
    };
    \addlegendentry{$1$ precise calculation at $16012$ steps};
    \addplot[color=tumOrange] coordinates {
      (0.0510, 1.4045) (0.0714, 1.2303) (0.0918, 1.1173) (0.1122, 1.0367) (0.1327, 0.9759) (0.1531, 0.9281) (0.1735, 0.8896) (0.1939, 0.8579) (0.2143, 0.8315) (0.2347, 0.8092)
    };
    \addlegendentry{CDF of $B\left(0.579,0.579\right)$};
  \end{axis}
\end{tikzpicture}

    \caption{Calculated distribution of the limit in the $2$ walkers, moving alternately case compared to Beta distribution}
    \label{fig:errw_segment_twolimit_alt_detail}
  \end{center}
\end{figure}

\begin{figure}[H]
  \begin{center}
    \begin{tikzpicture}
  \begin{axis}[
    xlabel={Limit of normalized left edge weight},
    ylabel={Frequency},
    width=12cm,
    grid=major,
    legend pos=north west,
    x tick label style={
      /pgf/number format/precision=2,
      /pgf/number format/fixed}
  ]
    \addplot[color=tumBlue] coordinates {
      (0.0510, 1.3815) (0.0714, 1.1668) (0.0918, 1.0764) (0.1122, 0.9802) (0.1327, 0.9402) (0.1531, 0.8880) (0.1735, 0.8619) (0.1939, 0.8412) (0.2143, 0.8159) (0.2347, 0.7834)
    };
    \addlegendentry{$1000001$ simulations at $4003$ steps};
    \addplot[color=tumOrange] coordinates {
      (0.0510, 1.4214) (0.0714, 1.2350) (0.0918, 1.1148) (0.1122, 1.0296) (0.1327, 0.9656) (0.1531, 0.9155) (0.1735, 0.8751) (0.1939, 0.8421) (0.2143, 0.8145) (0.2347, 0.7913)
    };
    \addlegendentry{CDF of $B\left(0.553,0.554\right)$};
  \end{axis}
\end{tikzpicture}

    \caption{Simulated distribution of the limit in the $2$ walkers, random selection case compared to Beta distribution}
    \label{fig:errw_segment_twolimit_detail}
  \end{center}
\end{figure}
}{}

\section{A Finite Number of Walkers on $\bbZ$}
\label{ssec:multiple_walkers_z}

So far, we have looked at the linearly edge-reinforced random walk with multiple walkers only on a very simple graph and
only with $2$ walkers. We now consider a finite number $K$ of edge-reinforced random walkers on $\bbZ$. As before,
the transition probabilities depend on the edge weights $w\left(n, z\right) > 0$ for $n \geq 0, z \in \bbZ$
where we will abuse notation and use $z$ to denote the edge from $z$ to $z+1$.
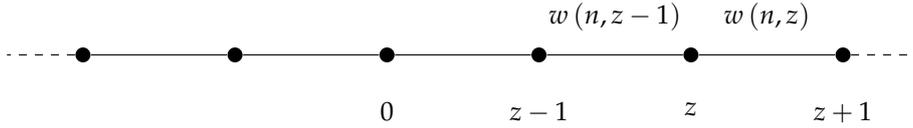
\begin{figure}[H]
  \begin{center}
    \begin{tikzpicture}
      \node[circle,fill=black,inner sep=0.7mm] (N0) at (0, 0) {};
      \node[circle,fill=black,inner sep=0.7mm] (N1) at (2, 0) {};
      \node[circle,fill=black,inner sep=0.7mm] (N2) at (4, 0) {};
      \node[circle,fill=black,inner sep=0.7mm] (N3) at (6, 0) {};
      \node[circle,fill=black,inner sep=0.7mm] (Nm1) at (-2, 0) {};
      \node[circle,fill=black,inner sep=0.7mm] (Nm2) at (-4, 0) {};
      \node[below=5mm] at (N0) {$0$};
      \node[below=5mm] at (N1) {$z-1$};
      \node[below=5mm] at (N2) {$z$};
      \node[below=5mm] at (N3) {$z+1$};
      \draw (Nm2) -- (Nm1) -- (N0) -- (N1) -- node[above=2mm] {$w\left(n,z-1\right)$} (N2) -- node[above=2mm] {$w\left(n,z\right)$} (N3);
      \draw[dashed] (N3) -- (7, 0);
      \draw[dashed] (Nm2) -- (-5, 0);
      \node (aN2l) at (3.9, 0.5) {};
      \node (aN2r) at (4.1, 0.5) {};
    \end{tikzpicture}
    \caption{Edge weights on $\bbZ$ at time $n$}
    \label{fig:errw_z_multiple}
  \end{center}
\end{figure}

\subsection{Model and Results}

We have $K$ sequences (for the $K$ walkers)
$\big(X^{\left(i\right)}_n\big)_{n \geq 0}$ (where $1 \leq i \leq K$) of interacting nearest neighbour processes with the following
dynamics, adapted from \autoref{def:errw}:

If $\calG_n$ denotes $\sigma\left(\big\{ X^{\left(i\right)}_{m} : 0 \leq m \leq n, 1 \leq i \leq K \big\} \cup \left\{w\left(m, z\right) : 0 \leq m \leq n, z \in \bbZ \right\}\right)$
i.e.~the history of the random walkers and edge weights up to and including time $n$, then we define,
conditional on $\calG_n$, the following transition probabilities:
\begin{itemize}
  \item At every time step, the random walker $i$ ($1 \leq i \leq K$) which is going to jump is chosen uniformly at random
  (independently of $\calG_n$) amongst the $K$ walkers.
  \item If the chosen random walker is at position $z$ (i.e.~$X^{\left(i\right)}_n = z$), then he jumps
  \begin{itemize}
    \item to the right (i.e.~$X^{\left(i\right)}_{n + 1} = z + 1$) with probability $\frac{w\left(n, z\right)}{w\left(n, z - 1\right) + w\left(n, z\right)}\,,$
    \item to the left (i.e.~$X^{\left(i\right)}_{n + 1} = z - 1$) with probability $\frac{w\left(n, z - 1\right)}{w\left(n, z - 1\right) + w\left(n, z\right)}\,,$
  \end{itemize}
  i.e.~the jump probabilities are proportional to the corresponding edge weights.
  \item If $z^{\ast}$ is the traversed edge ($z^{\ast} = z$ if the walker jumps to the right, $z^{\ast} = z - 1$
  if he jumps to the left), then for $y \neq z^\ast$, $w\left(n + 1, y\right) = w\left(n, y\right)$
  and $w\left(n + 1, z^{\ast}\right) \geq w\left(n, z^{\ast}\right)$, i.e.~the weight of the
  traversed edge may be increased according to some reinforcement scheme.
  
  We consider schemes where
  the increment $w\left(n + 1, z^{\ast}\right) - w\left(n, z^{\ast}\right)$ only depends on
  $z^{\ast}$, and the number of times the edge was crossed up to time $n$. 
  In other words, $w\left(n, z\right)$
  can still be written in terms of the weight function $W_z\left(k\right)$ as in \autoref{def:errw}.
  \item The initial edge weights can be chosen arbitrarily, but all of them must be strictly positive.
  \item The initial positions of the $K$ walkers can be chosen arbitrarily.
\end{itemize}

We call a walker recurrent if he visits all integers infinitely often, and we say that he has finite range if he only visits finitely many integers. We have the following main result:
\begin{theorem}
  \label{thm:allRecurrentOrAllFiniteRange}
  Assume the edge-reinforced random walk with $K$ walkers starts with an initial configuration of the
  weights $w\left(0, z\right)$ such that all but finitely many of them are $1$.
  Then, we have the following:
  \begin{align*}
    \Prb{\forall i: X^{\left(i\right)} \textrm{ is recurrent}}
    + \Prb{\forall i: X^{\left(i\right)} \textrm{ has finite range}} &= 1\, . \qedhere
  \end{align*}
\end{theorem}

We can further characterize these two possible cases for a particular class of reinforcement schemes, called ``sequence type'', which means that
all initial edge weights are $1$, and that the edge weights are then increased by a fixed sequence of increments $\mathbf{a} = \left(a_k\right)_{k \geq 1}$ which do not
depend on the location of the edge.
\begin{theorem}
  \label{thm:seqtypefinrange}
  Consider the edge-reinforced random walk with $K$ walkers and a reinforcement scheme of sequence type. Set
  \begin{align*}
    \phi\left(\mathbf{a}\right) = \sum_{k = 1}^{\infty} \left(1 + \sum_{l = 1}^k a_l\right)^{-1}\,.
  \end{align*}
  If $\phi\left(\mathbf{a}\right) = \infty$, then all walkers are recurrent a.s., and if $\phi\left(\mathbf{a}\right) < \infty$,
  then all walkers have finite range a.s.
\end{theorem}

\subsection{Recurrence or Finite Range on $\bbZ$}

\begin{definition}
  \label{def:rectrans_multiple}
  For each $i\in \{1, \ldots, K\}$, we say that walker $i$ 
  \begin{itemize}
    \item is \textbf{transient}, if he visits every integer only finitely often, that is, every integer
    appears only finitely often in the sequence $\left(X^{\left(i\right)}_n\right)_{n \geq 0}$
    \item is \textbf{recurrent}, if he visits every integer infinitely often, that is, every integer
    appears infinitely often in the sequence $\big(X^{\left(i\right)}_n\big)_{n \geq 0}$
    \item has \textbf{finite range}, if he only visits finitely many integers, that is, the number of distinct
    integers appearing in the sequence $\left(X^{\left(i\right)}_n\right)_{n \geq 0}$ is finite \qedhere
  \end{itemize}
\end{definition}

Our main result shows that either all walkers are recurrent or all walkers have finite range (\autoref{thm:allRecurrentOrAllFiniteRange}). This was already
known for a single walker. The proof for the single walker case has to be adapted and additional steps are needed to show that all walkers
show the same behavior: we will see that it cannot be the case that one walker has finite range, while another is recurrent. 
Of course, this is very plausible.

The proof of \autoref{thm:allRecurrentOrAllFiniteRange} follows the same strategy as \cite{errwz}, with some changes in the details. More precisely, \autoref{lem:retRootOrFinManyNew}
uses the methods from \cite{errwz} and modifies them to fit the case with multiple walkers. As a result, we see that every walker
either reaches $0$ or visits only finitely many nodes which have not been visited before by any other walker. \autoref{lem:labelExchange}
is new for the case with multiple walkers and shows that any two walkers which meet infinitely often are either both recurrent or do both have finite range.
Combining the lemmas in the proof of \autoref{thm:allRecurrentOrAllFiniteRange} allows us to conclude
that the behavior of all walkers is the same: either all of the walkers are recurrent, or all have finite range.

\begin{lemma}
  \label{lem:retRootOrFinManyNew}
  Assume the edge-reinforced random walk with $K$ walkers starts with an initial configuration of the
  weights $w\left(0, z\right)$ such that all but finitely many of them are $1$.
  Then, for every random walker $i$ ($1 \leq i \leq K$) 
  we have the following:
  \begin{align*}
    &\Prb{X^{\left(i\right)}_n = 0 \textrm{ for some }n \geq 0} \\
    +\; &\bbP\left[ X^{\left(i\right)}_n \neq 0\textrm{ for all }n \geq 0 \textrm{ and } X^{\left(i\right)} \textrm{ only visits finitely many nodes}\right. \\
    &\hphantom{\bbP[}\left. \vphantom{X^{\left(i\right)}_n}\qquad\quad \textrm{which have not been visited before by any other walker} \right] = 1\,. \qedhere
  \end{align*}
\end{lemma}

\begin{lemma}
  \label{lem:labelExchange}
  Assume the edge-reinforced random walk with $K \geq 2$ walkers starts with an arbitrary initial configuration of the
  weights $w\left(0, z\right)$.
  Assume further that $X^{\left(1\right)}$ and $X^{\left(2\right)}$ meet infinitely often. Then, almost surely,
  \begin{enumerate}[(i)]
    \item \label{lem:labelExchangeBothInf} if one of the walkers $X^{\left(1\right)}$ and $X^{\left(2\right)}$ does
    not have finite range, then both $X^{\left(1\right)}$ and $X^{\left(2\right)}$ do not have finite range.
    \item \label{lem:labelExchangeBothVisInf} if some integer $z$ is visited infinitely often by one of the walkers $X^{\left(1\right)}$
    or $X^{\left(2\right)}$, then both $X^{\left(1\right)}$ and $X^{\left(2\right)}$ visit $z$ infinitely often.
    \item \label{lem:labelExchangeBothRec} if every integer is visited infinitely often by one of the walkers $X^{\left(1\right)}$
    or $X^{\left(2\right)}$, then both $X^{\left(1\right)}$ and $X^{\left(2\right)}$ are recurrent. \qedhere
  \end{enumerate}
\end{lemma}

\begin{oproof}[\autoref{thm:allRecurrentOrAllFiniteRange}]
  We have to show the following: if at least one of the walkers does not have finite range, then,
  almost surely, all of them are recurrent. We look at the (random) set
  \begin{align*}
    A_{\textrm{rec}} := \left\{ i : 1 \leq i \leq K \textrm{ and } X^{\left(i\right)} \textrm{ is recurrent} \right\}\,.
  \end{align*}
  The complement can be written as
  \begin{align*}
    A_{\textrm{rec}}^{\textrm{c}} &= \left\{1, \ldots, K\right\} \setminus A_{\textrm{rec}} \\
    &= \bigcup_{z \in \bbZ} \left\{ i : 1 \leq i \leq K \textrm{ and } X^{\left(i\right)} \textrm{ visits } z \textrm{ only finitely often} \right\} \\
    &= \bigcup_{z \in \bbZ} \bigcup_{n \geq 0} \left\{ i : 1 \leq i \leq K \textrm{ and } X^{\left(i\right)} \textrm{ does not visit } z \textrm{ at any time }m \geq n \right\} \,.
  \end{align*}
  Since at time $n$, all but finitely many edge weights are still $1$, we can apply \autoref{lem:retRootOrFinManyNew} to the random walk from time $n$ onwards.
  By relabeling, \autoref{lem:retRootOrFinManyNew} can also be applied to any node $z$ and not just $0$, and with arbitrary initial walker position. Hence, a.s.
  \begin{align*}
    &\left\{ i : 1 \leq i \leq K \textrm{ and } X^{\left(i\right)} \textrm{ does not visit } z \textrm{ at any time }m \geq n \right\} \\
    = &\left\{ i : 1 \leq i \leq K \textrm{, } X^{\left(i\right)} \textrm{ does not visit } z \textrm{ at any time }m \geq n \textrm{ and } X^{\left(i\right)} \textrm{ only visits}\right.\\
    &\hphantom{\{}\left.\vphantom{X^{\left(i\right)}}\textrm{finitely many nodes which have not been visited before by any other walker}\right\}\,.
  \end{align*}
  Since $A_{\textrm{rec}}^{\textrm{c}}$ is a union of sets of this type, we conclude that, a.s.
  \begin{align*}
    A_{\textrm{rec}}^{\textrm{c}} \subseteq &\left\{i : 1 \leq i \leq K \textrm{ and }X^{\left(i\right)} \textrm{ only visits finitely many nodes}\right.\\
    &\hphantom{\{}\left.\vphantom{X^{\left(i\right)}}\textrm{which have not been visited before by any other walker}\right\}\,.
  \end{align*}
  Now, $A_{\textrm{rec}}^{\textrm{c}} = \left\{1, \ldots, K\right\}$ implies that only finitely many nodes are visited overall,
  so the event $\big\{\forall i: X^{\left(i\right)} \textrm{ has finite range}\big\}$ would occur. Hence, it suffices to show that
  $\left\{\forall i: X^{\left(i\right)} \textrm{ is recurrent}\right\}$ occurs whenever $A_{\textrm{rec}}^{\textrm{c}} \neq \left\{1, \ldots, K\right\}$ a.s.

  \ifthenelse{\boolean{extendedversion}}{
  \begin{figure}[H]
    \begin{center}


      \caption{Illustration: proof that either all walkers are recurrent or all have finite range}
      \label{fig:rec_proof}
    \end{center}
    Above, two possible behaviors of the walker $i_2$ (which are both, in fact, almost surely
    impossible) are shown: either $i_2$ has finite range (indicated by the \textcolor{tumGreen}{light green} background),
    but then it would meet $i_1$ infinitely often,
    or $i_2$ does not have finite range, but then $i_3$ would have to ``free the path'' for $i_2$
    (indicated by the \textcolor{tumOrange}{light orange} background) and $i_2$ would meet $i_3$ infinitely often. The relevant
    meeting points are circled.
  \end{figure}
  }{}

  $A_{\textrm{rec}}^{\textrm{c}} \neq \left\{1, \ldots, K\right\}$ implies that there is at least one recurrent walker, say $X^{\left(i_1\right)}$.
  Now, take, for a contradiction, any walker $X^{\left(i_2\right)}$ with $i_2 \in A_{\textrm{rec}}^{\textrm{c}}$. If $i_2$ has finite range, then 
  $i_2$
  will meet the recurrent walker $i_1$ infinitely often (recall that only walker moves in each time step!), but this is a contradiction to \autoref{lem:labelExchange}
  \ref{lem:labelExchangeBothInf}. So $i_2$ must have infinite range. Since $i_2$ is not recurrent, there must be some integer $z$
  which is never visited by $i_2$, w.l.o.g.~we assume that $i_2$ only visits nodes to the right of $z$. The walkers in $A_{\textrm{rec}}^{\textrm{c}}$, including $i_2$, only
  visit finitely many nodes not visited before by any other walker, so at least one walker $i_3$ in $A_{\textrm{rec}}$ must be to
  the right of $i_2$ infinitely often in order to ``free the path'' for $i_2$ (recall that we consider the case where $i_2$ has infinite range
  and must thus go infinitely far to the right). As the walkers in $A_{\textrm{rec}}$
  are all recurrent, and as $i_2$ only visits nodes to the right of $z$, this implies that $i_2$ meets
  this walker $i_3$ infinitely often, which is a contradiction to \autoref{lem:labelExchange} \ref{lem:labelExchangeBothRec}.
  Hence, there is a contradiction if we assume that a walker $i_2 \in A_{\textrm{rec}}^{\textrm{c}}$ exists.\\

  We conclude that $A_{\textrm{rec}}^{\textrm{c}} = \varnothing$ (which implies that all walkers are recurrent) a.s.~if
  $A_{\textrm{rec}}^{\textrm{c}} \neq \left\{1, \ldots, K\right\}$.
\end{oproof}

\begin{oproof}[\autoref{lem:retRootOrFinManyNew}]
  We follow the proof of \cite[Lemma 3.0]{errwz}. Consider a fixed random walker $i$ and assume w.l.o.g. $X^{\left(i\right)}_0 \geq 0$. We now define:
  \begin{align*}
    F\left(n, z\right) &:= \begin{cases}
      \sum_{y=0}^{z-1} \frac{1}{w\left(n, y\right)} & \textrm{ if } z > 0 \\
      0 & \textrm{ if } z \leq 0\,,
    \end{cases}
    \qquad \qquad \tau^{\left(i\right)} := \inf \left\{ n \geq 0: X^{\left(i\right)}_n \leq 0 \right\}\,, \\
    M^{\left(i\right)}_n &:= F\left(n \land \tau^{\left(i\right)}, X^{\left(i\right)}_{n \land \tau^{\left(i\right)}}\right) \\
    H^{\left(i\right)}_n &:= M^{\left(i\right)}_n + \sum_{m=1}^n
    \left(\frac{1}{w\left(m-1,X^{\left(i\right)}_{m-1}\right)} - \frac{1}{w\left(m,X^{\left(i\right)}_{m-1}\right)}\right)
    \cdot \mathbbm{1}_{X^{\left(i\right)}_m > X^{\left(i\right)}_{m-1}, m\leq\tau^{\left(i\right)}}\,, \\
    &\hphantom{\;:=\;} + \sum_{m = 1}^n \; \sum_{j = 1, j \neq i}^K \; \sum_{z = 0}^{\infty}
    \left(\frac{1}{w\left(m-1,z\right)} - \frac{1}{w\left(m,z\right)}\right)
    \cdot \underbrace{\mathbbm{1}_{\left\{X^{\left(j\right)}_{m-1}, X^{\left(j\right)}_m\right\} = \left\{z, z + 1\right\}, m\leq\tau^{\left(i\right)}, z < X^{\left(i\right)}_m}}_{= 1 \textrm{ for at most one pair of }j,z}\,.
  \end{align*}
  $M^{\left(i\right)}_n$ is nonnegative by definition of $F$, and $H^{\left(i\right)}_n \geq M^{\left(i\right)}_n \geq 0$
  since edge weights can only increase and therefore, all terms in the sums in the definition of $H^{\left(i\right)}_n$
  are nonnegative. $H^{\left(i\right)}_n$ is a martingale: setting
  \begin{align*}
    d^{\left(i\right)}_n :=
    H^{\left(i\right)}_n - H^{\left(i\right)}_{n-1}
    &= \underbrace{M^{\left(i\right)}_n - M^{\left(i\right)}_{n-1}}_{:= e^{\left(i\right)}_n}
    + \underbrace{\left(\frac{1}{w\left(n-1,X^{\left(i\right)}_{n-1}\right)} - \frac{1}{w\left(n,X^{\left(i\right)}_{n-1}\right)}\right)
    \cdot \mathbbm{1}_{X^{\left(i\right)}_n > X^{\left(i\right)}_{n-1}, n\leq\tau^{\left(i\right)}}}_{:= f^{\left(i\right)}_n} \\
    &\hphantom{\;=\;}+ \underbrace{\sum_{j = 1, j \neq i}^K \; \sum_{z = 0}^{\infty}
    \left(\frac{1}{w\left(n-1,z\right)} - \frac{1}{w\left(n,z\right)}\right)
    \cdot \mathbbm{1}_{\left\{X^{\left(j\right)}_{n-1}, X^{\left(j\right)}_n\right\} = \left\{z, z + 1\right\}, n\leq\tau^{\left(i\right)}, z < X^{\left(i\right)}_n}}_{:= g^{\left(i\right)}_n}\,,
  \end{align*}
  we have to show that $\Ex{d^{\left(i\right)}_n \midd| \calG_{n-1}} = 0$. We have:
  \begin{itemize}
    \item if $n - 1 \geq \tau^{\left(i\right)}$, then $d^{\left(i\right)}_n = 0$. Hence, it suffices to consider the case $X^{\left(i\right)}_{n-1} = z > 0$
    and $\tau^{\left(i\right)} \geq n$.
    \item with probability $\frac{1}{K}$, the walker $i$ jumps at time $n-1$. In this case, $g^{\left(i\right)}_n = 0$
    since no other walker can jump and the indicator variable in $g^{\left(i\right)}_n$ is therefore $0$.
    If he jumps to the right (with probability
    $\frac{1}{K} \cdot \frac{w\left(n - 1, z\right)}{w\left(n - 1, z - 1\right) + w\left(n - 1, z\right)}$), then
    $e^{\left(i\right)}_n = w\left(n, z\right)^{-1}$ and $f^{\left(i\right)}_n = w\left(n-1,z\right)^{-1} - w\left(n,z\right)^{-1}$, hence $d^{\left(i\right)}_n = w\left(n-1,z\right)^{-1}$.
    If he jumps left (with probability
    $\frac{1}{K} \cdot \frac{w\left(n - 1, z - 1\right)}{w\left(n - 1, z - 1\right) + w\left(n - 1, z\right)}$), then
    $e^{\left(i\right)}_n = -w\left(n - 1, z - 1\right)^{-1}$ and $f^{\left(i\right)}_n = 0$, hence $d^{\left(i\right)}_n = -w\left(n - 1, z - 1\right)^{-1}$.
    \item with probability $\frac{K-1}{K}$, the walker $i$ does not jump. In this case, $f^{\left(i\right)}_n = 0$
    since the indicator variable in $f^{\left(i\right)}_n$ is therefore $0$.
    The value of $M^{\left(i\right)}_{n-1}$ now changes
    (that is, $e^{\left(i\right)}_n \neq 0$) if one of the other $K-1$ walkers crosses one of the edges between the nodes $0$ and $z$.
    At the same time, $g^{\left(i\right)}_n \neq 0$ only in this exact case. Now assume the walker $j \neq i$
    crosses the edge $y$ with $0 \leq y < z$. Then $e^{\left(i\right)}_n = \frac{1}{w\left(n,y\right)} - \frac{1}{w\left(n-1,y\right)}$ and
    $g^{\left(i\right)}_n = \frac{1}{w\left(n-1,y\right)} - \frac{1}{w\left(n,y\right)}$, hence $d^{\left(i\right)}_n = 0$.
    \item conditioned on $X^{\left(i\right)}_{n-1} = z > 0$ and $\tau^{\left(i\right)} \geq n$ (both events measurable w.r.t.~$\calG_{n-1}$), we can therefore conclude
    \begin{align*}
      \Ex{d^{\left(i\right)}_n \midd| \calG_{n-1}} &= \frac{1}{K} \cdot \frac{1}{w\left(n - 1, z - 1\right) + w\left(n - 1, z\right)} \cdot \left(
      \frac{w\left(n - 1, z\right)}{w\left(n - 1, z\right)} - \frac{w\left(n - 1, z - 1\right)}{w\left(n - 1, z - 1\right)}\right) \\
      &= \frac{1}{K} \cdot \frac{1}{w\left(n - 1, z - 1\right) + w\left(n - 1, z\right)} \cdot \left(1 - 1\right) = 0\,.
    \end{align*}
    By the same arguments, but only considering $e^{\left(i\right)}_n$, we can show that $M^{\left(i\right)}_n$
    is a supermartingale.
  \end{itemize}
  As a nonnegative martingale, $H^{\left(i\right)}_n$ converges almost surely.
  
  We just showed this for all walkers $i$ ($1 \leq i \leq K$) with $X^{\left(i\right)}_0 \geq 0$.
  Further observe that for such a walker $i$, we have, on the event $B^{\left(i\right)}_n = \left\{X^{\left(i\right)}_n > X^{\left(i\right)}_{n-1}, n \leq \tau^{\left(i\right)}, w\left(n-1, X^{\left(i\right)}_{n-1}\right) = 1\right\}$,
  that $e^{\left(i\right)}_n = w\left(n, X^{\left(i\right)}_{n-1}\right)^{-1}$, $f^{\left(i\right)}_n = w\left(n - 1, X^{\left(i\right)}_{n-1}\right)^{-1} - w\left(n, X^{\left(i\right)}_{n-1}\right)^{-1}$, $g^{\left(i\right)}_n = 0$ and hence $d^{\left(i\right)}_n = 1$.
  Thus, by convergence, only a finite number of the events $B^{\left(i\right)}_n$ can occur
  for every such walker $i$.

  Now define $\Gamma$ to be the set of edges between two nonnegative integers to the right of the integer
  $\max\left\{ X^{\left(i\right)}_0 : 1 \leq i \leq K \right\}$ for which the initial weight
  was $1$ (all but finitely many edges meet the latter criterion), and further define the event
  \begin{align*}
    D_n &= \left\{ \exists j: X^{\left(j\right)}_0 \geq 0 \textrm{, an edge in } \Gamma \textrm{ is crossed between time } n-1 \textrm{ and } n \textrm{ for the first time by}\right.\\
    &\hphantom{\;=\;\{} \left.\vphantom{X^{\left(j\right)}_0}\textrm{any walker, the crossing walker is }j \textrm{ and } n \leq \tau^{\left(j\right)} \right\}\,.
  \end{align*}
  Clearly, $D_n \subseteq B^{\left(i\right)}_n$ for some random walker $i$ with $X^{\left(i\right)}_0 \geq 0$, hence only a finite number of
  the events $D_n$ can occur.

  Now the proof cannot be continued along \cite[Lemma 3.0]{errwz} since the walkers starting to the left
  of $0$ and the walkers which reach $0$ can
  later cross edges to the right of $0$ without triggering $D_n$ and the other walkers can then follow
  them without triggering $D_n$. So, we only proved that walkers which never go to $0$ and start to
  the right of $0$ cannot
  visit infinitely many edges which have not been visited before by any other walker.
\end{oproof}

\ifthenelse{\boolean{extendedversion}}{
\begin{figure}[H]
  \begin{center}
    \begin{tikzpicture}[scale=0.9]
  \draw[dashed] (-4.32,0) -- (-3.6,0);
  \draw[dashed] (10.8,0) -- (11.52,0);
  \draw (-3.6,0) -- node[above] {$1$} (-2.4,0) -- node[above] {$1$} (-1.2,0) -- node[above] {$1$} node[below=3mm] {$M^{\left(3\right)}_0\;=$} (0,0) node[above=1mm,sectionblue] {$0$} -- node[above] {$1$} node[below=3mm] {$\vphantom{M^{\left(m\right)}_n}1$} (1.2,0) -- node[above] {$1$} node[below=3mm] {$\vphantom{M^{\left(m\right)}_n}+\;1$} (2.4,0) -- node[above] {$1$} node[below=3mm] {$\vphantom{M^{\left(m\right)}_n}+\;1$} (3.6,0) -- node[above] {$1$} node[below=3mm] {$\vphantom{M^{\left(m\right)}_n}+\;1$} (4.8,0) -- node[above] {$1$} node[below=3mm] {$\vphantom{M^{\left(m\right)}_n}+\;1$} (6,0) -- node[above] {$1$} node[below=3mm] {$\vphantom{M^{\left(m\right)}_n}=\;5$} (7.2,0) -- node[above] {$1$} (8.4,0) -- node[above] {$1$} (9.6,0) -- node[above] {$1$} (10.8,0);
  \node[circle,fill=black,inner sep=0.7mm] at (-3.6,0) {};
  \node[circle,fill=black,inner sep=0.7mm] at (-2.4,0) {};
  \node[circle,draw=black,inner sep=1.1mm] at (-2.4,0) {};
  \node[circle,fill=black,inner sep=0.7mm] at (-1.2,0) {};
  \node[circle,fill=sectionblue,inner sep=0.7mm] at (0,0) {};
  \node[circle,fill=black,inner sep=0.7mm] at (1.2,0) {};
  \node[circle,fill=black,inner sep=0.7mm] at (2.4,0) {};
  \node[circle,draw=black,inner sep=1.1mm] at (2.4,0) {};
  \node[circle,fill=black,inner sep=0.7mm] at (3.6,0) {};
  \node[circle,fill=black,inner sep=0.7mm] at (4.8,0) {};
  \node[circle,fill=sectionblue,inner sep=0.7mm] at (6,0) {};
  \node[circle,draw=sectionblue,inner sep=1.1mm] at (6,0) {};
  \node[circle,fill=black,inner sep=0.7mm] at (7.2,0) {};
  \node[circle,fill=black,inner sep=0.7mm] at (8.4,0) {};
  \node[circle,draw=black,inner sep=1.1mm] at (8.4,0) {};
  \node[circle,fill=black,inner sep=0.7mm] at (9.6,0) {};
  \node[circle,fill=black,inner sep=0.7mm] at (10.8,0) {};
\end{tikzpicture}
 \\[1em]
    \begin{tikzpicture}[scale=0.9]
  \draw[line width=2pt] (6,0) -- (7.2,0);
  \draw[dashed] (-4.32,0) -- (-3.6,0);
  \draw[dashed] (10.8,0) -- (11.52,0);
  \draw (-3.6,0) -- node[above] {$1$} (-2.4,0) -- node[above] {$1$} (-1.2,0) -- node[above] {$1$} node[below=3mm] {$M^{\left(3\right)}_1\;=$} (0,0) node[above=1mm,sectionblue] {$0$} -- node[above] {$1$} node[below=3mm] {$\vphantom{M^{\left(m\right)}_n}1$} (1.2,0) -- node[above] {$1$} node[below=3mm] {$\vphantom{M^{\left(m\right)}_n}+\;1$} (2.4,0) -- node[above] {$1$} node[below=3mm] {$\vphantom{M^{\left(m\right)}_n}+\;1$} (3.6,0) -- node[above] {$1$} node[below=3mm] {$\vphantom{M^{\left(m\right)}_n}+\;1$} (4.8,0) -- node[above] {$1$} node[below=3mm] {$\vphantom{M^{\left(m\right)}_n}+\;1$} (6,0) -- node[above] {$2$} node[below=3mm] {$\vphantom{M^{\left(m\right)}_n}+\;\frac{1}{2}$} (7.2,0) -- node[above] {$1$} node[below=3mm] {$\vphantom{M^{\left(m\right)}_n}=\;\frac{11}{2}$} (8.4,0) -- node[above] {$1$} (9.6,0) -- node[above] {$1$} (10.8,0);
  \node[circle,fill=black,inner sep=0.7mm] at (-3.6,0) {};
  \node[circle,fill=black,inner sep=0.7mm] at (-2.4,0) {};
  \node[circle,draw=black,inner sep=1.1mm] at (-2.4,0) {};
  \node[circle,fill=black,inner sep=0.7mm] at (-1.2,0) {};
  \node[circle,fill=sectionblue,inner sep=0.7mm] at (0,0) {};
  \node[circle,fill=black,inner sep=0.7mm] at (1.2,0) {};
  \node[circle,fill=black,inner sep=0.7mm] at (2.4,0) {};
  \node[circle,draw=black,inner sep=1.1mm] at (2.4,0) {};
  \node[circle,fill=black,inner sep=0.7mm] at (3.6,0) {};
  \node[circle,fill=black,inner sep=0.7mm] at (4.8,0) {};
  \node[circle,fill=black,inner sep=0.7mm] at (6,0) {};
  \node[circle,fill=sectionblue,inner sep=0.7mm] at (7.2,0) {};
  \node[circle,draw=sectionblue,inner sep=1.1mm] at (7.2,0) {};
  \node[circle,fill=black,inner sep=0.7mm] at (8.4,0) {};
  \node[circle,draw=black,inner sep=1.1mm] at (8.4,0) {};
  \node[circle,fill=black,inner sep=0.7mm] at (9.6,0) {};
  \node[circle,fill=black,inner sep=0.7mm] at (10.8,0) {};
\end{tikzpicture}
 \\[1em]
    \begin{tikzpicture}[scale=0.9]
  \draw[line width=2pt] (1.2,0) -- (2.4,0);
  \draw[dashed] (-4.32,0) -- (-3.6,0);
  \draw[dashed] (10.8,0) -- (11.52,0);
  \draw (-3.6,0) -- node[above] {$1$} (-2.4,0) -- node[above] {$1$} (-1.2,0) -- node[above] {$1$} node[below=3mm] {$M^{\left(3\right)}_2\;=$} (0,0) node[above=1mm,sectionblue] {$0$} -- node[above] {$1$} node[below=3mm] {$\vphantom{M^{\left(m\right)}_n}1$} (1.2,0) -- node[above] {$2$} node[below=3mm] {$\vphantom{M^{\left(m\right)}_n}+\;\frac{1}{2}$} (2.4,0) -- node[above] {$1$} node[below=3mm] {$\vphantom{M^{\left(m\right)}_n}+\;1$} (3.6,0) -- node[above] {$1$} node[below=3mm] {$\vphantom{M^{\left(m\right)}_n}+\;1$} (4.8,0) -- node[above] {$1$} node[below=3mm] {$\vphantom{M^{\left(m\right)}_n}+\;1$} (6,0) -- node[above] {$2$} node[below=3mm] {$\vphantom{M^{\left(m\right)}_n}+\;\frac{1}{2}$} (7.2,0) -- node[above] {$1$} node[below=3mm] {$\vphantom{M^{\left(m\right)}_n}=\;5$} (8.4,0) -- node[above] {$1$} (9.6,0) -- node[above] {$1$} (10.8,0);
  \node[circle,fill=black,inner sep=0.7mm] at (-3.6,0) {};
  \node[circle,fill=black,inner sep=0.7mm] at (-2.4,0) {};
  \node[circle,draw=black,inner sep=1.1mm] at (-2.4,0) {};
  \node[circle,fill=black,inner sep=0.7mm] at (-1.2,0) {};
  \node[circle,fill=sectionblue,inner sep=0.7mm] at (0,0) {};
  \node[circle,fill=black,inner sep=0.7mm] at (1.2,0) {};
  \node[circle,draw=black,inner sep=1.1mm] at (1.2,0) {};
  \node[circle,fill=black,inner sep=0.7mm] at (2.4,0) {};
  \node[circle,fill=black,inner sep=0.7mm] at (3.6,0) {};
  \node[circle,fill=black,inner sep=0.7mm] at (4.8,0) {};
  \node[circle,fill=black,inner sep=0.7mm] at (6,0) {};
  \node[circle,fill=sectionblue,inner sep=0.7mm] at (7.2,0) {};
  \node[circle,draw=sectionblue,inner sep=1.1mm] at (7.2,0) {};
  \node[circle,fill=black,inner sep=0.7mm] at (8.4,0) {};
  \node[circle,draw=black,inner sep=1.1mm] at (8.4,0) {};
  \node[circle,fill=black,inner sep=0.7mm] at (9.6,0) {};
  \node[circle,fill=black,inner sep=0.7mm] at (10.8,0) {};
\end{tikzpicture}
 \\[1em]
    \begin{tikzpicture}[scale=0.9]
  \draw[line width=2pt] (7.2,0) -- (8.4,0);
  \draw[dashed] (-4.32,0) -- (-3.6,0);
  \draw[dashed] (10.8,0) -- (11.52,0);
  \draw (-3.6,0) -- node[above] {$1$} (-2.4,0) -- node[above] {$1$} (-1.2,0) -- node[above] {$1$} node[below=3mm] {$M^{\left(3\right)}_3\;=$} (0,0) node[above=1mm,sectionblue] {$0$} -- node[above] {$1$} node[below=3mm] {$\vphantom{M^{\left(m\right)}_n}1$} (1.2,0) -- node[above] {$2$} node[below=3mm] {$\vphantom{M^{\left(m\right)}_n}+\;\frac{1}{2}$} (2.4,0) -- node[above] {$1$} node[below=3mm] {$\vphantom{M^{\left(m\right)}_n}+\;1$} (3.6,0) -- node[above] {$1$} node[below=3mm] {$\vphantom{M^{\left(m\right)}_n}+\;1$} (4.8,0) -- node[above] {$1$} node[below=3mm] {$\vphantom{M^{\left(m\right)}_n}+\;1$} (6,0) -- node[above] {$2$} node[below=3mm] {$\vphantom{M^{\left(m\right)}_n}+\;\frac{1}{2}$} (7.2,0) -- node[above] {$2$} node[below=3mm] {$\vphantom{M^{\left(m\right)}_n}=\;5$} (8.4,0) -- node[above] {$1$} (9.6,0) -- node[above] {$1$} (10.8,0);
  \node[circle,fill=black,inner sep=0.7mm] at (-3.6,0) {};
  \node[circle,fill=black,inner sep=0.7mm] at (-2.4,0) {};
  \node[circle,draw=black,inner sep=1.1mm] at (-2.4,0) {};
  \node[circle,fill=black,inner sep=0.7mm] at (-1.2,0) {};
  \node[circle,fill=sectionblue,inner sep=0.7mm] at (0,0) {};
  \node[circle,fill=black,inner sep=0.7mm] at (1.2,0) {};
  \node[circle,draw=black,inner sep=1.1mm] at (1.2,0) {};
  \node[circle,fill=black,inner sep=0.7mm] at (2.4,0) {};
  \node[circle,fill=black,inner sep=0.7mm] at (3.6,0) {};
  \node[circle,fill=black,inner sep=0.7mm] at (4.8,0) {};
  \node[circle,fill=black,inner sep=0.7mm] at (6,0) {};
  \node[circle,fill=sectionblue,inner sep=0.7mm] at (7.2,0) {};
  \node[circle,draw=sectionblue,inner sep=1.1mm] at (7.2,0) {};
  \node[circle,draw=black,inner sep=1.5mm] at (7.2,0) {};
  \node[circle,fill=black,inner sep=0.7mm] at (8.4,0) {};
  \node[circle,fill=black,inner sep=0.7mm] at (9.6,0) {};
  \node[circle,fill=black,inner sep=0.7mm] at (10.8,0) {};
\end{tikzpicture}
 \\[1em]
    \begin{tikzpicture}[scale=0.9]
  \draw[line width=2pt] (6,0) -- (7.2,0);
  \draw[dashed] (-4.32,0) -- (-3.6,0);
  \draw[dashed] (10.8,0) -- (11.52,0);
  \draw (-3.6,0) -- node[above] {$1$} (-2.4,0) -- node[above] {$1$} (-1.2,0) -- node[above] {$1$} node[below=3mm] {$M^{\left(3\right)}_4\;=$} (0,0) node[above=1mm,sectionblue] {$0$} -- node[above] {$1$} node[below=3mm] {$\vphantom{M^{\left(m\right)}_n}1$} (1.2,0) -- node[above] {$2$} node[below=3mm] {$\vphantom{M^{\left(m\right)}_n}+\;\frac{1}{2}$} (2.4,0) -- node[above] {$1$} node[below=3mm] {$\vphantom{M^{\left(m\right)}_n}+\;1$} (3.6,0) -- node[above] {$1$} node[below=3mm] {$\vphantom{M^{\left(m\right)}_n}+\;1$} (4.8,0) -- node[above] {$1$} node[below=3mm] {$\vphantom{M^{\left(m\right)}_n}+\;1$} (6,0) -- node[above] {$3$} node[below=3mm] {$\vphantom{M^{\left(m\right)}_n}=\;\frac{9}{2}$} (7.2,0) -- node[above] {$2$} (8.4,0) -- node[above] {$1$} (9.6,0) -- node[above] {$1$} (10.8,0);
  \node[circle,fill=black,inner sep=0.7mm] at (-3.6,0) {};
  \node[circle,fill=black,inner sep=0.7mm] at (-2.4,0) {};
  \node[circle,draw=black,inner sep=1.1mm] at (-2.4,0) {};
  \node[circle,fill=black,inner sep=0.7mm] at (-1.2,0) {};
  \node[circle,fill=sectionblue,inner sep=0.7mm] at (0,0) {};
  \node[circle,fill=black,inner sep=0.7mm] at (1.2,0) {};
  \node[circle,draw=black,inner sep=1.1mm] at (1.2,0) {};
  \node[circle,fill=black,inner sep=0.7mm] at (2.4,0) {};
  \node[circle,fill=black,inner sep=0.7mm] at (3.6,0) {};
  \node[circle,fill=black,inner sep=0.7mm] at (4.8,0) {};
  \node[circle,fill=sectionblue,inner sep=0.7mm] at (6,0) {};
  \node[circle,draw=sectionblue,inner sep=1.1mm] at (6,0) {};
  \node[circle,fill=black,inner sep=0.7mm] at (7.2,0) {};
  \node[circle,draw=black,inner sep=1.1mm] at (7.2,0) {};
  \node[circle,fill=black,inner sep=0.7mm] at (8.4,0) {};
  \node[circle,fill=black,inner sep=0.7mm] at (9.6,0) {};
  \node[circle,fill=black,inner sep=0.7mm] at (10.8,0) {};
\end{tikzpicture}

    \caption{Evolution of the supermartingale during 4 steps of the ERRW with multiple walkers}
    \label{fig:errw_multiple_mart}
  \end{center}
  In the figure above, we consider $K = 4$ walkers whose positions are indicated by the circled nodes. The value
  of $M_n^{\left(3\right)}$ is shown for the first $4$ steps, and the corresponding walker $X^{\left(3\right)}$
  is highlighted in \textcolor{sectionblue}{blue}.
\end{figure}
}{}

\begin{oproof}[\autoref{lem:labelExchange}]
  The proof idea is the following: whenever $X^{\left(1\right)}$ and $X^{\left(2\right)}$ meet, we
  can randomly exchange their labels, i.e.~we can randomly decide whether we want to rename $X^{\left(1\right)}$
  to $X^{\left(2\right)}$ and vice versa, and the law of the edge-reinforced random walk with the
  two walkers is invariant under such relabelings because the only distinguishing feature of a random
  walker is his position. But now, to construct counterexamples to the two statements in \autoref{lem:labelExchange},
  we would have to choose a fixed labeling for infinitely many times at which the walkers meet. But
  if we randomize the labeling with a sequence of independent Bernoulli random variables, then the probability of
  choosing a certain fixed labeling at infinitely many points in the sequence is $0$, and since the
  law was invariant under random relabeling, it follows that the probability of any such counterexample
  is $0$. We continue with the formal proof.

  Set $\tau_1 := \inf\left\{ n \geq 0 : X^{\left(1\right)}_n = X^{\left(2\right)}_n \right\}$
  and $\tau_{m + 1} := \inf\left\{ n > \tau_m : X^{\left(1\right)}_n = X^{\left(2\right)}_n \right\}$.
  If $X^{\left(1\right)}, X^{\left(2\right)}$ meet infinitely
  often, then $\forall n: \tau_n < \infty$, but the construction also works if this is not the case. Let $\left(\omega_m\right)_{m \geq 1}$ be a sequence of
  iid random variables with $\Prb{\omega_m = 1} = \frac{1}{2} = \Prb{\omega_m = 0}$ (the $\omega_m$
  are also independent of $\calG_n$ for all $n$, i.e.~independent of the edge-reinforced random walk).
  Define $\widetilde{X}^{\left(1\right)}_n$ and $\widetilde{X}^{\left(2\right)}_n$ as follows (with $\omega_0 = 0$ and $\tau_0 = -1$):
  \begin{align*}
    \widetilde{X}^{\left(1\right)}_n := X^{\left(1,\omega\right)}_n &= \sum_{m \geq 0} \left( \left(1 - \omega_m\right) X^{\left(1\right)}_n + \omega_m X^{\left(2\right)}_n \right) \cdot
    \mathbbm{1}_{\tau_m < n \leq \tau_{m+1}}\,, \\
    \widetilde{X}^{\left(2\right)}_n := X^{\left(2,\omega\right)}_n &= \sum_{m \geq 0} \left( \left(1 - \omega_m\right) X^{\left(2\right)}_n + \omega_m X^{\left(1\right)}_n \right) \cdot
    \mathbbm{1}_{\tau_m < n \leq \tau_{m+1}}\,.
  \end{align*}
  Note that the sums consist of a single term. $\omega_m = 1$ means that we switch the labels of $X^{\left(1\right)}$
  and $X^{\left(2\right)}$ during the time interval $\left(\tau_m,\tau_{m+1}\right]$.
  
  \ifthenelse{\boolean{extendedversion}}{
  \begin{figure}[H]
    \begin{center}
      \begin{tikzpicture}[scale=0.8]
  \draw[-{Latex[length=1.8mm,width=1.8mm]},tumGray] (0,0) -- (15.5,0) node[right] {time};
  \draw[sectionblue,line width=1pt] (0,-1) -- (0.5,-1) -- (1,-1) -- (1.5,-1) -- (2,-1) -- (2.5,-1) -- (3,-1.5) -- (3.5,-2) -- (4,-2.5) -- (4.5,-2.5) -- (5,-2.5) -- (5.5,-2) -- (6,-2.5) -- (6.5,-2.5) -- (7,-2.5) -- (7.5,-2) -- (8,-1.5) -- (8.5,-1.5) -- (9,-1.5) -- (9.5,-1.5) -- (10,-1.5) -- (10.5,-1.5) -- (11,-1) -- (11.5,-1) -- (12,-1) -- (12.5,-1) -- (13,-0.5) -- (13.5,-1) -- (14,-1.5) -- (14.5,-1) -- (15,-0.5) node[right] {$X^{\left(1\right)}$};
  \draw[tumOrange,line width=1pt] (0,1) -- (0.5,0.5) -- (1,0) -- (1.5,-0.5) -- (2,-1) -- (2.5,-0.5) -- (3,-0.5) -- (3.5,-0.5) -- (4,-0.5) -- (4.5,-1) -- (5,-1.5) -- (5.5,-1.5) -- (6,-1.5) -- (6.5,-1) -- (7,-1.5) -- (7.5,-1.5) -- (8,-1.5) -- (8.5,-2) -- (9,-1.5) -- (9.5,-1) -- (10,-0.5) -- (10.5,-1) -- (11,-1) -- (11.5,-1.5) -- (12,-2) -- (12.5,-1.5) -- (13,-1.5) -- (13.5,-1.5) -- (14,-1.5) -- (14.5,-1.5) -- (15,-1.5) node[right] {$X^{\left(2\right)}$};
  \draw[-{Latex[length=1.8mm,width=1.8mm]}] (0,-3) -- (0,1.5) node[above] {$\bbZ$};
  \draw (0,0) -- (-0.2,0) node[left] {$0$};
  \draw (0,-2.5) -- (-0.2,-2.5);
  \draw (0,-2) -- (-0.2,-2);
  \draw (0,-1.5) -- (-0.2,-1.5);
  \draw (0,-1) -- (-0.2,-1);
  \draw (0,-0.5) -- (-0.2,-0.5);
  \draw (0,0) -- (-0.2,0);
  \draw (0,0.5) -- (-0.2,0.5);
  \draw (0,1) -- (-0.2,1);
  \node[circle,fill=black,inner sep=0.7mm] at (2,-1) {};
  \node[circle,fill=black,inner sep=0.7mm] at (8,-1.5) {};
  \node[circle,fill=black,inner sep=0.7mm] at (9,-1.5) {};
  \node[circle,fill=black,inner sep=0.7mm] at (11,-1) {};
  \node[circle,fill=black,inner sep=0.7mm] at (14,-1.5) {};
  \node[below=2mm] at (2,-1) {$\tau_1 = 4$};
  \node[above=2mm] at (8,-1.5) {$\tau_2 = 16$};
  \node[below=2mm] at (9,-1.5) {$\tau_3 = 18$};
  \node[above=2mm] at (11,-1) {$\tau_4 = 22$};
  \node[below=2mm] at (14,-1.5) {$\tau_5 = 28$};
\end{tikzpicture}
 \\[1em]
      \begin{tikzpicture}[scale=0.8]
  \fill[tumGray!20!white] (2,-3) rectangle ++(6,4.5);
  \fill[tumGray!20!white] (9,-3) rectangle ++(2,4.5);
  \fill[tumGray!20!white] (14,-3) rectangle ++(1.15,4.5);
  \draw[-{Latex[length=1.8mm,width=1.8mm]},tumGray] (0,0) -- (15.5,0) node[right] {time};
  \draw[sectionblue!50!white,line width=4pt] (0,-1) -- (0.5,-1) -- (1,-1) -- (1.5,-1) -- (2,-1) -- (2.5,-0.5) -- (3,-0.5) -- (3.5,-0.5) -- (4,-0.5) -- (4.5,-1) -- (5,-1.5) -- (5.5,-1.5) -- (6,-1.5) -- (6.5,-1) -- (7,-1.5) -- (7.5,-1.5) -- (8,-1.5) -- (8.5,-1.5) -- (9,-1.5) -- (9.5,-1) -- (10,-0.5) -- (10.5,-1) -- (11,-1) -- (11.5,-1) -- (12,-1) -- (12.5,-1) -- (13,-0.5) -- (13.5,-1) -- (14,-1.5) -- (14.5,-1.5) -- (15,-1.5) node[right] {$\widetilde{X}^{\left(1\right)}$};
  \draw[tumOrange!50!white,line width=4pt] (0,1) -- (0.5,0.5) -- (1,0) -- (1.5,-0.5) -- (2,-1) -- (2.5,-1) -- (3,-1.5) -- (3.5,-2) -- (4,-2.5) -- (4.5,-2.5) -- (5,-2.5) -- (5.5,-2) -- (6,-2.5) -- (6.5,-2.5) -- (7,-2.5) -- (7.5,-2) -- (8,-1.5) -- (8.5,-2) -- (9,-1.5) -- (9.5,-1.5) -- (10,-1.5) -- (10.5,-1.5) -- (11,-1) -- (11.5,-1.5) -- (12,-2) -- (12.5,-1.5) -- (13,-1.5) -- (13.5,-1.5) -- (14,-1.5) -- (14.5,-1) -- (15,-0.5) node[right] {$\widetilde{X}^{\left(2\right)}$};
  \draw[sectionblue,line width=1pt] (0,-1) -- (0.5,-1) -- (1,-1) -- (1.5,-1) -- (2,-1) -- (2.5,-1) -- (3,-1.5) -- (3.5,-2) -- (4,-2.5) -- (4.5,-2.5) -- (5,-2.5) -- (5.5,-2) -- (6,-2.5) -- (6.5,-2.5) -- (7,-2.5) -- (7.5,-2) -- (8,-1.5) -- (8.5,-1.5) -- (9,-1.5) -- (9.5,-1.5) -- (10,-1.5) -- (10.5,-1.5) -- (11,-1) -- (11.5,-1) -- (12,-1) -- (12.5,-1) -- (13,-0.5) -- (13.5,-1) -- (14,-1.5) -- (14.5,-1) -- (15,-0.5);
  \draw[tumOrange,line width=1pt] (0,1) -- (0.5,0.5) -- (1,0) -- (1.5,-0.5) -- (2,-1) -- (2.5,-0.5) -- (3,-0.5) -- (3.5,-0.5) -- (4,-0.5) -- (4.5,-1) -- (5,-1.5) -- (5.5,-1.5) -- (6,-1.5) -- (6.5,-1) -- (7,-1.5) -- (7.5,-1.5) -- (8,-1.5) -- (8.5,-2) -- (9,-1.5) -- (9.5,-1) -- (10,-0.5) -- (10.5,-1) -- (11,-1) -- (11.5,-1.5) -- (12,-2) -- (12.5,-1.5) -- (13,-1.5) -- (13.5,-1.5) -- (14,-1.5) -- (14.5,-1.5) -- (15,-1.5);
  \draw[-{Latex[length=1.8mm,width=1.8mm]}] (0,-3) -- (0,1.5) node[above] {$\bbZ$};
  \draw (0,0) -- (-0.2,0) node[left] {$0$};
  \draw (0,-2.5) -- (-0.2,-2.5);
  \draw (0,-2) -- (-0.2,-2);
  \draw (0,-1.5) -- (-0.2,-1.5);
  \draw (0,-1) -- (-0.2,-1);
  \draw (0,-0.5) -- (-0.2,-0.5);
  \draw (0,0) -- (-0.2,0);
  \draw (0,0.5) -- (-0.2,0.5);
  \draw (0,1) -- (-0.2,1);
  \node[circle,fill=black,inner sep=0.7mm] at (2,-1) {};
  \node[circle,fill=black,inner sep=0.7mm] at (8,-1.5) {};
  \node[circle,fill=black,inner sep=0.7mm] at (9,-1.5) {};
  \node[circle,fill=black,inner sep=0.7mm] at (11,-1) {};
  \node[circle,fill=black,inner sep=0.7mm] at (14,-1.5) {};
  \node[below=2mm] at (2,-1) {$\tau_1 = 4$};
  \node[above=2mm] at (8,-1.5) {$\tau_2 = 16$};
  \node[below=2mm] at (9,-1.5) {$\tau_3 = 18$};
  \node[above=2mm] at (11,-1) {$\tau_4 = 22$};
  \node[below=2mm] at (14,-1.5) {$\tau_5 = 28$};
  \node at (1,1.25) {$\omega_0 = 0$};
  \node at (5,1.25) {$\omega_1 = 1$};
  \node at (8.5,1.25) {$\omega_2 = 0$};
  \node at (10,1.25) {$\omega_3 = 1$};
  \node at (12.5,1.25) {$\omega_4 = 0$};
  \node at (14.58,1.25) {$\omega_5 = 1$};
\end{tikzpicture}

      \caption{Label exchange lemma}
      \label{fig:label_exchange}
    \end{center}
    The figure above illustrates how the ``label exchange'' of the two walkers works. Two sample paths
    for the two walkers are drawn, together with meeting points and the result of the label exchange.
  \end{figure}
  }{}
  
  If we consider
  $\left(X^{\left(1\right)},X^{\left(2\right)}\right)$ and $\left(\widetilde{X}^{\left(1\right)},\widetilde{X}^{\left(2\right)}\right)$
  as sequences of pairs of integers, then we have
  \begin{align}
    \left(\widetilde{X}^{\left(i\right)}_n\right)_{1\leq i\leq 2, n\geq 0} &\overset{\textrm{d}}{=}
    \left(X^{\left(i\right)}_n\right)_{1\leq i\leq 2, n\geq 0}\,.
    \label{eq:eqInDist}
  \end{align}
  The equality in distribution follows from the above-mentioned invariance of the law of the random walk
  under relabelings at meeting times which is quite intuitive, and could be proved formally by looking at cylinder events,
  for example.\\
\noindent
  We now show that any counterexamples to statements \ref{lem:labelExchangeBothInf}, \ref{lem:labelExchangeBothVisInf} or \ref{lem:labelExchangeBothRec}
  have probability $0$:
  \begin{enumerate}[(i)]
    \item Let $A$ be the event that one of the walkers $X^{\left(1\right)},X^{\left(2\right)}$ has finite range while the other one has infinite
    range, and that they meet infinitely often. It suffices to show that $\Prb{A} = 0$. Denote by $\bP$ the probability measure induced
    by the edge reinforced random walk alone and by $\bQ$ the probability measure induced by the
    sequence $\left(\omega_m\right)_{m \geq 1}$ alone. Then, by \eqref{eq:eqInDist}, we have
    \begin{align*}
      \Prb{A} = &\int \int \mathbbm{1}_{B} \dx{\bQ} \dx{\bP} \\
      \textrm{where } B := &\left\{\textrm{one of }\widetilde{X}^{\left(1\right)},\widetilde{X}^{\left(2\right)}\textrm{ has finite range while}\right.\\
      & \left.\hphantom{\;\{}\vphantom{\widetilde{X}^{\left(1\right)}}\textrm{the other has infinite range, they meet infinitely often}\right\}\,.
    \end{align*}
    We have to show that the inner integral is $0$ almost surely with respect to $\bP$. Consider
    fixed walker sequences $X^{\left(1\right)}$ and $X^{\left(2\right)}$. If one of
    $\widetilde{X}^{\left(1\right)},\widetilde{X}^{\left(2\right)}$ should have finite range while the
    other has infinite range, then, by definition of $\widetilde{X}^{\left(1\right)}$ and $\widetilde{X}^{\left(2\right)}$, at
    least one of $X^{\left(1\right)},X^{\left(2\right)}$ must have infinite range. Of course, by definition,
    we also have that $\widetilde{X}^{\left(1\right)},\widetilde{X}^{\left(2\right)}$ meet infinitely often
    if, and only if, $X^{\left(1\right)},X^{\left(2\right)}$ meet infinitely often. Hence, the indicator variable in the
    integral above can only be $1$ in the case where one of the walkers $X^{\left(1\right)},X^{\left(2\right)}$ has infinite range and the two walkers
    meet infinitely often, so we only need to show that in this particular case, the inner integral is
    still $0$ almost surely.

    Assume $X^{\left(1\right)}$ does not have finite range (w.l.o.g.). Then, for every $n$, one can
    find $m$ such that between times $\tau_m$ and $\tau_{m+1}$ (all $\tau_m$ are finite if the two
    walkers meet infinitely often), $X^{\left(1\right)}$ visits a node at distance at least $n$ from
    the integer $0$. Call these times $\tau_{m_n}$ with $m_n$ strictly increasing in $n$ (w.l.o.g.).
    
    Now consider the walkers $\widetilde{X}^{\left(1\right)},\widetilde{X}^{\left(2\right)}$. One of
    them can have finite range only if the following holds. The same argument works for both walkers,
    we do it here for $\widetilde{X}^{\left(1\right)}$ w.l.o.g. $\widetilde{X}^{\left(1\right)}$
    can only have finite range if there exists $N$ such that for all $n \geq N$ we have $\omega_{m_n} = 1$.
    Assume to the contrary that no such $N$ exists. Then we can find arbitrarily large $n$ such that
    $\omega_{m_n} = 0$ which means that the labels of $X^{\left(1\right)}$ and $X^{\left(2\right)}$
    are \emph{not} exchanged in the interval $\left(\tau_{m_n},\tau_{m_n+1}\right]$. Since $X^{\left(1\right)}$ visits a node at distance at
    least $n$ from $0$ in this time interval, the same holds then for $\widetilde{X}^{\left(1\right)}$,
    so $\widetilde{X}^{\left(1\right)}$ would not have finite range.

    But the probability that the sequence $\omega_{m_n}$ is $1$ for all $n \geq N$ is $0$ for any $N$ (since the
    choice of $m_n$ only depends on the edge-reinforced random walk, i.e.~is independent of the $\omega_m$,
    and since the probability of $\omega$ being constantly $1$ on any fixed infinite subset of the
    integers is $0$ by the choice of $\omega$). Hence, the probability that such $N$ exists is $0$,
    and therefore the probability that $\widetilde{X}^{\left(1\right)}$ has finite range is $0$ as well,
    and the same arguments give that the probability for $\widetilde{X}^{\left(2\right)}$ having finite
    range is $0$ as well (both with respect to the measure $\bQ$).

    So the indicator variable in the integral above is $0$ almost surely w.r.t.~$\bQ$, and hence the
    inner integral is always $0$, which implies that the outer integral is also $0$ and hence $\Prb{A} = 0$.
    \item Similar to (i): Let $A$ now be the event that the integer $z$ is visited infinitely often by at least one of
    the walkers $X^{\left(1\right)},X^{\left(2\right)}$, that they meet infinitely often, and that one of
    them does \emph{not} visit $z$ infinitely often. Then, we have again:
    \begin{align*}
      \Prb{A} = &\int \int \mathbbm{1}_{B} \dx{\bQ} \dx{\bP} \\
      \textrm{where } B := &\left\{z\textrm{ visited }\infty\textrm{ often by at least one of }\widetilde{X}^{\left(1\right)},\widetilde{X}^{\left(2\right)}\textrm{,}\right.\\
      & \left.\hphantom{\;\{}\vphantom{\widetilde{X}^{\left(1\right)}}\textrm{they meet }\infty\textrm{ often, one of them visits }z\textrm{ only finitely often}\right\}\,.
    \end{align*}
    We see that the indicator variable can be $1$ only if at least one of $X^{\left(1\right)},X^{\left(2\right)}$
    visits $z$ infinitely often, and w.l.o.g.~assume that this holds for $X^{\left(1\right)}$. As before,
    we can construct a stricly increasing sequence $m_n$ such that in the time interval $\left(\tau_{m_n},\tau_{m_n+1}\right]$,
    $X^{\left(1\right)}$ visits $z$. Again as before, one of $\widetilde{X}^{\left(1\right)},\widetilde{X}^{\left(2\right)}$,
    take $\widetilde{X}^{\left(1\right)}$ w.l.o.g., can visit $z$ only finitely often only if
    $\omega_{m_n} = 1$ for all $n \geq N$ for some $N$, an event which has again probability $0$ w.r.t.~$\bQ$.
    \item Apply \ref{lem:labelExchangeBothVisInf} to every integer $z$. \qedhere
  \end{enumerate}
\end{oproof}

\subsection{Sequence-Type Reinforcement}

We now consider the following class of reinforcement schemes, called ``sequence type'' in \cite{errwz}.
Let $\mathbf{a} = \left(a_k\right)_{k \geq 1}$ be a sequence of nonnegative numbers. The reinforcement is
of sequence type if $W_z\left(k\right) = 1 + \sum_{l = 1}^{k} a_l$ (see \autoref{def:errw}) for all edges $z$.
In words, we require that all initial edge weights are $1$ (this is an additional property called ``initially fair''
in \cite{errwz}), and that the edge weights are then increased by a fixed sequence of increments which do not
depend on the location of the edge.

We want to show \autoref{thm:seqtypefinrange}, i.e.~if
\begin{align*}
  \phi\left(\mathbf{a}\right) = \sum_{k = 1}^{\infty} \left(1 + \sum_{l = 1}^k a_l\right)^{-1}
\end{align*}
then all walkers are recurrent a.s.~if $\phi\left(\mathbf{a}\right) = \infty$, and all walkers have finite range a.s.~otherwise.
This result is not surprising as the case $K=1$ is known due to \cite{errwz}.

\begin{oproof}[\autoref{thm:seqtypefinrange}]
  \begin{itemize}
    \item Case $\phi\left(\mathbf{a}\right) < \infty$: let $z_0 = \max\left\{ X_0^{\left(i\right)} : 1 \leq i \leq K \right\} + 1$. Consider $z \geq z_0$. Define a stopping time
    by $T_z := \inf\left\{ n \geq 0: X_n^{\left(i\right)} = z \textrm{ for some } i \right\}$. We want to bound
    $\Prb{T_{z + 2} < \infty \midd| T_z < \infty}$ uniformly away from $1$. Once this is proven, we can conclude
    that $\Prb{\forall z \geq z_0: T_z < \infty} = 0$ and therefore, by \autoref{thm:allRecurrentOrAllFiniteRange},
    all walkers must have finite range.

    Call $E_z$ the event that any walker reaching $z$ at some point in time only traverses
    the edge $\left\{z, z+1\right\}$ forever afterwards. For the uniform upper bound, first note that
    \begin{align*}
      \Prb{T_{z + 2} < \infty \midd| T_z < \infty} \leq 1 - \Prb{E_z}\,.
    \end{align*}
    Hence we want to a uniform lower bound for $\Prb{E_z}$. Consider the following situation. A given number of walkers are currently
    located at the two neighboring nodes $z$ and $z + 1$ which are incident to the edge $e = \left\{z, z+1\right\}$, which has been traversed $k$ times so far.
    Set $\alpha_k := 1 + \sum_{l = 1}^k a_l$ (this is the weight of edge $e$). $\phi\left(\mathbf{a}\right) < \infty$
    implies $\sum_{k \geq 0} \alpha_k^{-1} < \infty$. We do not assume anything about the weights
    of the adjacent edges at this point and call them $w_0$ and $w_1$.

    \begin{figure}[H]
      \begin{center}
        \begin{tikzpicture}
          \node[circle,fill=black,inner sep=0.7mm] (N0) at (0, 0) {};
          \node[circle,fill=black,inner sep=0.7mm] (N1) at (2, 0) {};
          \node[circle,fill=black,inner sep=0.7mm] (N2) at (4, 0) {};
          \node[circle,fill=black,inner sep=0.7mm] (N3) at (6, 0) {};
          \node[below=1mm] at (N1) {$z$};
          \node[below=1mm] at (N2) {~~~$z + 1$};
          \draw (N0) -- node[above=2mm] {$w_0$} (N1) -- node[above=2mm] {$\alpha_k$} node[below] {edge $e$} (N2) -- node[above=2mm] {$w_1$} (N3);
          \draw[dashed] (N3) -- (7, 0);
          \draw[dashed] (N0) -- (-1, 0);
          \stickman{(2,0.7)}{sectionblue}
          \stickman{(3.75,0.7)}{tumGreen}
          \stickman{(4.25,0.7)}{tumOrange}
        \end{tikzpicture}
        \caption{Some walkers are located at neighboring nodes}
        \label{fig:single_edge_trav}
      \end{center}
    \end{figure}
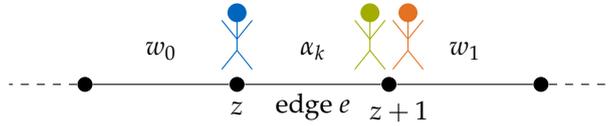
    We now give a lower bound for the probability that the edge $e$ will be traversed in the next step, if one of the walkers
    at nodes $z, z + 1$ is chosen to move:
    \begin{align*}
      \Prb{e\textrm{ traversed in next step}} &\geq \min\left\{\frac{\alpha_k}{w_0 + \alpha_k}, \frac{\alpha_k}{w_1 + \alpha_k}\right\}
      = \frac{\alpha_k}{w_{\textrm{m}} + \alpha_k} = \frac{1}{w_{\textrm{m}}\alpha_k^{-1} + 1} \\
      \textrm{where } w_{\textrm{m}} &= \max\left\{w_0, w_1\right\}\,.
    \end{align*}
    Using that $\frac{1}{x} \geq e^{1 - x}$, we further get
    \begin{align*}
      \Prb{e\textrm{ traversed in next step}} &\geq \exp\left(-w_{\textrm{m}}\alpha_k^{-1}\right) = \left(e_{\textrm{m}}\right)^{\alpha_k^{-1}} \\
      \textrm{where } e_{\textrm{m}} &= \exp\left(-w_{\textrm{m}}\right) > 0\,.
    \end{align*}

    We give a lower bound for $\Prb{E_z}$ as follows. Define $E_z^{\left(n\right)}$ as the event that up to time $n$,
    the behavior of the walkers is consistent with the event $E_z$, i.e.~up to time $n$, any walker which reached $z$ only traversed
    the edge $e = \left\{z, z + 1\right\}$ afterwards. Then $E_z = \bigcap_{n \geq 1} E_z^{\left(n\right)}$, and the situation looks
    as follows if $E_z^{\left(n\right)}$ occurs:

    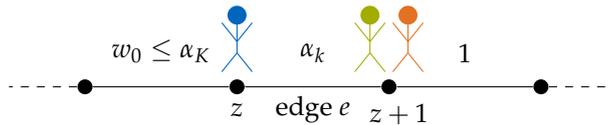
\begin{figure}[H]
      \begin{center}
        \begin{tikzpicture}
          \node[circle,fill=black,inner sep=0.7mm] (N0) at (0, 0) {};
          \node[circle,fill=black,inner sep=0.7mm] (N1) at (2, 0) {};
          \node[circle,fill=black,inner sep=0.7mm] (N2) at (4, 0) {};
          \node[circle,fill=black,inner sep=0.7mm] (N3) at (6, 0) {};
          \node[below=1mm] at (N1) {$z$};
          \node[below=1mm] at (N2) {~~~$z + 1$};
          \draw (N0) -- node[above=2mm] {$w_0 \leq \alpha_K$} (N1) -- node[above=2mm] {$\alpha_k$} node[below] {edge $e$} (N2) -- node[above=2mm] {$1$} (N3);
          \draw[dashed] (N3) -- (7, 0);
          \draw[dashed] (N0) -- (-1, 0);
          \stickman{(2,0.7)}{sectionblue}
          \stickman{(3.75,0.7)}{tumGreen}
          \stickman{(4.25,0.7)}{tumOrange}
        \end{tikzpicture}
        \caption{The walkers' behavior is still consistent with the event $E_z$}
        \label{fig:consistent_behavior}
      \end{center}
    \end{figure}
    Note:
    \begin{itemize}
      \item $z$ is to the right of the initial walker positions, so the weight of $\left\{z + 1, z + 2\right\}$
      must be $1$ as long as the walkers' behavior is consistent with $E_z$, since they arrive at $z$ from the left
      and are trapped in $e = \left\{z, z + 1\right\}$ if $E_z$ occurs.
      \item The weight $w_0$ of the edge $\left\{z - 1, z\right\}$ is at most $\alpha_K$ if the behavior is consistent:
      the edge can have been traversed at most $K$ times, since there are only $K$ walkers, and since a walker can
      never go back to $z - 1$ if he reaches $z$ and $E_z$ occurs.
      \item The weight of $e$ depends on the number of traversals of $e$, which can be arbitrary. We assume $k$
      traversals, so the edge weight is $\alpha_k$.
    \end{itemize}
    We have $w_{\textrm{m}} = \max\left\{w_0, w_1\right\} \leq \max\left\{\alpha_K, 1\right\} = \alpha_K$, so
    $e_{\textrm{m}} \geq \exp\left(-\alpha_K\right)$ in any such situation.
    \begin{align*}
      \Prb{E_z} = \prod_{n \geq 1} \Prb{E_z^{\left(n\right)}\midd|E_z^{\left(n - 1\right)}}
      \textrm{ since }E_z^{\left(n\right)} \subseteq E_z^{\left(n - 1\right)}
      \textrm{ and where }\Prb{E_z^{\left(0\right)}} = 1\,.
    \end{align*}
    But $\Prb{E_z^{\left(n\right)}\midd|E_z^{\left(n - 1\right)}}$ is the probability
    that $e$ is traversed if one of the walkers incident to $e$ is selected to move.
    Thus
    \begin{align*}
      \prod_{n \geq 1} \Prb{E_z^{\left(n\right)}\midd|E_z^{\left(n - 1\right)}}
      &\geq \prod_{k \geq 0} \exp\left(-\alpha_K\right)^{\alpha_k^{-1}}
      = \exp\left(-\alpha_K\right)^{\sum_{k \geq 0} \alpha_k^{-1}} > 0,
    \end{align*}
    since $\sum_{k \geq 0} \alpha_k^{-1} < \infty$. Note that the bound does not depend on $z$, so this is
    indeed the desired uniform lower bound.

    \item Case $\phi\left(\mathbf{a}\right) = \infty$: we set again $\alpha_k := 1 + \sum_{l = 1}^k a_l$, so we have
    $\sum_{k \geq 0} \alpha_k^{-1} = \infty$ in this case.
    We want to prove recurrence. By \autoref{thm:allRecurrentOrAllFiniteRange}, it suffices to show that there
    is at least one node which is visited infinitely often by at least one of the walkers. Assume for a contradiction that not every node is
    visited infinitely often and let $z$ be the largest node which is visited infinitely often. Consider the last time $n$ at which a node
    to the right of $z$ is occupied (visited) by any of the walkers. We will now look at the random walk from time $n + 1$ onwards, and we will show that the probability that
    $z + 1$ is never visited again is $0$, which is a contradiction. 
    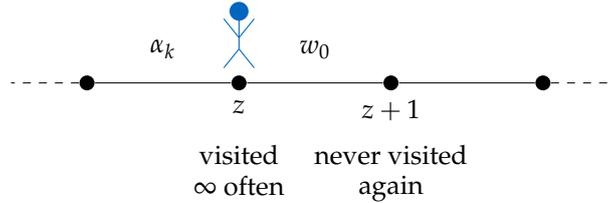
\begin{figure}[H]
      \begin{center}
        \begin{tikzpicture}
          \node[circle,fill=black,inner sep=0.7mm] (N0) at (0, 0) {};
          \node[circle,fill=black,inner sep=0.7mm] (N1) at (2, 0) {};
          \node[circle,fill=black,inner sep=0.7mm] (N2) at (4, 0) {};
          \node[circle,fill=black,inner sep=0.7mm] (N3) at (6, 0) {};
          \node[below=1mm] at (N1) {$z$};
          \node[below=7mm] at (N1) {visited};
          \node[below=11mm] at (N1) {$\infty$ often};
          \node[below=1mm] at (N2) {$z + 1$};
          \node[below=7mm] at (N2) {never visited};
          \node[below=11mm] at (N2) {again};
          \draw (N0) -- node[above=2mm] {$\alpha_k$} (N1) -- node[above=2mm] {$w_0$} (N2) -- (N3);
          \draw[dashed] (N3) -- (7, 0);
          \draw[dashed] (N0) -- (-1, 0);
          \stickman{(2,0.7)}{sectionblue}
        \end{tikzpicture}
        \caption{$z$ is visited infinitely often, but not its right neighbor}
        \label{fig:neighbor_vis_inf}
      \end{center}
    \end{figure}
    We assume that $\left\{z - 1, z\right\}$ has been traversed $k$ times and therefore has weight $\alpha_k$, and we call
    $w_0$ the arbitrary weight of $\left\{z, z + 1\right\}$ at time $n + 1$.
    Let $\tau_1, \tau_2, \tau_3, \ldots$ be the times at which a walker located at $z$ is selected to move after time $n$.
    There will be infinitely many such times since there are infinitely many visits to $z$. We call $t_m$ the number of
    traversals of $\left\{z - 1, z\right\}$ at time $\tau_m$. Note that
    $t_m \leq k + K - 1 + 2\left(m - 1\right)$. The reason for this is as follows: at time $n$, there was still one walker to the
    right of $z$ by definition of the time $n$, and he must have moved to $z$ at time $n + 1$. Hence, there is at least one
    walker at $z$ at time $n + 1$. Before a walker at $z$ is selected to move, the weight of $\left\{z - 1, z\right\}$ could still
    be increased by walkers moving to $z$ from the left, but since there are only $K$ walkers, the number of traversals can increase
    by at most $K - 1$. At subsequent times $\tau_m$, the only other possibility for the weight of $\left\{z - 1, z\right\}$ to
    increase is that a walker leaves $z$ and then (possibly) returns, adding an additional two traversals. This leads to the given upper bound.
    Then, at each time $\tau_m$, we have (recall that at time $\tau_m$, a walker at $z$ was selected to move)
    \begin{align*}
      \Prb{\textrm{the selected walker moves to }z - 1} &= \frac{\alpha_{t_m}}{\alpha_{t_m} + w_0}\\
      &\leq \frac{\alpha_{k + K + 2m - 3}}{\alpha_{k + K + 2m - 3} + w_0} \\
      &\hskip-0.25cm \implies \\
      \Prb{\textrm{at all times }\tau_m\textrm{, the selected walker moves to }z - 1} &\leq \prod_{m \geq 1} \frac{\alpha_{\left(k + K - 3\right) + 2m}}{\alpha_{\left(k + K - 3\right) + 2m} + w_0} \\
      &\hskip-0.4cm \underset{l = k + K - 3}{\overset{\circledast}{\leq}} \exp\left(- \sum_{m \geq 1} \frac{w_0}{\alpha_{l+2m} + w_0}\right) \\
      \textrm{where }\circledast\textrm{ holds since } \forall x \in \bbR: x &\leq \exp\left(x - 1\right)
    \end{align*}
    Since $\alpha_k \geq 1$ for all $k$, we have that $\alpha_{l+2m} + w_0 \leq \left(1 + w_0\right)\alpha_{l+2m}$, and we can analyze the sum as follows:
    \begin{align*}
      \sum_{m \geq 1} \frac{w_0}{\alpha_{l+2m} + w_0} &\geq \frac{w_0}{1 + w_0} \cdot \sum_{m \geq 1} \frac{1}{\alpha_{l+2m}}
      \overset{\alpha_m \textrm{ increasing in }m}{\geq} \frac{w_0}{1 + w_0} \cdot \frac{1}{2} \cdot \sum_{m \geq 2} \frac{1}{\alpha_{l+m}} \\
      &= \frac{w_0}{2\left(1 + w_0\right)} \cdot \sum_{m \geq l + 2} \alpha_m^{-1} = \infty\,.
    \end{align*}
    Therefore,
    \begin{align*}
      \Prb{\textrm{at all times }\tau_m\textrm{, the selected walker moves to }z - 1}
      &\leq \exp\left(- \sum_{m \geq 1} \frac{w_0}{\alpha_{l+2m} + w_0}\right) = 0\,.
    \end{align*}
    Hence, the probability that $z + 1$ is never visited again after time $n$ is $0$, which concludes the proof. \qedhere
  \end{itemize}
\end{oproof}

Nina Gantert: Technical University of Munich, School of Computation, Information and Technology, Department of Mathematics, Boltzmannstr. 3, 85748 Garching, Germany. \href{mailto:gantert@ma.tum.de}{gantert@ma.tum.de}

Fabian Michel: Technical University of Munich, School of Computation, Information and Technology, Department of Mathematics, Boltzmannstr. 3, 85748 Garching, Germany. E-mail: \href{mailto:fabian.michel@tum.de}{fabian.michel@tum.de}

Guilherme H.~de Paula Reis: IMPA -- Instituto Nacional de Matem\'atica Pura e Aplicada, Estrada Dona Castorina 110, 22460-320 Rio de Janeiro, Brazil. \href{mailto:guilherme.reis.mat@gmail.com}{guilherme.reis.mat@gmail.com}

\end{document}